\documentclass[nohyperref]{article}

\usepackage{microtype}
\usepackage{graphicx}
\usepackage{booktabs} 

\usepackage{hyperref}

 \usepackage[accepted]{icml2022}

\usepackage{amsmath}
\usepackage{amssymb}
\usepackage{mathtools}
\usepackage{amsthm}

\usepackage[capitalize,noabbrev]{cleveref}

\usepackage{url}
\usepackage{algorithm}
\usepackage{algorithmic}
\usepackage{comment}
\usepackage{wrapfig,lipsum}
\usepackage{enumitem}
\usepackage{multicol}
\usepackage{caption}
\usepackage{subcaption}
\usepackage{amsfonts,bm}
\usepackage{multirow}

\theoremstyle{plain}

\theoremstyle{definition}

\theoremstyle{remark}

\usepackage[textsize=tiny]{todonotes}

\icmltitlerunning{Constrained Discrete Black-Box Optimization using MILP}

\begin{document}

\twocolumn[
\icmltitle{Constrained Discrete Black-Box Optimization using Mixed-Integer Programming}



\icmlsetsymbol{equal}{*}

\begin{icmlauthorlist}
\icmlauthor{Theodore P. Papalexopoulos}{mit,google}
\icmlauthor{Christian Tjandraatmadja}{google}
\icmlauthor{Ross Anderson}{google}
\icmlauthor{Juan Pablo Vielma}{google}
\icmlauthor{David Belanger}{google}
\end{icmlauthorlist}

\icmlaffiliation{mit}{Operations Research Center, Massachusetts Institute of Technology, Cambridge MA, USA}
\icmlaffiliation{google}{Google Research, Cambridge MA, USA}

\icmlcorrespondingauthor{Theodore P. Papalexopoulos}{tedpapalex@gmail.com}

\icmlkeywords{Black-Box Optimization, Neural Networks, Mixed-Integer Programming}

\vskip 0.3in
]



\printAffiliationsAndNotice{}  

\begin{abstract}
Discrete black-box optimization problems are challenging for model-based optimization (MBO) algorithms, such as Bayesian optimization, due to the size of the search space and the need to satisfy combinatorial constraints.  In particular, these methods require repeatedly solving a complex discrete global optimization problem in the inner loop, where popular heuristic inner-loop solvers introduce approximations and are difficult to adapt to combinatorial constraints. In response, we propose NN+MILP, a general discrete MBO framework using piecewise-linear neural networks as surrogate models and mixed-integer linear programming (MILP) to optimize the acquisition function. MILP provides optimality guarantees and a versatile declarative language for domain-specific constraints. We test our approach on a range of unconstrained and constrained problems, including DNA binding, constrained binary quadratic problems from the MINLPLib benchmark, and the NAS-Bench-101 neural architecture search benchmark. NN+MILP surpasses or matches the performance of black-box algorithms tailored to the constraints at hand, with global optimization of the acquisition problem running in a few minutes using only standard software packages and hardware.
\end{abstract}

\section{Introduction}
\label{intro}

The problem of optimizing an expensive black-box function $f: \Omega \mapsto \mathbb{R}$ over a discrete, constrained domain arises in numerous application domains, e.g.\ neural architecture search~\citep{zoph2016neural}, program synthesis~\citep{summers1977methodology, biermann1978inference}, small-molecule design~\citep{elton2019deep}, and protein design~\citep{yang2019machine}. In such resource-constrained settings, it is desirable to develop algorithms that exploit known combinatorial structure in $\Omega$ to search the space more efficiently.

Model-based Black-box Optimization (MBO), a popular paradigm that includes Bayesian Optimization as a special case, iteratively refines a function approximator $\hat{f}(x) \approx f(x)$ and selects new points to query by optimizing an \textit{acquisition function} $a(x)$ derived from a point estimate or posterior distribution over $\hat{f}$ (Section~\ref{sec:background_MBO}). This \textit{inner-loop optimization} problem is assumed to be easier than the original, since, for example, $a(x)$ is less expensive to query than $f(x)$ or provides ``white-box'' properties such as gradients.

There is a vast literature addressing the challenges of applying MBO in practice. We focus on two of these: first, optimizing $a(x)$ may itself be a computationally-difficult optimization problem; second, in many applications, practitioners are confronted by additional constraints on $x$. For example, in neural architecture search, $x$ might represent a computation graph that must be both connected and acyclic. Due to the difficulty in optimizing the acquisition function over a combinatorial domain, most approaches resort to heuristic inner-loop solvers, which often need to be specialized to the problem at hand to ensure feasibility, e.g., evolutionary solvers with custom mutation operators.

To address the challenge of inner-loop optimization, we introduce a general framework for discrete, constrained MBO, \textit{NN+MILP}, that \textit{exactly} solves the acquisition problem using mixed-integer linear programming (MILP). Crucially, by framing the inner-loop optimization as an MILP, our approach can flexibly incorporate a wide variety of logical, combinatorial, and polyhedral constraints on the domain, which need only be provided in a \textit{declarative} sense. 

Using MILP in the inner loop does restrict the functional form of $\hat{f}$ (or the acquisition function based on it), but it supports any piecewise linear function. In particular, we employ the class of neural network (NN) approximators with ReLU activation functions due to their scalability and accuracy in practice, and because we can draw on recent work improving the performance of MILP for optimizing such NNs with respect to their inputs~\citep{anderson2020strong}. For us, MILP is practical to use in the inner loop because the dimensionality of typical black-box optimization problems is orders of magnitude smaller than those usually considered by MILP solvers. Our contributions are as follows:

\begin{itemize}[leftmargin=*]
    \item We introduce \textit{NN+MILP}, an MBO framework for discrete black-box problems with NN surrogates and exact optimality guarantees for solving the acquisition problem.
    \item We show that \textit{NN+MILP} matches or surpasses the performance of strong MBO baselines based on problem-specific evolutionary algorithms on a wide range of synthetic and real-world discrete black-box problems.
    \item We observe in our experiments that the runtime of MILP is practical for use with black-box problems of real-world scale, often solving the inner acquisition problem in seconds using standard packages and hardware. 
    \item We test our algorithm on a range of constrained binary quadratic problems from the MINLPLib benchmark, to highlight MILP's flexible declarative language for problem-specific constraints. 
    \item We use the NAS-Bench-101 neural architecture search benchmark as a case study, presenting a novel MILP formulation of its graph-structured domain. 
\end{itemize}

\section{Background and Related Work}
\label{background}

\subsection{Model-Based Black-Box Optimization}
\label{sec:background_MBO}

Model-based Black-box Optimization (MBO) is a broad family of methods that includes Bayesian optimization as a special case~\citep{mockus1978application,Jones1998efficient,hutter2011sequential,snoek2012practical,Shahriari2015taking}. As depicted in Algorithm~\ref{alg:generic_MBO}, the method proposes $x_t$ at iteration $t$ using three steps. First, the user performs inference over a \textit{surrogate model} $\hat{f}$ to approximate $f$ using the data previously collected from the black-box function. Here, $\texttt{fit}()$ may return a point estimate for $\hat{f}$, a posterior distribution over $\hat{f}$, or a posterior predictive distribution. Next, an \textit{acquisition function} $a(x)$ based on $\hat{f}(x)$ is posed that quantifies the quality of new points to query. Finally, $x_t$ is selected as the best point found by solving the \textit{acquisition problem}, where an \textit{inner-loop solver} (approximately) optimizes $a(x)$. The acquisition problem is typically designed such that it is more approachable than directly solving the original problem. For example, $a(x)$ may be orders of magnitude less expensive to evaluate or have a tractable functional form. Practitioners can encode prior knowledge about the structure of $f$ via a choice of inductive bias for $\hat{f}$, e.g., a suitable Gaussian Process kernel or neural-network architecture.

Bayesian optimization performs Bayesian inference over $\hat{f}$ and employs an acquisition function that accounts for uncertainty in $\hat{f}$. Doing so provides principled mechanisms for balancing exploration and exploitation~\citep{mockus1978application, srinivas2010gaussian} and is particularly important in early rounds of optimization when models are fit on limited data. We refer to our method as an instance of MBO, not Bayesian optimization, because it does not assume formal Bayesian inference for $\hat{f}$. Gaussian processes (GPs) are often used for $\hat{f}$ in Bayesian optimization, since they provide closed-form posterior inference, naturally adjust their expressivity as the dataset grows, and users can inject domain knowledge via a choice of kernel~\citep{rasmussen2006gp,oh2019combinatorial}. On the other hand, neural networks provide a practical alternative~\citep{snoek2015scalable,hernandez2017parallel}, since they often scale more gracefully, either computationally or statistically, to large datasets or high-dimensional domains.

\begin{figure}[!t]
\begin{minipage}[t]{0.46\textwidth}
        \begin{algorithm}[H]
           \caption{MBO}
           \label{alg:generic_MBO}
        \begin{algorithmic}
           \STATE {\bfseries Input:} hypothesis class $\mathcal{F}$, budget $N$, initial dataset $\mathcal{D}_{n} = \{x_i, f(x_i)\}_{i=1}^n$, optimization domain $\Omega$
           \FOR{$t = n + 1$ to $t = N$}
             \STATE $P(\hat{f}_t) \leftarrow \texttt{fit}(\mathcal{F}, \mathcal{D}_{t-1})$ \;
             \STATE $a(x) \leftarrow \texttt{get\_acquisition\_function}(P(\hat{f}_t))$ 
             \STATE $x_t \leftarrow \texttt{inner\_loop\_solver}(a(x), \Omega)$
             \STATE $\mathcal{D}_t \leftarrow \mathcal{D}_{t-1} \cup \{x_t, f(x_t)\}$
           \ENDFOR
           \STATE \textbf{return} $\arg \max_{(x_t, y_t) \in \mathcal{D}_N} y_t$
        \end{algorithmic}
        \end{algorithm}
\end{minipage}
\hfill
\begin{minipage}[t]{0.46\textwidth}
    \begin{algorithm}[H]
       \caption{NN+MILP}
       \label{alg:nn_milp}
    \begin{algorithmic}
       \STATE {\bfseries Input:} hypothesis class $\mathcal{F}$, budget $N$, initial dataset $\mathcal{D}_{n} = \{x_i, f(x_i)\}_{i=1}^n$, MILP domain formulation $\mathcal{M}_\Omega$
       \FOR{$t = n + 1$ to $t = N$}
         \STATE $\hat{f}_t \leftarrow \texttt{fit}(\mathcal{F}, \mathcal{D}_{t-1})$ \; \hfill (\ref{sec:algo_approx})
         \STATE $\mathcal{M}_t \leftarrow \texttt{build\_milp}(\hat{f}_t,\; \mathcal{M}_\Omega, \; \mathcal{D}_{t-1})$ \hfill (\ref{sec:algo_opt})\
         \STATE $x_t \leftarrow \texttt{optimize}(\mathcal{M}_t)$ \; \hfill (generic MILP solver)
         \STATE $\mathcal{D}_t \leftarrow \mathcal{D}_{t-1} \cup \{x_t, f(x_t)\}$
       \ENDFOR
       \STATE \textbf{return} $\arg \max_{(x_t, y_t) \in \mathcal{D}_N} y_t$
    \end{algorithmic}
    \end{algorithm} 
\end{minipage}
\end{figure}

\subsection{Solving the Discrete MBO Acquisition Problem}
\label{sec:inner_problem_rw}
In general, the inner-loop problem is itself a non-trivial global optimization problem. Prior work on discrete MBO has mainly employed local search solvers, such as evolutionary search, with limited guarantees \citep{hutter2011sequential,muller2016miso,oh2019combinatorial,kandasamy2020tuning}. A key advantage of such solvers is that they treat $a(x)$ as a black box, which provides practitioners with freedom when designing application-specific surrogate models. On the other hand, certain choices of surrogate model and acquisition function lead to acquisition problems that can be (approximately) solved using specialized combinatorial solvers~\citep{baptista2018bayesian,deshwal2020mercer}, mixed-integer nonlinear programming (MINLP) ~\citep{costa2018rbfopt, kim2020surrogate}, or continuous optimization solvers~\citep{bliek2021black}. 

Therefore, practitioners must decide between either introducing difficult-to-analyze approximations due to inexact heuristic solvers or using tractable surrogate models that may be mis-specified for the application domain. This serves as a key motivation for our work: we seek to enable practitioners to employ broad families of surrogate models and exactly solve the acquisition problem with reasonable computational overhead in practice. 

\subsection{Constrained MBO}
\label{sec:background_constrained_MBO}
In many applications, $x$ is subject to non-trivial structural constraints.
Prior work has largely focused on the case where determining whether $x$ is feasible requires evaluating an expensive, perhaps noisy, black-box function $h(x)$ with cost comparable to $f(x)$~\citep{schonlau1998global,gelbart2014bayesian,hernandez2016general,ariafar2019admmbo,letham2019constrained}. Here, standard acquisition functions can be extended to account for an additional classifier $\hat{h}(x)$ trained to predict $h(x)$. 

Problems with inexpensive white-box $h(x)$ can be tackled using these approaches for black-box constraints, but doing so may lead to slower optimization and may query $f(x)$ at invalid $x$, which can be unsafe when performing physical experiments~\citep{berkenkamp2016safe}. Instead, the inner-loop solver can be modified directly to guarantee feasibility, e.g., by using rejection sampling~\citep{shi2020learned,kandasamy2020tuning}. If using local search algorithms, the solver would need to be customized for each family of constraints, a task usually left to the user.  Prior work employing MINLP solvers addresses white-box constraints either by adding a penalty for constraint violation~\citep{costa2018rbfopt} or in small-scale settings~\citep{kim2020surrogate}. Conversely, our approach unifies both the surrogate model and domain within the same declarative constraint framework (MILP), and thus allows for exact optimization over general combinatorial domains with minimal algorithmic effort on the part of the user.

\subsection{Mixed Integer Linear Programming}
\label{sec:background_mip}

Mixed Integer Linear Programming (MILP) seeks to maximize a linear function over a set of decision variables, some of which may be integral, subject to linear inequality constraints. Decades of development have allowed MILP to have a significant impact in a wide range of  applications due to its better-than-expected computational performance~\citep{Juenger2010}. Indeed, while MILP problems are computationally hard (NP-complete), they are routinely solved (to global or near-global optimality) in production environments thanks to state-of-the-art solvers that nearly double their machine-independent performance every year \citep{Achterberg2013,bixby2012brief}. 

A notable aspect of MILP is that it provides a simple yet extremely versatile declarative language for white-box constraints.  It is well known that linear inequalities over integer variables can be used to easily build \emph{pure-integer} formulations for logical constraints and combinatorial optimization problems 
\citep{williams2013model, schrijver2003combinatorial, wolsey1999integer}. In addition, using both integer and continuous variables leads to \emph{mixed-integer} formulations that can combine polyhedral and logical constraints \citep{jeroslow1989logic, pochet2006production, Mixed-Integer-Linear-Programming-Formulation-Techniques}. 

Particularly interesting to our proposed approach are MILP formulations for piecewise-linear functions \citep{huchette2017nonconvex, Mixed-Integer-Models-for-Nonseparable}. Specifically, our work leverages MILP formulations for trained neural networks with piecewise-linear activation functions such as ReLUs~\citep{anderson2020strong}. Optimizing over trained ReLU networks with MILP has been done in contexts such as neural network verification~\citep{cheng2017maximum, lomuscio2017approach, tjeng2017evaluating}, reinforcement learning~\citep{ryu2019caql, delarue2020reinforcement}, and analysis and exact compression of neural networks~\citep{serra2018bounding,serra2021scaling}. MILP has also been used to optimize ReLU network surrogates of simulation-based constraints~\citep{grimstad2019relu}, although their approach optimizes a single surrogate model once, unlike in ours.

\section{MILP for MBO}
\label{sec:algo}

We propose the \textit{NN+MILP} framework (Algorithm~\ref{alg:nn_milp}), which uses neural network surrogate models and solves the acquisition problem using MILP at every step. This provides practitioners with the flexibility to use a wide variety of models and leverage MILP's versatile declarative language to incorporate constraints. This section describes various design choices to make the approach practical.

\subsection{Problem Setting}
\label{sec:problem_setting}

Our goal is to find:
\begin{equation}
    \label{eq:bb_problem}
    x^* = \arg\max_{x \in \Omega} f(x),
\end{equation}
where $f: \Omega \mapsto \mathbb{R}$ is an expensive, noiseless black-box function and $\Omega \subseteq \Omega_1 \times \ldots \times \Omega_n$ is a domain on $n$ decision variables. We assume $\Omega$ can be described by an inexpensive function  $h_\Omega(x)$ indicating whether $x$ is in $\Omega$. Algorithms are allowed a fixed budget of $N$ sequential queries to $f$. $\mathcal{X}_t := \{x_i\}_{i=1}^{t}$ refers to the set of sampled points by iteration $t$, and $\mathcal{D}_t := \{x_i, y_i = f(x_i)\}_{i=1}^{t}$ includes corresponding rewards. We measure performance by the best reward in $\mathcal{D}_N$. Since $f$ is noiseless, it is advantageous for algorithms to avoid repeated evaluations of the same $x$.

We choose to focus on finite discrete sets $\Omega$ as we believe this is the area where MILP can provide the greatest benefit. As noted in Section~\ref{sec:background_mip}, there are many well-studied formulation techniques for $\Omega$ with combinatorial structure, such as directed graphs. More generally, such sets have a polynomially-sized MILP formulation whenever $h_\Omega(x)$ can be evaluated in polynomial time (e.g., \citet{yannakakis1991expressing}). Continuous and mixed-integer domains could be incorporated in our approach with some modifications (Section~\ref{sec:conclusion}), although they are outside the scope of this paper.

\subsection{Surrogate Model and Acquisition Function}
\label{sec:algo_approx}

For surrogate model $\hat{f}$, we allow any feedforward neural network with piecewise-linear activation functions, as they can be represented by MILP (Section~\ref{sec:algo_opt}). Though we focus on fully-connected ReLU networks, a range of such architectures (e.g., with convolutional or max-pooling layers) can be used to place problem-specific inductive bias on $\hat{f}$. 

In order to manage the tradeoff between exploration and exploitation, we employ a heuristic based on the well-established Thompson sampling approach~\citep{thompson1933likelihood,hernandez2017parallel,kandasamy2018parallelised}. In step $t$ of Thompson sampling, a model $\hat{f}(x)$ is sampled from the posterior $P(\hat{f} | \mathcal{D}_{t-1})$, and a greedy action is taken with respect to the model, i.e., $a(x) = \hat{f}(x)$. We approximate this by using an informal method to generate posterior samples that has been shown in prior work to perform well~\citep{lakshminarayanansimple, riquelme2018deep}: we train $\hat{f}(x)$ from scratch at each iteration using random parameter initialization and stochastic gradient descent. Our method is orthogonal to the choice of posterior sampling technique, though, and variational methods or MCMC could be used in the future. We also discuss alternative acquisition functions in Section~\ref{sec:conclusion}.

We select the capacity of the surrogate -- i.e., the number of layers and neurons in the network -- so as to balance expressivity and statistical/computational scalability. Given the relatively small number of dimensions and training points, particularly in early iterations, larger networks are likely to overfit, while also being more computationally expensive to optimize. We empirically find that small, single-layer networks often suffice in our setting, with larger networks not improving results significantly (see Section~\ref{sec:exp_constr}). While out of scope for this paper, we also note that, in general, the size of the surrogate could be gradually increased across iterations to reflect the larger number of training points.

We use a flattened one-hot encoding of $x$ for the input layer, and train each network $\hat{f}_t \in \mathcal{F}$ on $\mathcal{D}_{t-1}$ using $\ell_2$ loss. Before training, we re-scale the observed rewards in $\mathcal{D}_{t-1}$ to aid both in training and optimization. Poorly-scaled data may result in slower performance or small inaccuracies in MILP solvers~\citep{miltenberger2018exploring}.

\subsection{MILP Formulation of the Acquisition Problem}
\label{sec:algo_opt}

The inner-loop solver then seeks to find 
\begin{equation}
\label{eq:inner_problem}
x_t = \arg\max_{x \in \Omega \setminus \mathcal{X}_{t-1}} \hat{f}_t(x),
\end{equation}
where $\Omega$ is the feasible set for~\eqref{eq:bb_problem} and $\mathcal{X}_{t-1}$ is the set of points where the noiseless $f(x)$ has been queried already. The MILP formulation of \eqref{eq:inner_problem} is denoted by $\mathcal{M}_t$ and has the following three components:

\textbf{Domain}\; We use a one-hot encoding of decision variables $x$ (unless they are already binary), defining the binary decision vector $z$ with $z_{ij} \equiv \mathbb{I}\{x_i = j\}$ for $i \in [n], j \in \Omega_i$, and subject to linear constraints $\sum_{j \in \Omega_i} z_{ij} = 1 \;\forall i$. Integer domains with small range may be one-hot encoded; see Appendix~\ref{sec:appendix_binary_vs_integer} for a comparison between integer and one-hot encodings.

Additional constraints due to $\Omega$ are added as necessary, with form dependent on the application at hand. We assume that these are MILP-representable, which as noted in Section~\ref{sec:background_mip} could include a wide range of combinatorial, logical, and polyhedral constraints. We use $\mathcal{M}_\Omega$ to denote the domain formulation itself.

\textbf{No-good Constraints}\; A \textit{no-good constraint} is one that eliminates undesirable solutions from the domain. Here, we leverage the binary nature of $z$ to \textit{exactly} eliminate the set $\mathcal{X}_{t-1}$ from $\mathcal{M}_t$. For illustrative purposes, consider a single point $\bar{x} \in \Omega$ we wish to exclude from the acquisition problem's domain, and let $\bar{z}$ denote its one-hot encoding (or $\bar{x}$ itself if the problem is binary).

Then the constraint: 
\begin{equation}
    \sum_{i,j\ :\ \bar{z}_{ij} = 0} z_{ij} + \sum_{i,j\ :\ \bar{z}_{ij} = 1} (1 - z_{ij}) \geq 1
\end{equation}
enforces that any feasible $z$ has a Hamming distance of at least 1 from $\bar{z}$. As $z$ are binary, this effectively eliminates just the single point $\bar{z}$ from the feasible region. We therefore formulate ${\Omega \setminus \mathcal{X}_{t-1}}$ by including one such constraint for each $\bar{x} \in \mathcal{X}_{t-1}$. Note that the right-hand side can be tightened to 2 for one-hot encodings, and these no-good constraints do not extend naturally to continuous $x$ (Section~\ref{sec:conclusion}).

\textbf{Neural Network}\; We formulate the neural network by introducing auxiliary decision variables encoding the activation of each neuron for a given $z$. We present here the formulation for a single ReLU, commonly used throughout the literature (Section~\ref{sec:background_mip}), while noting that the full formulation is obtained by combining all ReLU formulations and matching their input and output variables according to the structure of the network. The overall MILP objective is the activation corresponding to the regressor's output neuron.

A ReLU neuron with vector input $x$ and scalar output $y$ has the piecewise-linear form $y = \max(0, w^{\top} x + b)$, where $w$ and $b$ are its weights and bias respectively. At optimization time, $w$ and $b$ are fixed, while $x$ and $y$ are represented by decision variables (also used as the inputs and outputs of other ReLUs according to the feedforward structure). To handle the non-linearity, we add a binary decision variable $\alpha$ that indicates whether the ReLU is active or not. We then write the following set of constraints to enforce that $y = 0$ when $\alpha = 0$ and $y = w^{\top} x + b$ when $\alpha = 1$:
\begin{align}
    0 \leq y \leq  &\ M\alpha \label{eq:relu_0}\\
    w^\top x + b \leq y \leq &\ w^\top x + b + M(1 - \alpha) \label{eq:relu_1}
\end{align}
where $M$ is a sufficiently large fixed value, such as an upper bound on the range of $y$. As $w$ and $b$ are fixed, values for $M$ can be computed in advance of the optimization, e.g., by propagating bounds from $\Omega$. Our experiments use a more advanced method to compute $M$, detailed in Appendix~\ref{sec:relu_formulation}.

\begin{figure*}[ht]
\begin{center}
\centerline{\includegraphics[width=\textwidth]{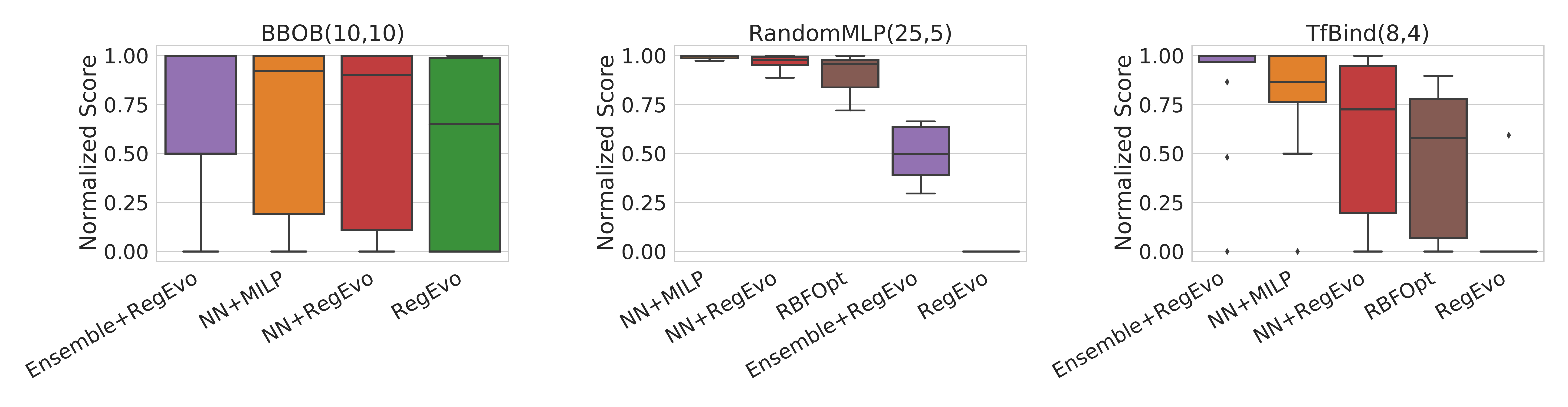}}
\vskip -0.1in
\centerline{\includegraphics[width=\textwidth]{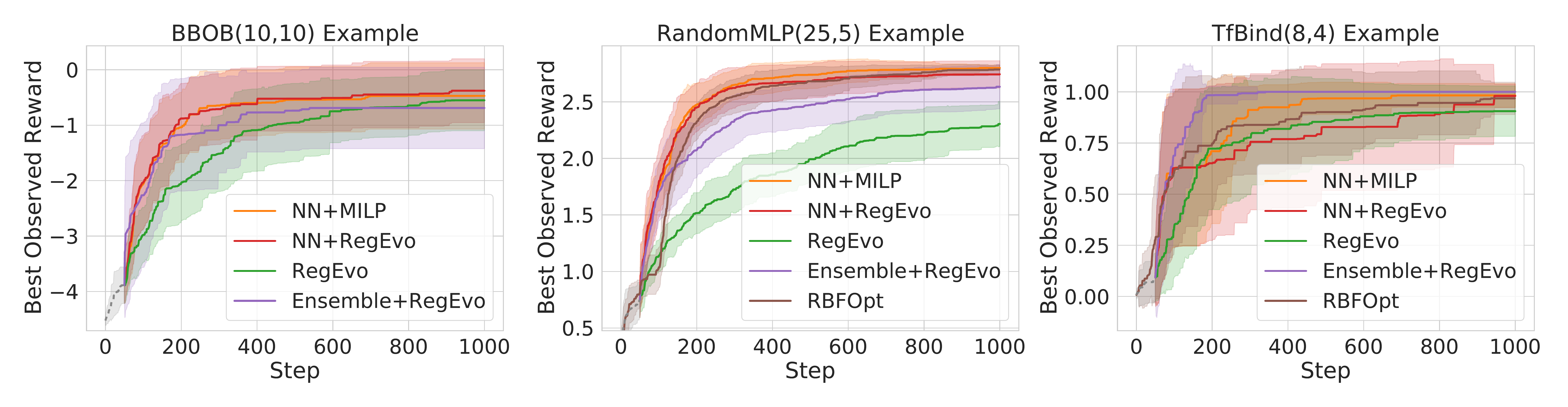}}
\vskip -0.1in
\caption{(Top) Distribution of algorithms' normalized scores (Section~\ref{sec:exp_tasks}) on unconstrained problems split by class. Higher is better. NN+MILP matches or outperforms NN+RegEvo on 22/30 problems. See Appendix~\ref{sec:appendix_results_unconstr} for an alternate plot where scores correspond to area under the best-observed reward curve (AUC). (Bottom) Best observed reward as a function of iteration for an example problem in each class, averaged over 20 trials (bands indicate $\pm 1$sd). Dashed grey lines in the first 50 steps indicate the initial randomly sampled dataset, common to all methods except RBFOpt.}
\label{fig:unconstrained_scores}
\end{center}
\end{figure*}

\subsection{Optimality Guarantees for MILP}
\label{sec:algo_guarantees}
The full acquisition problem formulation, denoted by $\mathcal{M}_t$, is passed to a generic MILP solver with fixed time limit. If the solver does not time out, it is guaranteed to have produced a global optimum of \eqref{eq:inner_problem}. Even if the solver times out, it will return the best feasible solution it found, plus an upper bound on the global optimal value. This bound can be used to evaluate the level of \emph{potential} sub-optimality of the feasible solution. Note that solvers often find an optimal solution before finding the upper bound that guarantees its optimality, so timing out do not imply sub-optimality. Finally, inner-loop optimality guarantees do not translate into guarantees for the overall black-box optimization, particularly when $f(x)$ does not belong to $\mathcal{F}$. However, they do provide a useful empirical tool for understanding the impact of exact inner-loop optimization (Section~\ref{sec:experiments}).

\section{Experiments}
\label{sec:experiments}

This section presents experimental results on a wide range of discrete black-box problems, with and without combinatorial constraints. We focus primarily on analyzing the effect of \textit{global} optimization of the acquisition function, by including controlled ablations of \textit{NN+MILP} where the inner-loop solver is replaced by an inexact evolutionary alternative. Depending on the problem, we also include independent baselines tailored to the application domain. 

In all experiments, we fix the surrogate model hypothesis class $\mathcal{F}$ to networks with a single, fully-connected hidden layer of 16 neurons. Models are trained with TensorFlow~\citep{abadi2016tensorflow}, using the ADAM optimizer. No hyper-parameter tuning is performed across problems. The MILP acquisition problem is solved with the Mixed-Integer Programming solver SCIP 7.0.1 \citep{gamrath2020scip} using default settings. While the acquisition problem is typically solved to optimality in seconds (Section~\ref{sec:exp_pract}), we set a time limit of 500s as a safeguard.  We use standard CPU machines with $\sim$1G RAM and $\leq 10$ cores.

\subsection{Benchmarking Tasks}
\label{sec:exp_tasks}

Unless otherwise stated, tasks' domains consist of discrete vectors of length $n$, with a common alphabet $\mathcal{A}$ for all elements. We consider four families of black-box objectives: 

\begin{itemize}[leftmargin=*]
    \item \textbf{RandomMLP} The output of a multi-layer perceptron operating on a one-hot encoding of the input. Notably, architectures have significantly more layers/parameters than the 16-neuron networks used as surrogates by \textit{NN+MILP}.
    \item \textbf{TfBind} Binding strength of a length-8 DNA sequence to a given transcription factor \citep{barrera2016survey}.
    \item \textbf{BBOB} Non-linear function from the continuous Black-Box Optimization Benchmarking library \citep{hansen2009real}, where each coordinate is uniformly discretized along its range. Despite the underlying continuous structure, inputs are treated as unordered and categorical.
    \item \textbf{Ising} The negative energy of fully-connected binary Ising Model with normally distributed pairwise potentials.
\end{itemize}
 
We use parentheses after the family name to denote dimensionality of a problem, e.g., RandomMLP(10,5) refers to a RandomMLP objective over a discrete domain with $n = 10$ and $|\mathcal{A}| = 5$. Appendix~\ref{sec:appendix_tasks} lists all functions considered, and provides further details on the BBOB discretization.

Algorithms are evaluated in terms of the best reward observed after 1000 queries, averaged over 20 trials per problem. Algorithms' performance is significantly influenced by the set of $x$ that are proposed in early iterations. Therefore, to reduce variance when comparing algorithms, we initialize each of the 20 trials with a different fixed dataset of 50 random points. To facilitate comparison across problems with different reward scales, the algorithms' average final rewards are min/max normalized within each problem. That is, the best (resp.\ worst) on-average algorithm for a given problem is assigned a score of one (resp.\ zero), and intermediate values express relative distance from these extremes. No hyper-parameter tuning was performed across problems.

\subsection{Unconstrained Optimization}
\label{sec:exp_unconstr}

Before considering problems with combinatorial white-box constraints, we first tackle simple problems with no additional constraints on the discrete domain, i.e.,\ $\Omega = \mathcal{A}^n$. This allows us to compare against general-purpose algorithms for unconstrained discrete black-box optimization. We vary the problem sizes over 30 functions, consisting of eight RandomMLP(25,5), ten BBOB(10,10) and twelve TfBind(8,4) targets (Appendix~\ref{sec:appendix_tasks}).

NN-MILP provides an analytical tool for understanding the relative impacts of the choice of surrogate model and whether the acquisition problem is solved to optimality. Doing so requires ablations that vary along two axes: the family of surrogate models and the inner-loop solver. Further configuration details are provided in Appendix~\ref{sec:appendix_solvers}.

\begin{itemize}[leftmargin=*]
\item \textbf{RegEvo} Local evolutionary search \citep{real2019regularized} using pointwise mutations of single parent sequences and crossover recombination of two parent sequences. 
\item \textbf{NN + RegEvo} An ablation of \textit{NN+MILP}, with the only difference being the use of \textit{RegEvo} in lieu of MILP for solving the acquisition problem. Here, the inner-loop solver is allowed 10k queries of the acquisition function batched over 100 rounds, and proposes the point it has visited with the highest acquisition function value. The surrogate model is fit exactly as in NN+MILP.
\item \textbf{Ensemble + RegEvo} A re-implementation of the `MBO' baseline from \citet{angermueller2020population}, using an ensemble of linear and random forest regressors as the surrogate, where hyper-parameters are dynamically selected at each iteration. The acquisition function is the ensemble mean and inner-loop optimization uses \textit{RegEvo}.
\item \textbf{RBFOpt} A competitive mixed-integer black-box optimization solver that uses the `Radial Basis Function method' as a surrogate model \citep{costa2018rbfopt}.
\end{itemize}

\begin{figure*}[ht]
\begin{center}
\includegraphics[width=0.33\textwidth]{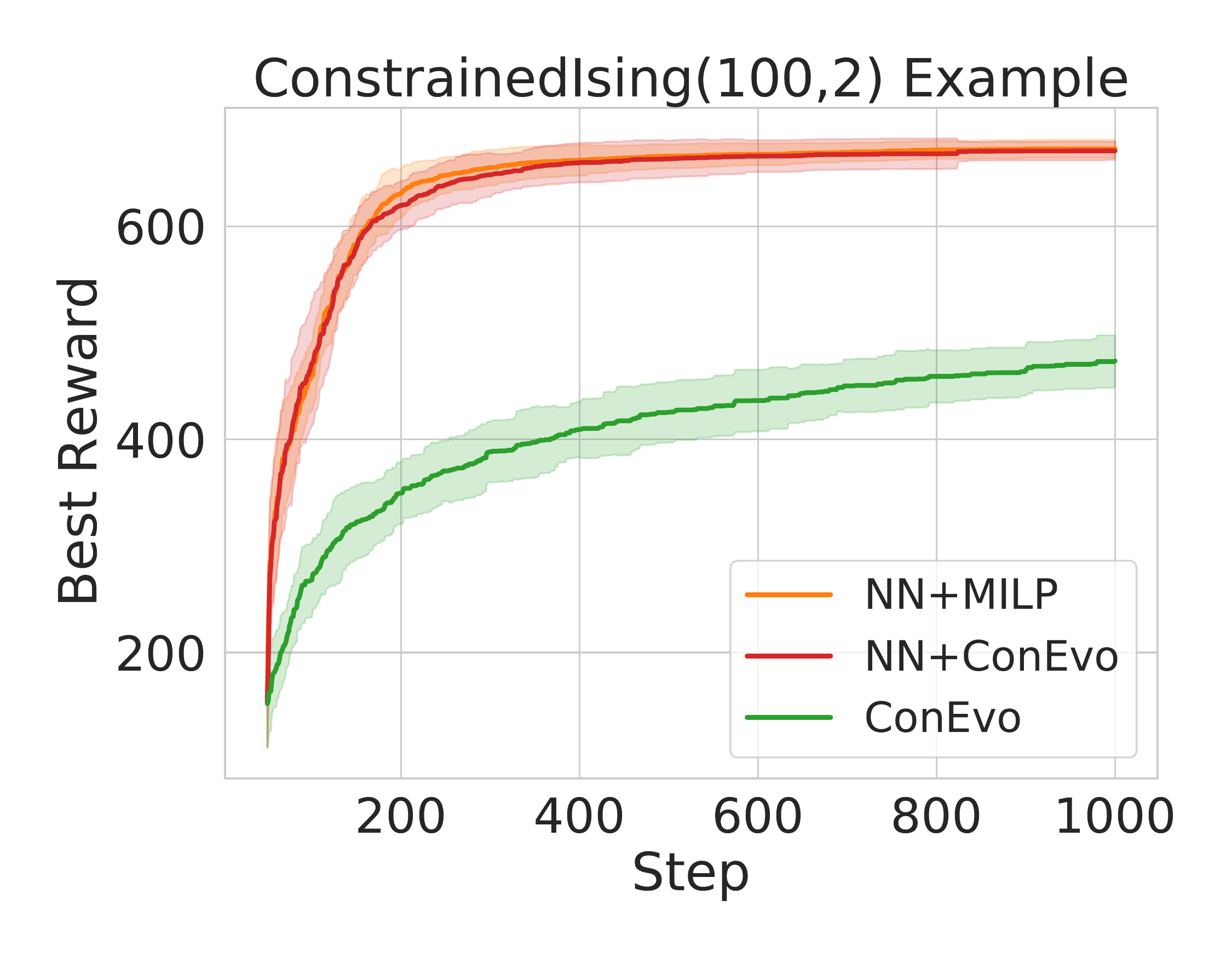}
\includegraphics[width=0.33\textwidth]{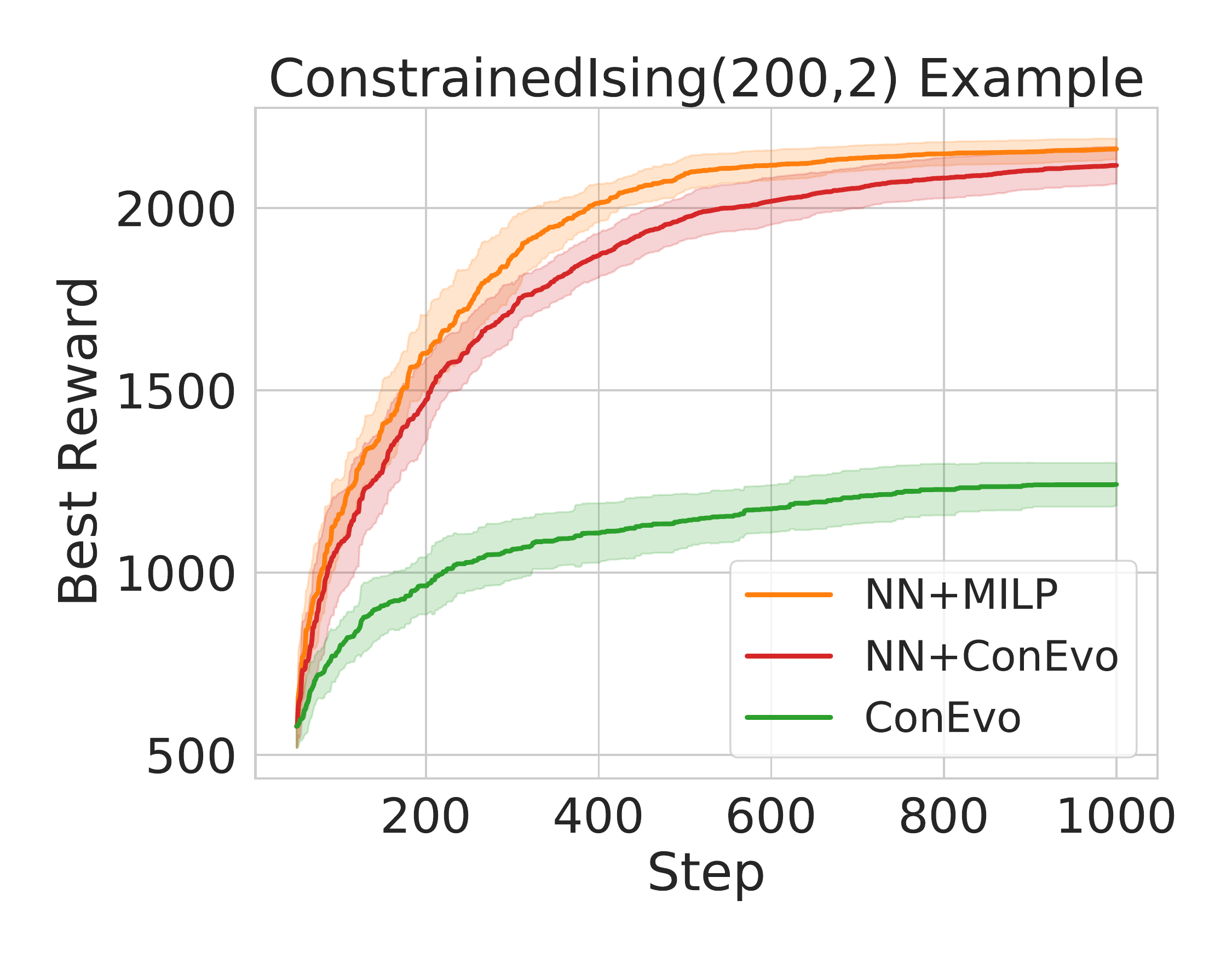}
\includegraphics[width=0.33\textwidth]{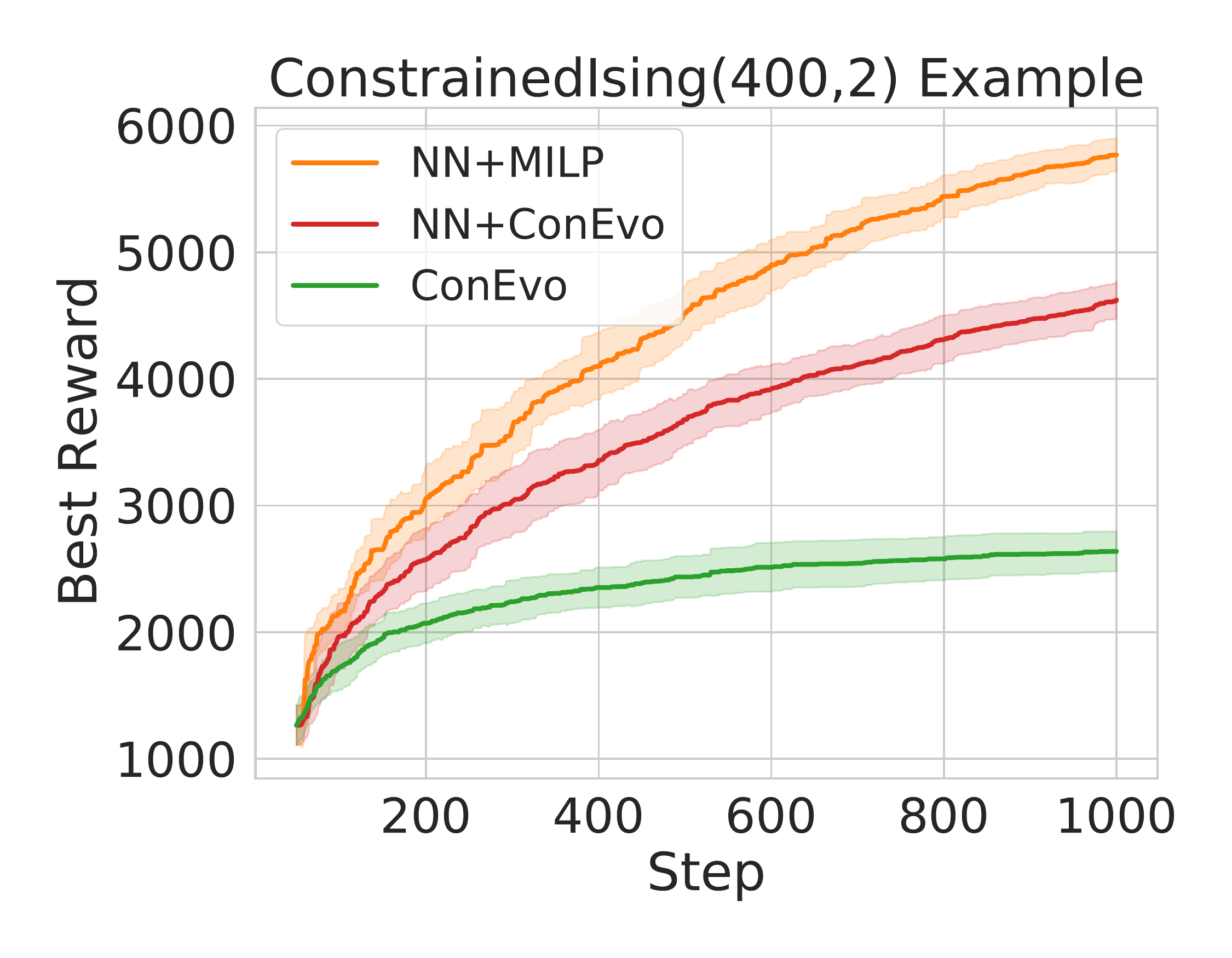}
\caption{Best observed reward as a function of iteration for an example constrained problem (Section~\ref{sec:exp_constr}) for each of n = 100, 200, and 400 (left-to-right). Lines and bands indicate the average and $\pm 1$ sd respectively, over 20 trials for $n=100$ and 10 trials for the rest. Distribution of normalized final scores and more examples can be found in Appendix~\ref{sec:appendix_results_constr}}.
\label{fig:constrained_examples}
\end{center}
\end{figure*}

Figure~\ref{fig:unconstrained_scores} plots the distribution of algorithms' scores for all unconstrained problems and an example reward curve from each class. We omit RBFOpt from the BBOB problems since it proposes the integer midpoint (rounded down) as part of its initialization, which is close to optimal by design (see Appendix~\ref{sec:appendix_tasks_bbob}). We observe that relative performance of algorithms varies significantly by objective family, with \textit{NN+MILP} performing well across the board. In particular, we wish to highlight the empirical benefits of global optimization of the acquisition function, as illustrated by the improved performance of \textit{NN+MILP} vs.\ \textit{NN+RegEvo}. The only difference between the two is the former's stronger optimality guarantees when solving the acquisition problem. We observe that \textit{NN+MILP} obtains a greater or equal score than its evolution-based counterpart in 22 of the 30 problems considered, and variance in its normalized scores is lower within a given objective family.

The comparison of \textit{NN+MILP} and \textit{Ensemble+RegEvo} solver is also instructive. Here, the primary difference is the hypothesis class $\mathcal{F}$. The strong performance of \textit{Ensemble+RegEvo} on TfBind, and to a lesser extent BBOB, suggests that ensembles of linear and tree-based regressors are better suited to approximate those black-box objectives. However, the combination of a single neural network surrogate and exact optimization yields comparable performance.

\subsection{Constrained Optimization}
\label{sec:exp_constr}

Next, problems are augmented with combinatorial constraints on the domain. We simulate fine-balance constraints in observational studies \citep{zubizarreta2018designmatch,Building-Representative-Matched-Samples}, where the same number of items must be selected from given sub-populations (e.g., sharing a common attribute). These simple, yet highly combinatorial, constraints allow for comparison with evolutionary algorithms that are designed to maintain feasibility with every mutation. 

We use a binary alphabet $\mathcal{A} = \{0, 1\}$ to indicate whether each of  $n$ items is selected. These have been partitioned into given equally-sized subsets $S_1, ..., S_{2k}$ for some integer $k$, and constraints enforce that the number of selected items is equal in pairs of subsets: $\sum_{i \in S_{2j-1}} z_{i1} = \sum_{i \in S_{2j}} z_{i1}$ for $j\in[k]$. $\text{Ising}(n,2)$ functions simulate the non-linear reward for a given selection. We create 30 problems by sampling 10 sets of Ising parameters for each of  $n \in \{100, 200, 400\}$, and setting $k=n/10$. See Appendix~\ref{sec:appendix_tasks} for details.

The following optimization approaches provide ablations to contrast declarative vs. procedural approaches to handling constraints. Configuration details are given in Appendix~\ref{sec:appendix_solvers}. 

\begin{itemize}[leftmargin=*]
    \item \textbf{ConEvo} \textit{RegEvo} with our own custom mutator that procedurally maintains feasibility. Paired subsets are mutated jointly, such that the number of changes in each pair is the same.
    \item \textbf{NN+ConEvo} An ablation of \textit{NN+MILP} where \textit{ConEvo} replaces MILP as the inner-loop solver. The inner-loop solver is allowed 10k queries of the acquisition function, batched over 100 rounds.
\end{itemize}

Figure~\ref{fig:constrained_examples} plots algorithms' best observed reward over time for a representative problem of each size. We observe that \textit{NN+MILP} and \textit{NN+ConEvo} significantly outperform \textit{ConEvo} for all problem sizes, owing to their ability to model the objective with a surrogate. The benefits of global optimization of the acquisition function are evident in the improved performance of \textit{NN+MILP} vis-a-vis \textit{NN+ConEvo} at larger scales; while the two model-based methods perform similarly for $n=100$, \textit{NN+MILP} improves considerably for $n=400$ and, to a lesser extent, $n=200$. We also note that neither method benefits significantly from using a larger surrogate network (see Appendix~\ref{sec:appendix_results_constr}), suggesting that the relatively small number of training points is a more significant bottleneck for approximation than surrogate capacity.

We emphasize that the two methods also differ in terms of ease of implementation. In particular, \textit{NN+MILP} required few extra lines of code to add subset-equality constraints to the existing MILP formulation, and could have just as easily been extended to other, possibly interacting, MILP-representable constraints. Conversely, \textit{NN+ConEvo} relied on a custom mutator tailored to the given structure, and may require significant reworking if other constraints are added.

\subsection{Practicality of MILP}
\label{sec:exp_pract}

\begin{figure}[ht]
    \centering
    \includegraphics[width=0.9\linewidth]{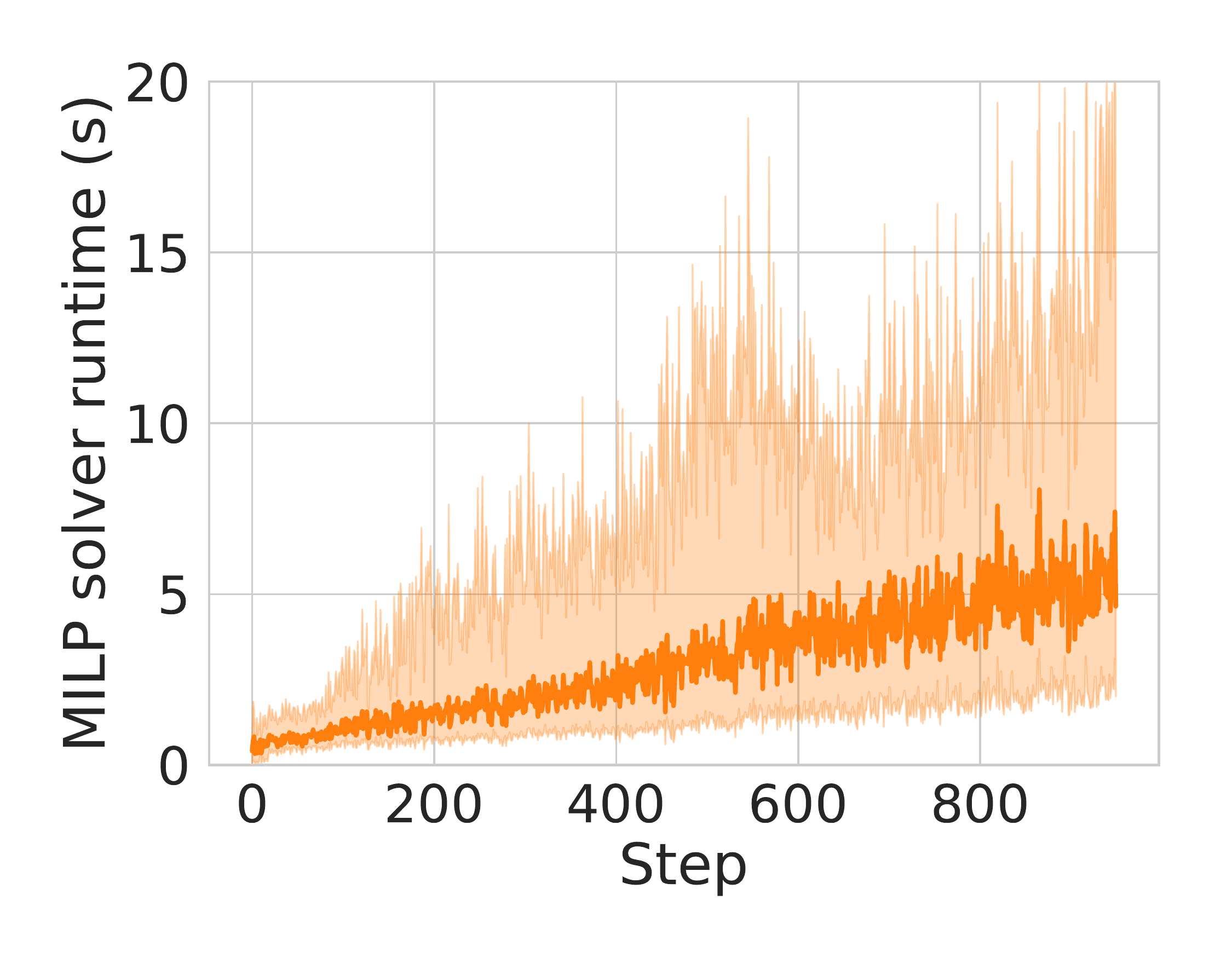}
    \caption{Distribution of MILP acquisition problem solve times as a function of iteration. Line and bands show the median and 5th/95th percentile range over all trials of all TfBind(8,4) problems (Section~\ref{sec:exp_unconstr}).}
    \label{fig:unconstrained_runtimes}
\end{figure}

Despite the computational complexity of the acquisition problem, MILP finds globally optimal solutions in seconds: the inner-loop optimization for \textit{NN+MILP} took $7.92 \pm 4.23$s (avg.\ $\pm$ sd) across all unconstrained experiments (Section~\ref{sec:exp_unconstr}) and $17.20 \pm 12.08$s for the largest ($n = 400$) constrained experiments (Section~\ref{sec:exp_constr}), never exceeding the time limit of 500s (i.e., all solutions were provably optimal). Note that \textit{NN+RegEvo/NN+ConEvo} are tuned to take comparable (or larger) time: $9.00 \pm 1.94$s and $56.80 \pm 15.18$s per step in the two settings above respectively. Figure~\ref{fig:unconstrained_runtimes} plots the distribution of MILP solve times for inner-loop optimization as a function of iteration for all TfBind problems, which seems to increase roughly linearly as no-good constraints are added. See Appendix~\ref{sec:appendix_results_pract} for other problem classes, and timing results when using larger surrogate networks.

\section{MINLPLib Case Study}
\label{sec:minlp}

Next, we apply our method to a class of linearly-constrained binary quadratic problems (BQPs) from MINLPLib \citep{minlplib} that contain practically-motivated constraints such as graph partitioning (\texttt{graphpart}), generalized assignment (\texttt{pb}), and shortest path (\texttt{qspp}). These problems are typically used to benchmark specialized white-box solvers that exploit known objective structure, but here we treat them as problems with black-box objectives and linear white-box constraints. While many black-box algorithms contain a sampling step that could be adapted to handle these constraints (e.g., \citet{oh2019combinatorial}), this may require specialization to each constraint type (much like \textit{ConEvo}), since finding feasible points through standard rejection sampling can be impractical (e.g., the chance is below $10^{-6}$ in \texttt{graphpart} instances). In contrast, we show that by using MILP to tackle constraints, our general method can often find feasible solutions that rival the best known solution $v^*$ found by a white-box solver (provided by MINLPLib).

We report the \textit{primal gap} of the best objective value $v$ found by \textit{NN+MILP} after 1000 steps, defined as $\frac{|v - v^*|}{max(|v|,|v^*|)}$, or 0 if $|v|=|v^*|=0$, or 1 if $v$ and $v^*$ have different signs \citep{berthold2013measuring}. MINLPLib contains 61 BQP problems with at least one linear constraint, with a number of binary variables ranging from 48 to 2203. We run 20 optimization trials per problem, each using a different initial dataset of 50 feasible points produced by solving MILPs with random objectives. We consider the same \textit{NN+MILP} configuration as in Section~\ref{sec:experiments}. Table~\ref{tab:minlp_summ} summarizes the proportion of all \textit{NN+MILP} trials achieving an optimality gap of 0\%, $\leq 1\%$ and $\leq 10\%$. Of note, we match $v^*$ in at least one trial for 20 of the 61 problems, spanning all classes of constraints, while we get within 10\% in an additional 25 problems. See Appendix~\ref{sec:app_minlplib} for details.

\begin{table}[]
\centering
\setlength{\tabcolsep}{4pt}
\caption{Proportion of MINLPLib trials where \textit{NN+MILP} achieved primal gap of 0\%, $\leq$1\%, $\leq$10\%, by problem class.}
\label{tab:minlp_summ}
\begin{tabular}{lccrrr}
\multicolumn{2}{c}{\multirow{2}{*}{\textbf{\begin{tabular}[c]{@{}c@{}}Class\\ (\# problems)\end{tabular}}}} & \multirow{2}{*}{\textbf{\begin{tabular}[c]{@{}c@{}}Range of\\ \# variables\end{tabular}}} & \multicolumn{3}{c}{\textbf{Proportion with gap}}                                             \\
\multicolumn{2}{c}{}                                                                                        &                                                                                           & \multicolumn{1}{c}{0\%} & \multicolumn{1}{c}{$\leq$ 1\%} & \multicolumn{1}{c}{$\leq$ 10\%} \\ \hline
\texttt{graphpart} & (31)  & [48,300]   & 20\% & 21\% & 59\% \\
\texttt{pb}        & (8)   & [525,600]  & 19\% & 38\% & 86\% \\
\texttt{qspp}      & (6)   & [180,420]  & 33\% & 72\% & 100\% \\
\texttt{other}     & (16)  & [50,2203]  &  3\% &  9\% & 27\% \\
\hline
\end{tabular}
\end{table}
\section{NAS-Bench-101 Case Study}
\label{sec:nas}

Finally, we use the NAS-Bench-101 \citep{nasbench101} neural architecture search (NAS) benchmark to illustrate the power of MILP's declarative constraint language in formulating complex combinatorial domains. The optimization domain consists of directed acyclic graphs (DAGs) representing the \textit{cell} in a neural architecture. Two nodes represent the input and output, and must be connected by a directed path, while the remaining nodes are each assigned to be 1x1 convolution, 3x3 convolution, or 3x3 max-pooling. Edges specify the flow of activations between nodes. The objective $f(x)$ is out-of-sample image classification accuracy. More details can be found in Appendix~\ref{sec:nasbench_formulation}.

We introduce a novel MILP formulation that precisely characterizes the set of valid NAS-Bench-101 cells. We use two sets of decision variables; the first set are binary and encode the upper-triangular adjacency matrix of a DAG with exactly $V$ nodes. The second set are a one-hot binary encoding of nodes' operations. Crucially, we introduce a new ``null'' operation, allowing the MILP to represent DAGs with fewer than $V$ nodes. Constraints enforce that all non-null nodes appear on a path from the input to output node, and that there exists at least one such path. A full formulation in terms of linear constraints appears in Appendix~\ref{sec:nasbench_formulation}, along with a variant to address certain graph isomorphisms.

\begin{figure}[ht]
    \centering
    \includegraphics[width=0.9\linewidth]{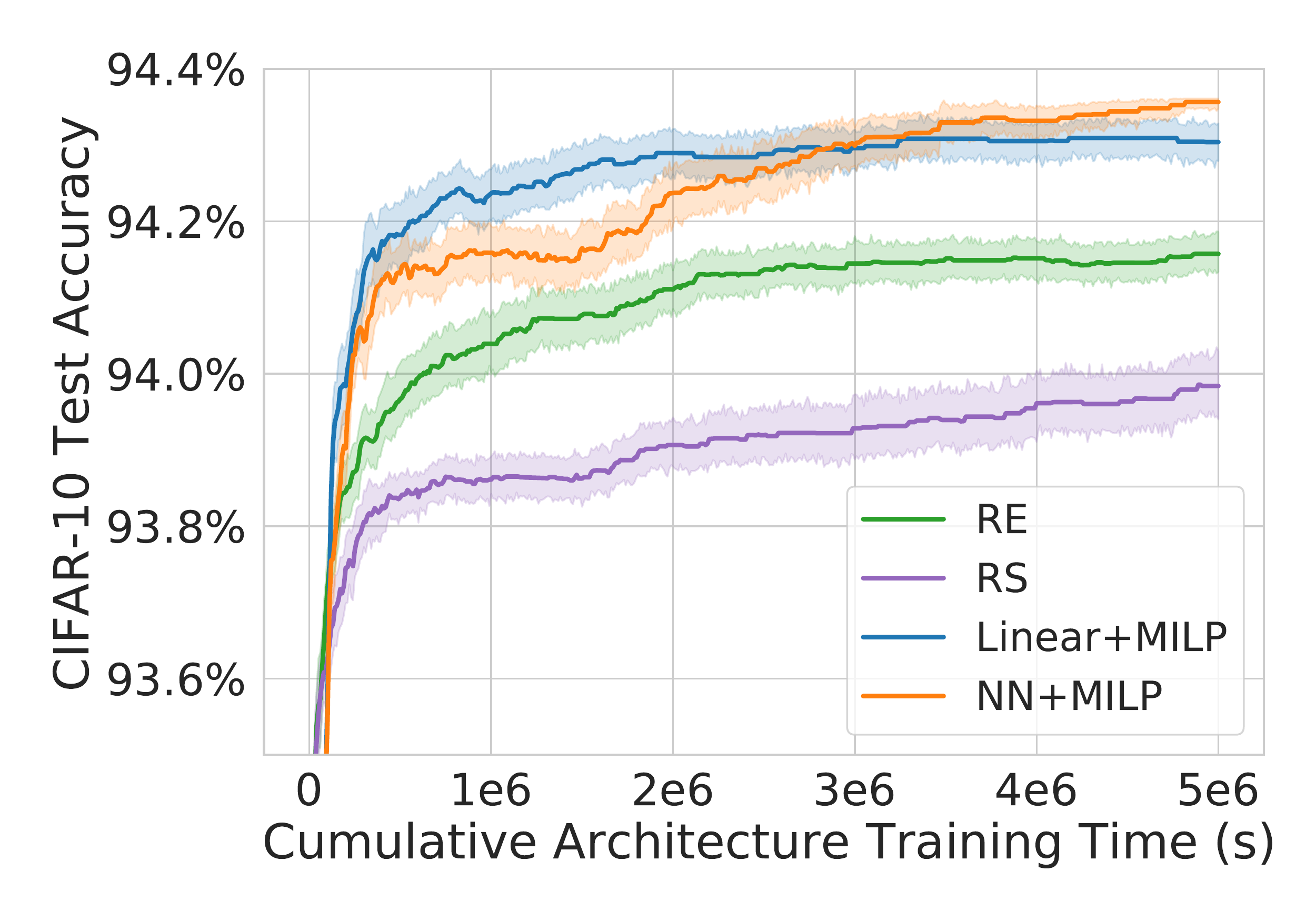}
    \caption{Test accuracy of algorithms' incumbent architecture as a function of cumulative training time on NAS-Bench-101, averaged over 100 trials. Bands indicate 95\% confidence interval for the mean.}
    \label{fig:nasbench_results}
\end{figure}

We use the same configuration of \textit{NN+MILP} as in Section~\ref{sec:experiments}, and include an ablation \textit{Linear+MILP} that replaces the surrogate by a linear model. The latter is trained on $\mathcal{D}_{t-1}$ with additional randomization provided by bootstrapping. Regularized evolution (\textit{RE}) and random search (\textit{RS}) baselines are from~\citet{nasbench101}. Figure~\ref{fig:nasbench_results} plots the out-of-sample accuracy of the proposed architecture with the highest observed validation accuracy (the ``incumbent'' architecture) vs.\ the cumulative architecture training time. \textit{NN+MILP}, despite its more general design, significantly outperforms \textit{RE}. Interestingly, \textit{Linear+MILP} outperforms \textit{NN+MILP} in early iterations, but is eventually overtaken. Future work could select among MILP-compatible models at each iteration.

\section{Conclusion and Future Work}
\label{sec:conclusion}

In this work we propose the \textit{NN+MILP} framework for  discrete MBO, using neural networks with ReLU activations for surrogate modeling and MILP to solve the acquisition problem. A major advantage of our method is its generality, using MILP's versatile declarative constraint language to address domains that might otherwise require specialized search algorithms for inner-loop optimization. Our experiments show that \textit{NN+MILP} performs well on a range of discrete black-box problems with practical computational overhead using standard packages and hardware. If there is a need for faster runtimes, one could devise problem-specific heuristics to use as warm-start for the acquisition MILP. 

MILP's versatility also suggests several interesting directions for future work. More complex acquisition functions could be considered to manage the exploration-exploitation trade-off as long as they remain MILP-representable, e.g., Expected Improvement defined over the posterior predictive distribution of an ensemble of ReLU networks. Alternatively, one could parameterize the Hamming Distance exclusion radius in the no-good constraints, which could be increased or decreased dynamically across iterations to encourage more exploration or exploitation respectively. For future applications to continuous or mixed-integer domains, the question arises as to how to best avoid redundant proposals given that no-good constraints cannot be applied as stated.

\bibliography{iclr2022_conference}
\bibliographystyle{iclr2022_conference}

\clearpage

\appendix

\section{Strengthening the MILP formulation for neural networks}
\label{sec:relu_formulation}

Here we discuss more advanced techniques for formulating the neural network surrogate model in the MILP problem. Recall the ReLU formulation constraints \eqref{eq:relu_0} and \eqref{eq:relu_1} from Section~\ref{sec:algo_opt}, except that we consider $M$ separately for each constraint:
\begin{align*}
    0 \leq y \leq  &\ M_0\alpha \tag{\ref{eq:relu_0}'}\label{eq:relu_0_m0}\\
    w^\top x + b \leq y \leq  &\ w^\top x + b + M_1(1 - \alpha),\tag{\ref{eq:relu_1}'}\label{eq:relu_1_m1}
\end{align*}
Here, we require a nonnegative value $M_0$ such that the right-hand side of~\eqref{eq:relu_0_m0} is greater or equal than a valid upper bound on $y$ when $\alpha = 1$. Similarly, $M_1$ must be a nonnegative value such that the right-hand side of~\eqref{eq:relu_1_m1} is greater or equal than zero when $\alpha = 0$. Therefore, we may choose $M_0$ to be any upper bound of $\max_{x \in \Omega'} w^\top x + b$ and $M_1$ to be any upper bound of $\max_{x \in \Omega'} -(w^\top x + b)$, where $\Omega'$ is the domain of the inputs of this ReLU, which depends on $\Omega$. The tighter these bounds are, the better the MILP performs.

Moreover, if we find negative $M_0$ or $M_1$, then we may (in fact, must) replace the formulation by $y = 0$ or $y = w^\top x + b$ respectively, since in these cases the ReLU is always inactive or active for any $x \in \Omega'$. This replacement must be done because the formulation assumes nonnegative $M_0$ and $M_1$ for feasibility.

The simplest way to compute $M_0$ and $M_1$ is to start from bounds in $\Omega$ and propagate them via interval arithmetic. For example, if $x \in [L, U]$, then $M_0$ can be set to $\sum_{i : w_i > 0} w_i U_i + \sum_{i : w_i < 0} w_i L_i + b$ and $M_1$ to $-(\sum_{i : w_i > 0} w_i L_i + \sum_{i : w_i < 0} w_i U_i + b)$. However, despite being fast, the drawback of this simple approach is that it does not take into account constraints on $\Omega$ or one-hot and no-good constraints. 

In our experiments, we compute $M_0$ and $M_1$ by solving the linear programming (LP) relaxations of $\max_{x \in \Omega'} w^\top x + b$ and $\max_{x \in \Omega'} -(w^\top x + b)$ respectively (i.e., without integrality constraints). We remark that for neurons in the same layer these LPs have the same constraints but different objectives, and thus we may take advantage of the warm starting functionality in LP solvers. While this requires solving two LPs per neuron, taking into account the constraints from $\Omega$ into the bounds often enable the overall MILP to be solved much faster.

The formulation can also be strengthened with cutting plane techniques \citep{anderson2020strong}, but they are not particularly beneficial for the small network sizes considered in this paper (at most two layers with 16 ReLUs each) and thus we do not add them. Future work could explore warm-starting the MILP solver using results from earlier MBO iterations or problem-specific heuristics.

\section{Benchmarking Tasks}
\label{sec:appendix_tasks}

This section details the black-box objective functions considered in both unconstrained (Section~\ref{sec:exp_unconstr}) and constrained (Section~\ref{sec:exp_constr}) experiments. Recall that all objective functions are defined over fixed-length discrete vectors of length $n$, with each element drawn from an alphabet $\mathcal{A}$ of fixed size.

\subsection{TfBind}
\label{sec:appendix_tasks_tfbind}

The objective function is given by the binding affinity of a length-8 DNA sequence to a particular transcription factor, characterized experimentally in the dataset described by \citet{barrera2016survey}. The problem size is thus fixed by the application at hand, with $n = 8$ and $|\mathcal{A}| = 4$ (each input element corresponding to a given DNA nucleotide).  We min/max-normalize the binding affinity values for each factor to the zero-one interval. We create 12 unconstrained problems (Section~\ref{sec:exp_unconstr}) using the following datasets: CRX\_R90W\_R1, CRX\_REF\_R1, FOXC1\_REF\_R1, GFI1B\_REF\_R1, HOXD13\_Q325R\_R1, HOXD13\_REF\_R1, NR1H4\_C144R\_R1, NR1H4\_REF\_R1, PAX4\_REF\_R1, PAX4\_REF\_R2, POU6F2\_REF\_R1, SIX6\_REF\_R1. Here, the 3 fields separated by underscores represent the transcription factor id, any mutations that have been made to the transcription factor, and the id of the experimental replicate used when collecting data. 

\subsection{RandomMLP}
\label{sec:appendix_tasks_randommlp}

The objective function is given by the output of a multi-layer perceptron (MLP) with randomly-sampled weights. Different functions are generated by varying the architecture type (described below) and random seed. All architectures employ a one-hot encoding of the inputs as the first layer. Weights are sampled using the default behavior of tf.keras.layers.Dense (\texttt{glorot\_uniform}).

We consider two architecture types, both utilizing more layers/parameters than the 16-neuron networks used by \textit{NN+MILP} (Section~\ref{sec:experiments}). The \textit{RandomFCC} architecture uses two fully-connected layers with 128 hidden units each, while the \textit{RandomCNN} architecture uses two convolutional layers each with 64 hidden units each, a kernel width of 13 and stride size of 1. We use a linear activation function for the output and ReLU activations for all intermediate layers. 

Unconstrained \textit{RandomMLP} problems (Section~\ref{sec:exp_unconstr}) all have size $n = 25$ and $|\mathcal{A}| = 5$ . Eight objective functions are created by varying the architecture type (FCC or CNN) and random seed (0, 13, 42, 77). 

\subsection{BBOB}
\label{sec:appendix_tasks_bbob}

The objective is given by a function from the continuous Black-Box Optimization Benchmarking library \citep{hansen2009real}. All BBOB functions are defined for a variable number of dimensions $n$ and the search domain is given as $[-5, 5]^n$, with the global optimum centered at zero. We normalize each function's output range by evaluating it at $30$ fixed points and dividing outputs by the median absolute deviation in those points' values. 

We discretize functions for our setting (Section~\ref{sec:problem_setting}) by defining a grid over the continuous search domain, adjusted so that the optimal solution exactly corresponds to a point in the grid. Concretely, we use a fixed alphabet $\mathcal{A} = \{1, \ldots, m\}$ for all coordinates, denoting the \textit{index} of one of $m$ allowed values for that coordinate. Allowed values for each coordinate are $m$ equally-spaced points in the range $[-5, 5]$, except for a point lying closest to zero which is overwritten to exactly equal that value. In this way, the optimum is guaranteed to lie on the discretized grid. Note that, despite the underlying continuous structure, all algorithms treat each dimension as an unordered, categorical variable.

For unconstrained \textit{BBOB} problems (Section~\ref{sec:exp_unconstr}), we select a diverse set of objectives by taking two functions from each of the five categories defined by the BBOB library:
\begin{enumerate}
\setlength{\itemsep}{-2pt}
    \item Separable functions: Sphere (SPHERE) and Ellipsoidal (ELLIPSOID\_SEPARABLE).
    \item Functions with low or moderate conditioning: Attractive Sector (ATTRACTIVE\_SECTOR) and Step Ellipsoidal (STEP\_ELLIPSOID).
    \item Functions with high conditioning and unimodal: Discus (DISCUS) and Bent Cigar (BENT\_CIGAR).
    \item Multi-modal functions with adequate global structure: Weierstrass (WEIERSTRASS) and Schaffers F7 (SCHAFFERS\_F7).
    \item Multi-model functions with weak global structure: Schwefel (SCHWEFEL) and Gallagher's Gaussian 21-hi Peaks (GALLAGHER\_21ME).
\end{enumerate}
We set the dimension for all of these to $n = 10$ and discretize as described above, using an alphabet of size $|\mathcal{A}| = 10$ for all coordinates. We purposefully use a relatively large alphabet to ensure that the discretization does not obscure any inherent variance across a given coordinate.

\subsection{Ising}
\label{sec:appendix_tasks_ising}

The objective computes the negative energy of fully-connected binary Ising Model with pairwise potentials drawn i.i.d. from a standard Gaussian. Binary decision variables (i.e., $\mathcal{A}=\{0,1\}$) represent the spins of $n$ particles in the system, which are treated as nodes in a fully-connected graph. Each edge of the graph is defined by a a $2 \times 2$ table of scores for each possible spin configuration of the nodes that are connected by the edge. All edge scores are drawn i.i.d. from a standard Gaussian, and the overall function is the sum of the scores over all edges. 

For the constrained experiments (Section~\ref{sec:exp_constr}) we create 30 problems by varying $n \in \{100, 200, 400\}$ and generating 10 different random instances of Ising model parameters for each $n$. These are combined with the subset-equality constraints (defined in Section~\ref{sec:exp_constr}), setting the number of paired subsets to $k=\frac{n}{10}$ (i.e., the cardinality of subsets is 5, regardless of $n$). 

\section{Baseline Optimization Algorithms}
\label{sec:appendix_solvers}

In this section we describe implementation and configuration details for all baseline optimization algorithms described in Section~\ref{sec:experiments}.

\subsection{NN+MILP}
\label{sec:appendix_solvers_nnmilp}

For our experiments, we implement our main algorithm (Section~\ref{sec:algo}) as follows: we use a fixed surrogate model hypothesis class $\mathcal{F}$ of networks with a single, fully-connected hidden layer of 16 neurons. Models are trained with TensorFlow~\citep{abadi2016tensorflow}, using the ADAM optimizer for $25K$ epochs with a batch size of 64 and no explicit regularization. We use a constant learning rate of $\alpha = 0.01$ and default decay parameters $(\beta _1, \beta_2) = (0.9, 0.999)$. No hyper-parameter tuning is performed across problems. Model training is randomized due to the random example ordering of SGD training and random parameter initialization. The MILP acquisition problem is solved with the Mixed-Integer Programming solver SCIP 7.0.1 \citep{gamrath2020scip} using default settings and a time limit of 500 seconds. In order to increase the diversity of trained models, we train each model from scratch at each iteration of optimization instead of fine-tuning a model from an earlier iteration. 
\subsection{RegEvo}
\label{sec:appendix_solvers_regevo}

We re-implement the local evolutionary search algorithm of \citet{real2019regularized}, and extend the set of mutation operators from just pointwise mutators to also include a crossover operation that re-combines two parent sequences. The algorithm proposes $x_{t+1}$ by selecting two parent sequences from the existing population, recombining them and mutating them. Parents are chosen by tournament selection, taking the two best samples from a randomly-selected subset of size $T$ of previously sampled points. The pool from which parents can be selected is limited to the $D$ most recently-proposed points (referred to as the ``alive population''), to avoid high-reward points from early rounds dominating the process. The selected parent sequences are recombined by copying them left-to-right, starting a pointer at one parent at switching reading to the other parent with a fixed cross-over probability $p_c$ after each copy. The resulting sequence is finally mutated by changing each position to a different token from $\mathcal{A}$ with a fixed probability $p_m$.

In the unconstrained experiments (Section~\ref{sec:exp_unconstr}), we use \textit{RegEvo} as the outer-loop optimization algorithm and set the tournament size to $T = 10$, the alive population size to $D = 100$, and the crossover/mutation probabilities to $(p_c, p_m) = (0.1, 0.1)$. 

\subsection{NN+RegEvo}
\label{sec:appendix_solvers_nnregevo}

This algorithm is an ablation of \textit{NN+MILP}, with the only difference being the use of \textit{RegEvo} in lieu of MILP to solve the acquisition problem at every iteration. A surrogate neural network $\hat{f}_t \in \mathcal{F}$ is trained as in \textit{NN+MILP}, and the acquisition function is $a(x) = \hat{f}_t(x)$. The problem of selecting $x_{t+1}$ is posed as a \textit{batched} optimization problem and solved by \textit{RegEvo}. 

More concretely, at iteration $t$, the acquisition function is evaluated for all points in the existing population $\mathcal{D}_t$ to generate the initial inner-loop population $\hat{\mathcal{D}}_t := \{x_i, a(x_i)\}$. This population is iteratively extended by generating candidate proposals with \textit{RegEvo} in batches of size $b$, and with rewards now corresponding to the value of the acquisition function rather than the original black-box function. That is, \textit{RegEvo} generates $b$ points by recombination/mutation of parents from $\hat{\mathcal{D}}_t$, which are evaluated on the acquisition function and added to the inner-loop population. The process repeats until a total of $B$ candidates have been generated, at which point the one with the highest acquisition function value (excluding any points already proposed) is proposed as $x_{t+1}$.

In the unconstrained experiments (Section~\ref{sec:exp_unconstr}), we use $\textit{NN+RegEvo}$ and set surrogate model hyper-parameters exactly as in \textit{NN+MILP} (Section~\ref{sec:appendix_solvers_nnmilp}). For the inner-loop optimizer, we set the total number of acquisition function evaluations to $B = 10,000$ and batch size to $b = 100$. The \textit{RegEvo} optimizer's hyper-parameters, defined in Section~\ref{sec:appendix_solvers_regevo}, are set to $T = 20$, $D = 1,000$ and $(p_c, p_m) = (0.2, 0.01)$. 

\subsection{Ensemble+RegEvo}
\label{sec:appendix_solvers_ensembleregevo}
We recreate the \textit{MBO} baseline of \citet{angermueller2020population}. Here, surrogate modeling proceeds by optimizing the hyper-parameters of a diverse set of regressor models through randomized search. Regressors are trained using the \texttt{scikit-learn} libary \citep{pedregosa2011scikit}, drawing from the following model classes (randomized search parameters are listed in parentheses):

\begin{itemize}
    \setlength{\itemsep}{0pt}
    \item LassoRegressor (alpha)
    \item RidgeRegressor (alpha)
    \item RandomForestRegressor (max\_depth, max\_features, n\_estimators)
    \item LGBMRegressor (learning\_rate, n\_estimators)
\end{itemize}

Each model is evaluated by an explained variance score using five-fold cross validation on the training set. All models with a score $\geq 0.4$ are used as an ensemble for the surrogate model, with their \textit{average} prediction serving as the acquisition function. The acquisition problem is solved by batched \textit{RegEvo} with a total of $B = 12,500$ acquisition function evaluations and a batch size of $b = 25$. The optimizer's hyper-parameters, defined in Section~\ref{sec:appendix_solvers_regevo}, are set to $T = 20$, $D = 1,000$ and $(p_c, p_m) = (0.2, 0.01)$. 

We use \textit{Ensemble+RegEvo} in both the unconstrained (Section~\ref{sec:exp_unconstr}) and constrained (Section~\ref{sec:exp_constr}) experiments. In the latter case, we use the algorithm as a baseline that makes use of the declarative definition of constraints; during training of the ensemble, infeasible points are assigned a highly negative reward (worse than any observed). In this way, the surrogate model might be expected to implicitly model infeasibility with low predictions which should be avoided by the inner-loop optimizer.

\subsection{RBFOpt}
\label{sec:appendix_solvers_rbfopt}

\emph{RBFOpt}~\citep{costa2018rbfopt} is a black-box optimization solver for mixed-integer unconstrained problems (i.e., with only bound constraints) that performs competitively with respect to other solvers of its type. It uses a Radial Basis Function as a surrogate model and includes a number of practical enhancements. It relies on a mixed-integer nonlinear programming (MINLP) solver, BONMIN~\citep{bonami2008algorithmic}, to optimize the inner loop problems. The MINLP solver could in theory incorporate constraints in a similar fashion as in our work, although this is not offered by the open-source implementation (aside from manually penalizing the objective function) and we expect it to not scale as well as a MILP solver in practice since MINLP is a significantly more difficult problem class than MILP.

We use \emph{RBFOpt} for our unconstrained experiments (Section~\ref{sec:exp_unconstr}), using the open-source implementation available at \url{https://github.com/coin-or/rbfopt}.
We leave all settings at their defaults, including building the initial set of points. By default, RBFOpt uses a one-hot encoding for the categorical variables, and for all problems we mark them as categorical. As we note in the main text, we omit the RBFOpt results for BBOB because RBFOpt proposes the midpoint of the integer representation (rounded down) as part of its initialization, which is close to the optimal solution. 

\subsection{ConEvo}
\label{sec:appendix_solvers_conevo}

In Section~\ref{sec:exp_constr} we introduce \textit{ConEvo}, a local evolutionary search algorithm that exploits the known combinatorial structure of the subset-equality constraints considered therein. The method selects just a single parent sequence (using the same tournament procedure as \textit{RegEvo}) and mutates it in a way that guarantees feasibility of the child sequence. We do not implement recombination of multiple parent sequences since they are not likely to maintain feasibility. We described the application-specific mutator below.

Recall that the domain encodes the selection or not of each of $n$ items using a binary alphabet $\mathcal{A} = \{0,1\}$. The items' indices are partitioned into disjoint, equally-sized subsets $S_1, \ldots, S_{2k}$ for some $k$ and the constraints enforce that the number of selected items should be the same in pairs of subsets; that is:

\begin{align*}
    \sum_{i \in S_{2j-1}} \mathbb{I}\{x_i = 1\} = \sum_{i \in S_{2j}} \mathbb{I}\{x_i = 1\} \quad \forall j \in [k]
\end{align*}

where we have used indicator notation and the original decision variables $x$ rather the one-hot encoding.

The mutator begins with a single parent sequence $x$, assumed feasible, and is given access to the item subsets $S_1, \ldots, S_{2k}$. Each pair of subsets $(S_{2j-1}, S_{2j})$ is mutated concurrently to create the child sequence $y$, ensuring that mutations to one subset are counter-balanced by mutations to the second. Concretely, one of the two subsets is chosen randomly to be the ``independent'' mutatee with equal probability. We denote the selected subset $\mathcal{I}^+$, and the other subset in the pair by $\mathcal{I}^-$. Each position $i \in \mathcal{I}^+$ of the child sequence is flipped from its parent value with some fixed probability $p_m$. We compute $c$, the \textit{net} number of 0-to-1 conversions in positions $\mathcal{I}^+$. If $c$ is positive (i.e., there were more 0-to-1 conversions than 1-to-0 conversions) then exactly $c$ indices are chosen randomly from $\{i \in \mathcal{I}^-: x_i = 0\}$, and also flipped in the child. If $c$ is negative, then $-c$ indices are selected randomly from $\{i \in \mathcal{I}^-: x_i = 1\}$ and flipped in the child. If $c$ is zero, the positions in $\mathcal{I}^-$ are left unchanged in the child. As a result, the subsets $S_{2j-1}$ and $S_{2j}$ retain exactly the same number of selected items in the mutated sequence.

In the constrained experiments (Section~\ref{sec:exp_constr}), we use \textit{ConEvo} as the outer-loop optimization algorithm setting the tournament size for selecting parent sequences to $T = 20$ and the mutation probability to $p_m = 0.05$.

\subsection{NN+ConEvo}
\label{sec:appendix_solvers_nnconevo}

This algorithm is an ablation of \textit{NN+MILP} used in Section~\ref{sec:exp_constr}, with the only difference being the use of \textit{ConEvo} in lieu of MILP to solve the \textit{constrained} acquisition problem at every iteration. A surrogate neural network $\hat{f}_t \in \mathcal{F}$ is trained as in \textit{NN+MILP}, and the acquisition function is $a(x) = \hat{f}_t(x)$. Selecting $x_{t+1}$ is posed as a \textit{batched} optimization problem and solved by \textit{ConEvo}. The batched optimization procedure is exactly as described for \textit{NN+RegEvo} (Section~\ref{sec:appendix_solvers_nnregevo}).

In the constrained experiments (Section~\ref{sec:exp_constr}), we use $\textit{NN+ConEvo}$ and set surrogate model hyper-parameters exactly as in \textit{NN+MILP} (Section~\ref{sec:appendix_solvers_nnmilp}). For the inner-loop optimizer, we set the total number of acquisition function evaluations to $B = 10,000$ and batch size to $b = 100$. The \textit{ConEvo} optimizer's hyper-parameters, defined in Section~\ref{sec:appendix_solvers_conevo}, are set to $T = 20$ and $p_m = 0.05$.

\section{Binary vs Integer Variables}
\label{sec:appendix_binary_vs_integer}

In this work, we focus on problems with binary domain formulations (e.g., one-hot encoding of categorical domains), and even problems with integer variables such as the discretized BBOB are binarized. Part of the reason is to allow no-good constraints as described in Section~\ref{sec:algo_opt}, but in addition we have experimentally observed that the method performs better when using a binary encoding instead of an integer one.

When running this algorithm for unconstrained (bounded) integer or continuous problems, we have informally observed that our method frequently proposes solutions where several variable values are at either their lower bound or upper bound. As a result, our method would underexplore solutions away from the boundary. A possible explanation for this is that feedforward ReLU networks tend to extrapolate linearly, and thus their optima may often lie on the boundary~\citep{xu2021neural}. In contrast, every feasible point of a binary problem lies on a corner of the 0-1 hypercube. A similar observation has been made in the context of IDONE~\citep{bliek2021black}, which also uses a ReLU-based surrogate model: encoding the Rosenbrock problem using binary variables improves the performance of the IDONE algorithm, although the opposite happens for a Bayesian optimization algorithm~\citep{karlsson2020continuous}.

\begin{figure*}[ht]
    \centering
    \includegraphics[width=0.48\linewidth]{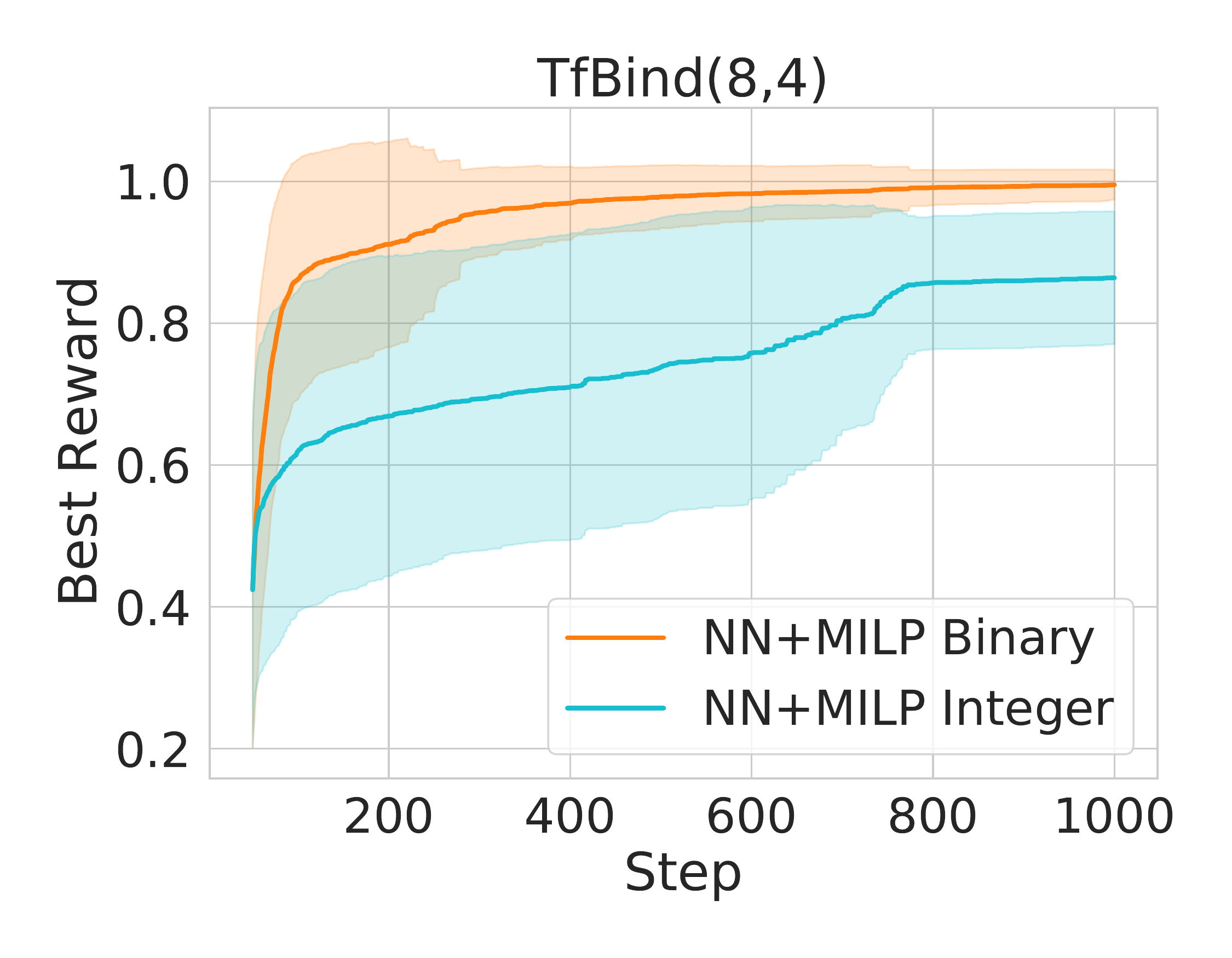}
    \includegraphics[width=0.48\linewidth]{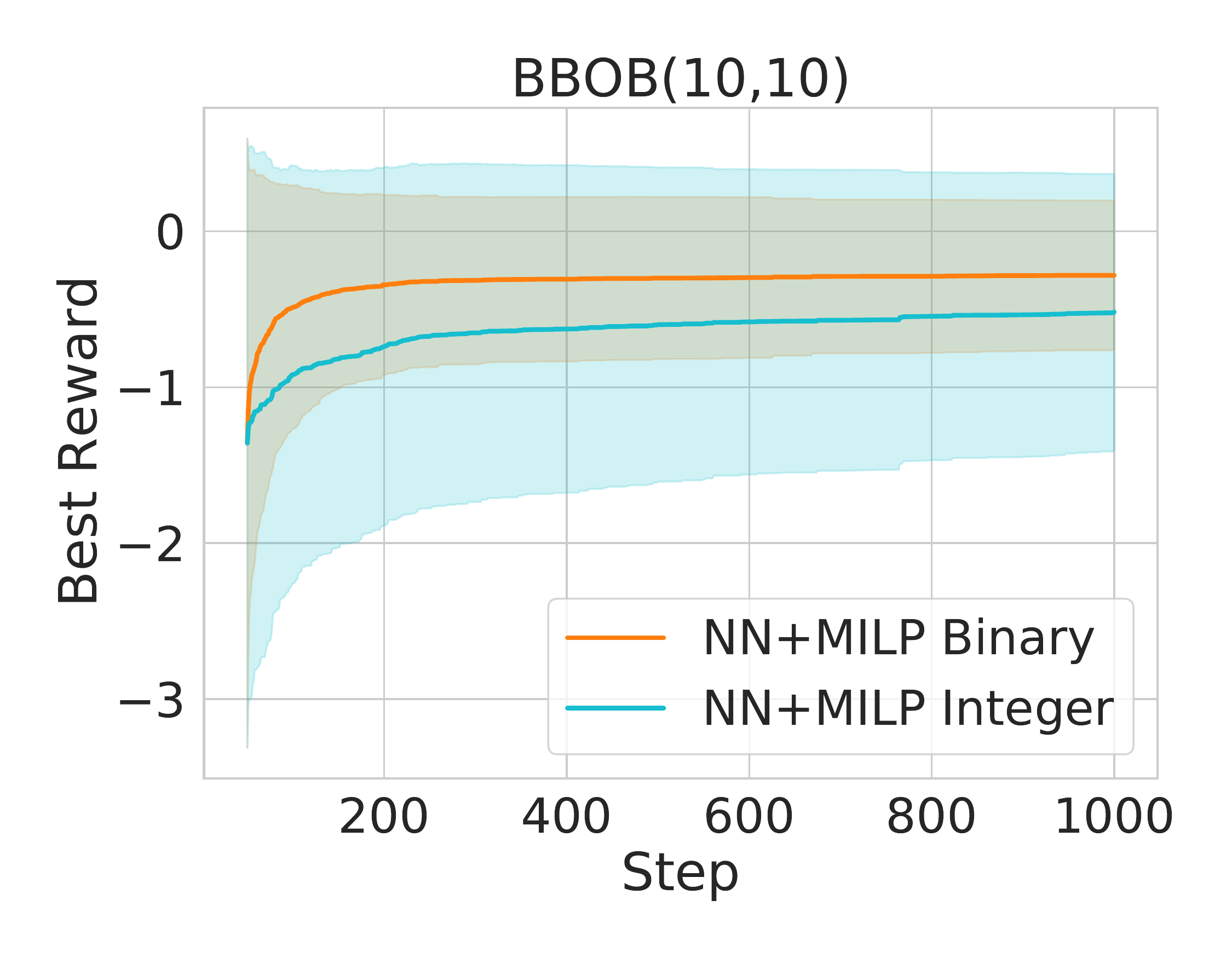}
	\caption{Best observed reward as a function of iteration for all TfBind (top) and BBOB (bottom) instances, comparing the use of binary and integer variables. For TfBind, the categorical variables are transformed to integer with an arbitrary ordering, and for BBOB, we use the given ordering of the problem. Note that the error region is large here since we aggregate all of the instances of each class.}
	\label{fig:integer_vs_binary}
\end{figure*}

We provide computational evidence of this behavior in Figure~\ref{fig:integer_vs_binary} for TfBind and BBOB instances, with the same experiment setup as Section~\ref{sec:exp_unconstr}. The binary variables are encoded as one-hot variables, whereas the integer variables follow an arbitrary ordering for TfBind and the problem ordering for BBOB.

\begin{figure*}[htbp]
\begin{center}
\centerline{\includegraphics[width=\textwidth]{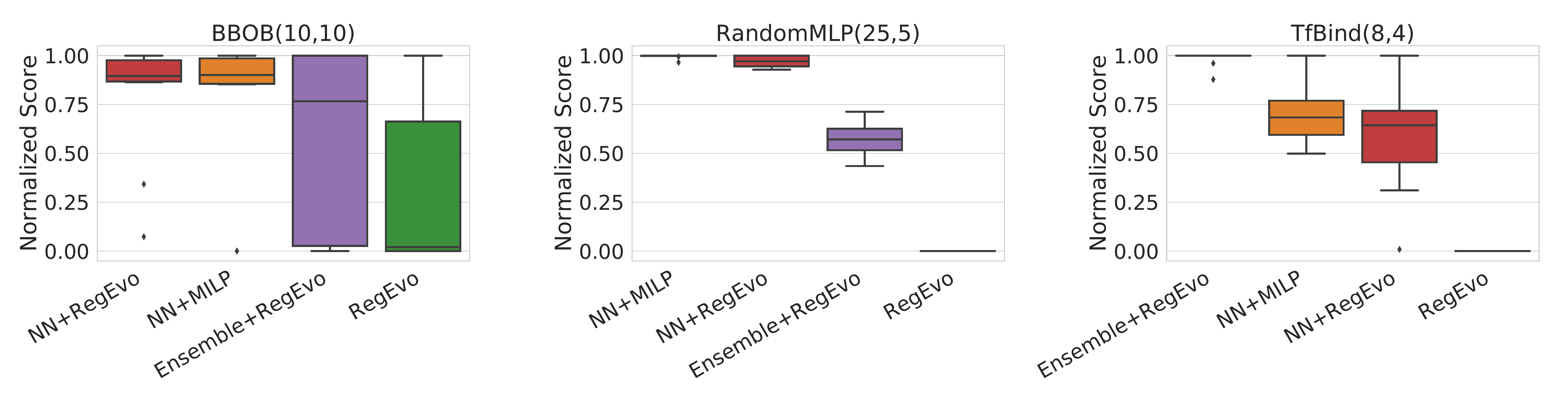}}
\caption{Distribution of algorithms' normalized AUC scores (Section~\ref{sec:appendix_results_unconstr}) on all unconstrained problems from Section~\ref{sec:exp_unconstr}, split by objective function class. Higher is better. We observe that relative performance of algorithms in terms of AUC is qualitatively similar as best-observed reward (Figure~\ref{fig:unconstrained_scores}).}
\label{fig:appendix_unconstrained_scores_auc}
\end{center}
\end{figure*}

\section{Additional Experiments}
\label{sec:appendix_results}

\subsection{Unconstrained Optimization}
\label{sec:appendix_results_unconstr}

\subsubsection{Normalized Area Under the Curve (AUC)}

While the best observed reward in $\mathcal{D}_N$ (i.e., after all evaluations) is the primary metric of comparison for algorithms per Section~\ref{sec:problem_setting}, it is also instructive to consider a measure of how fast algorithms converge to their best observed reward. To this end, we define an AUC metric that computes the \textit{area under the best observed reward curve}; higher values indicate that an algorithm found better points in earlier iterations. To facilitate comparison across problems, we min/max normalize algorithms' AUC scores within each problem exactly as we did for the best observed reward (Section~\ref{sec:exp_tasks}). That is, the best (resp. worst) on-average algorithm in terms of AUC is assigned a score of one (resp. zero) and intermediate values express relative distance from these extremes.

\begin{figure*}[ht]
\begin{center}
    \begin{multicols}{3}
    \includegraphics[width=\linewidth]{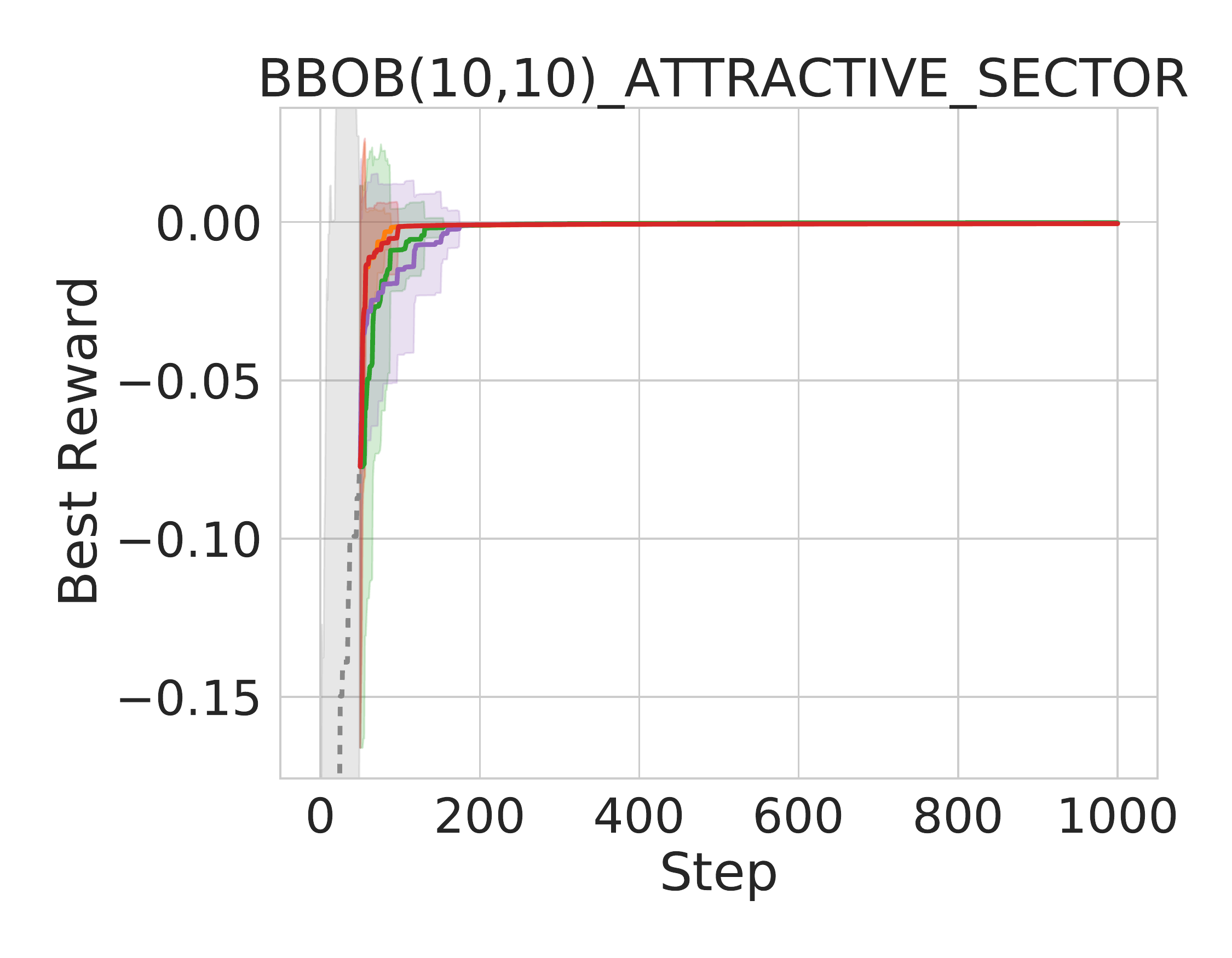}
    \includegraphics[width=\linewidth]{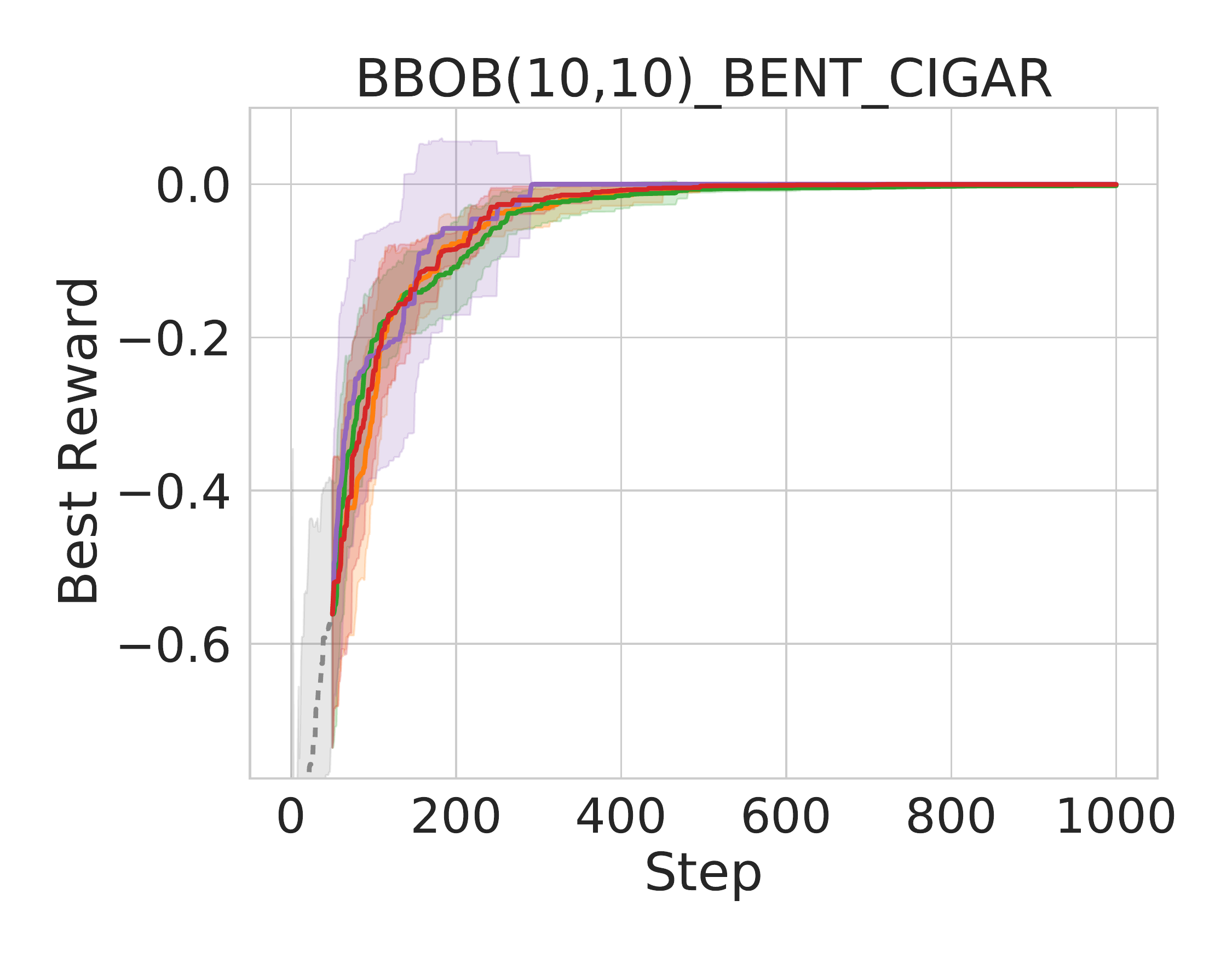}
    \includegraphics[width=\linewidth]{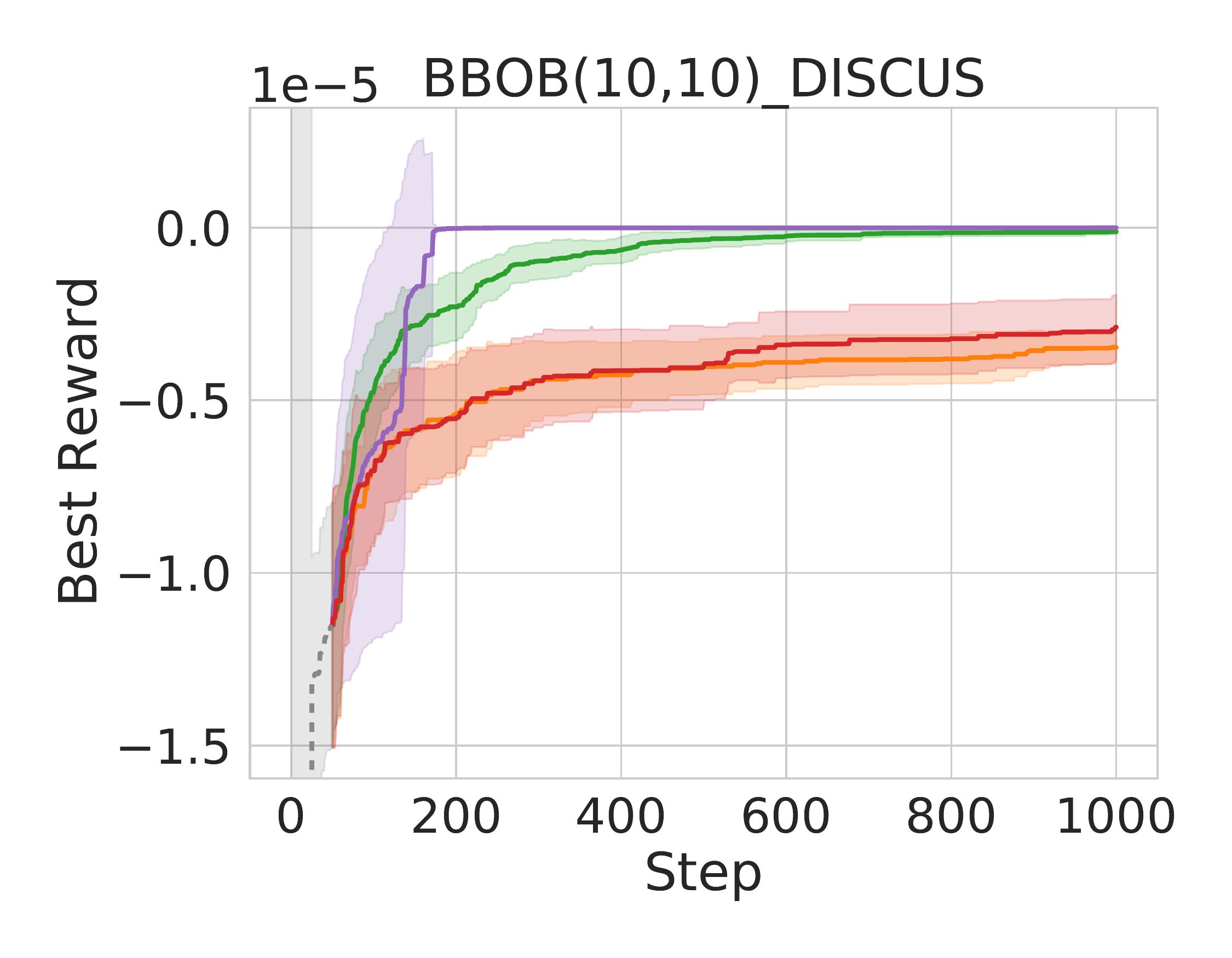}
    \end{multicols}
    \vspace{-30pt}
    \begin{multicols}{3}
    \includegraphics[width=\linewidth]{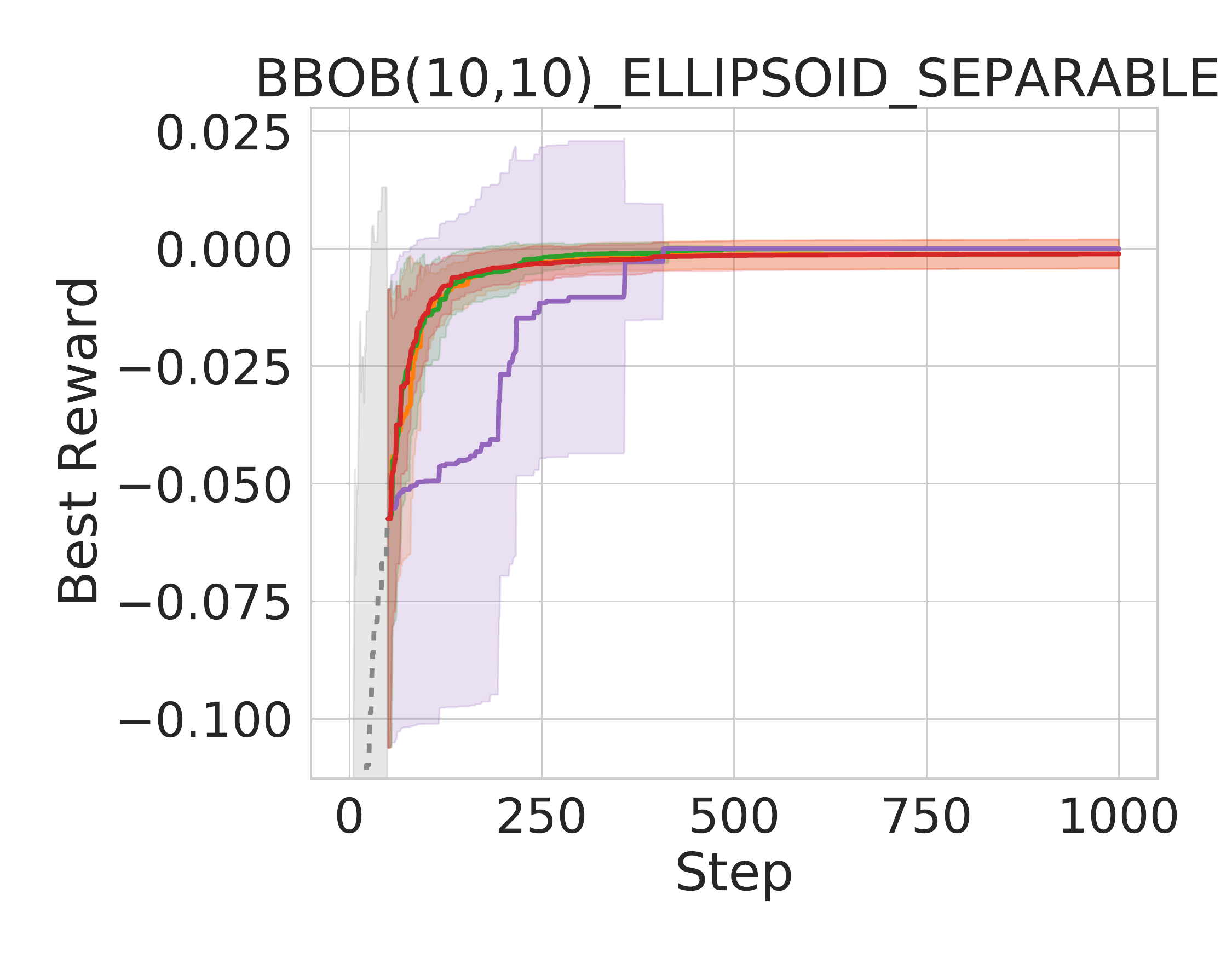}
    \includegraphics[width=\linewidth]{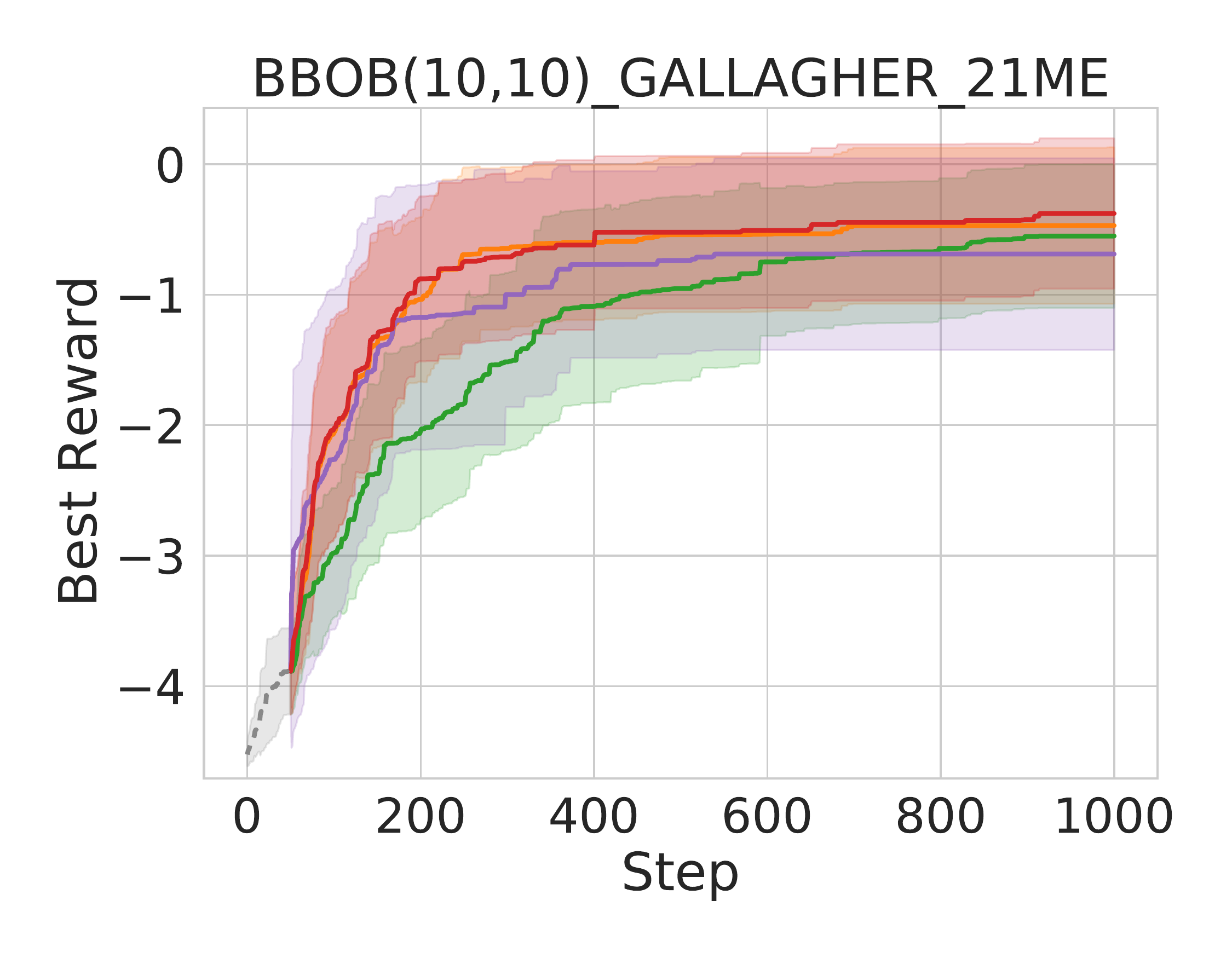}
    \includegraphics[width=\linewidth]{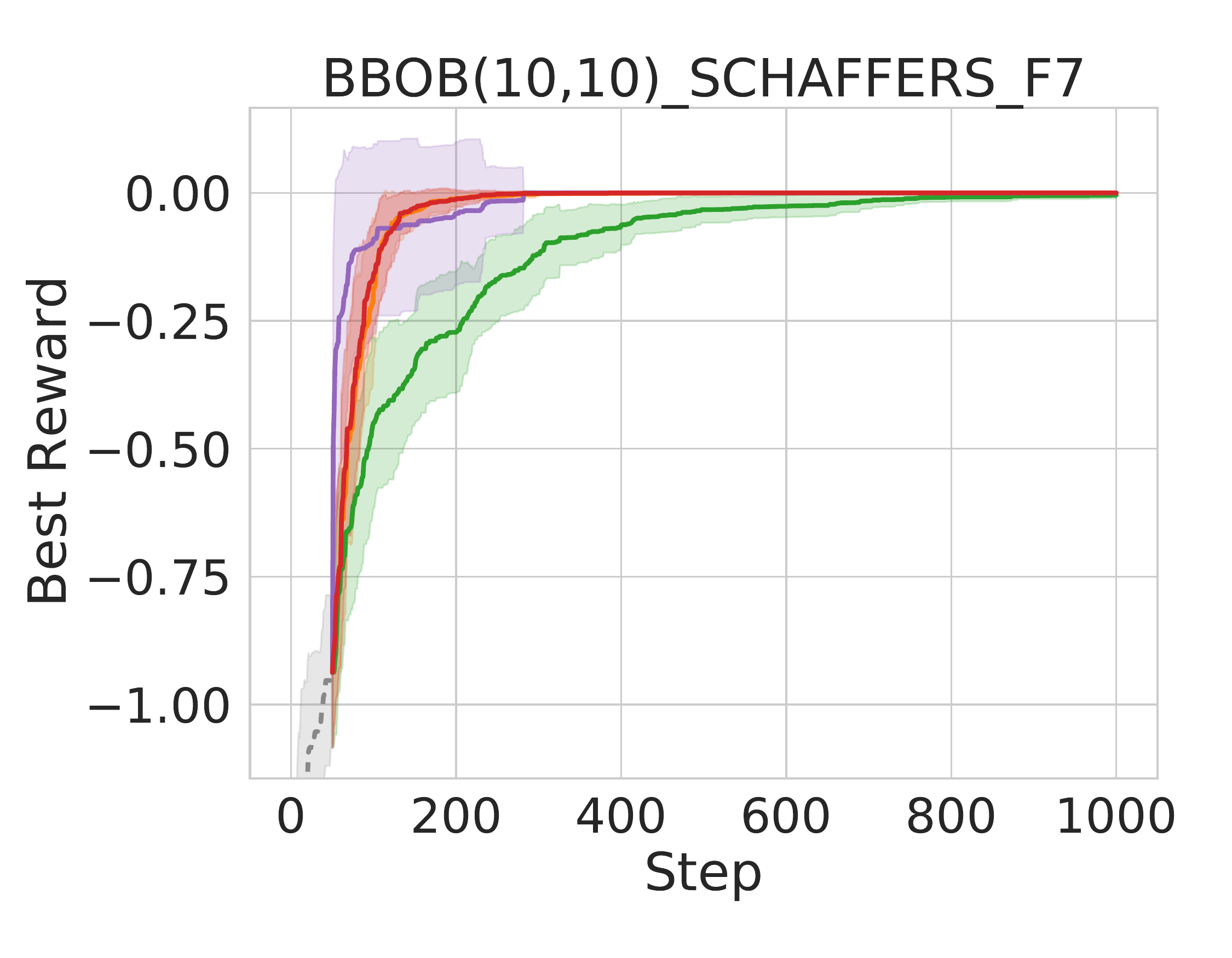}
    \end{multicols}
    \vspace{-30pt}
    \begin{multicols}{3}
    \includegraphics[width=\linewidth]{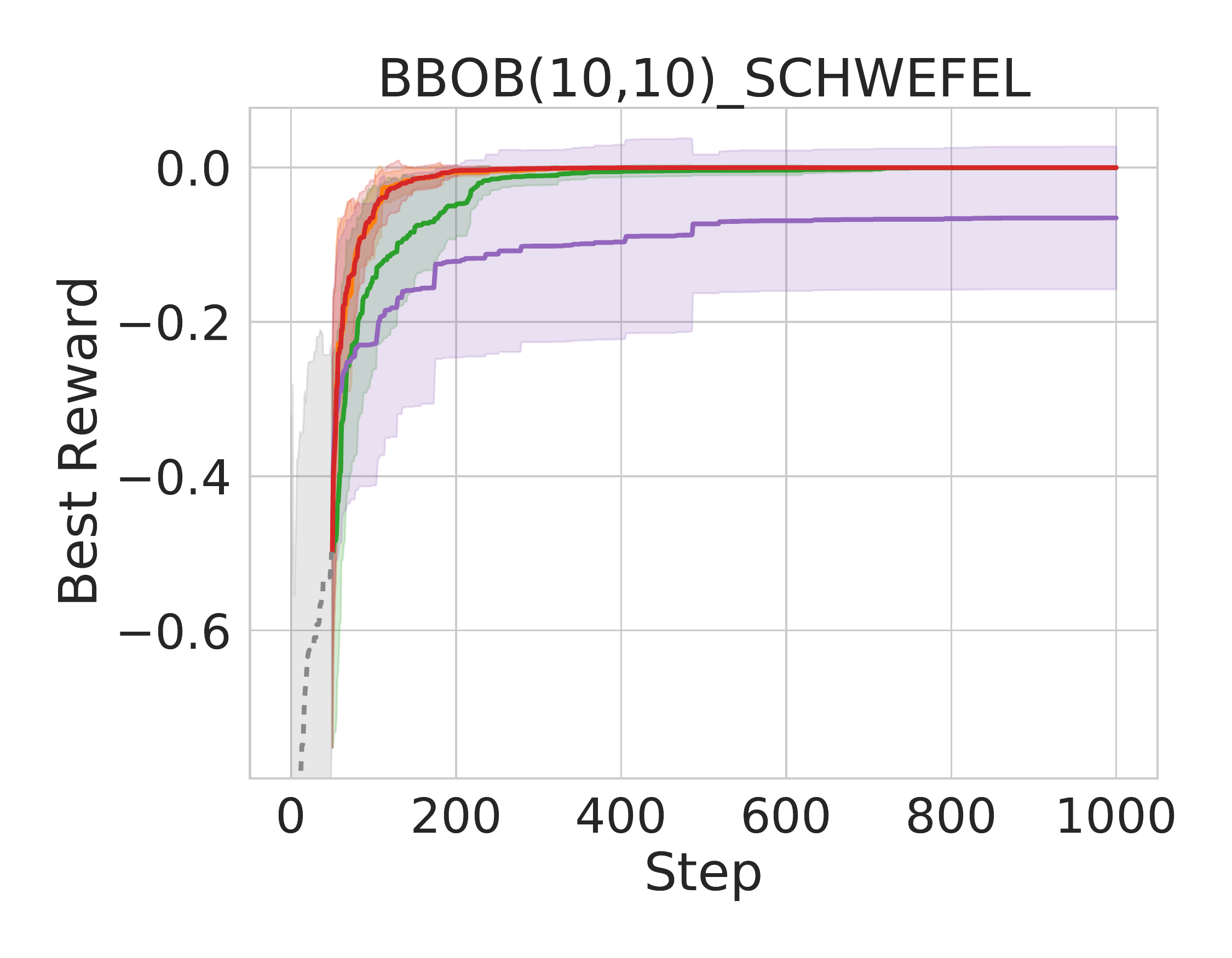}
    \includegraphics[width=\linewidth]{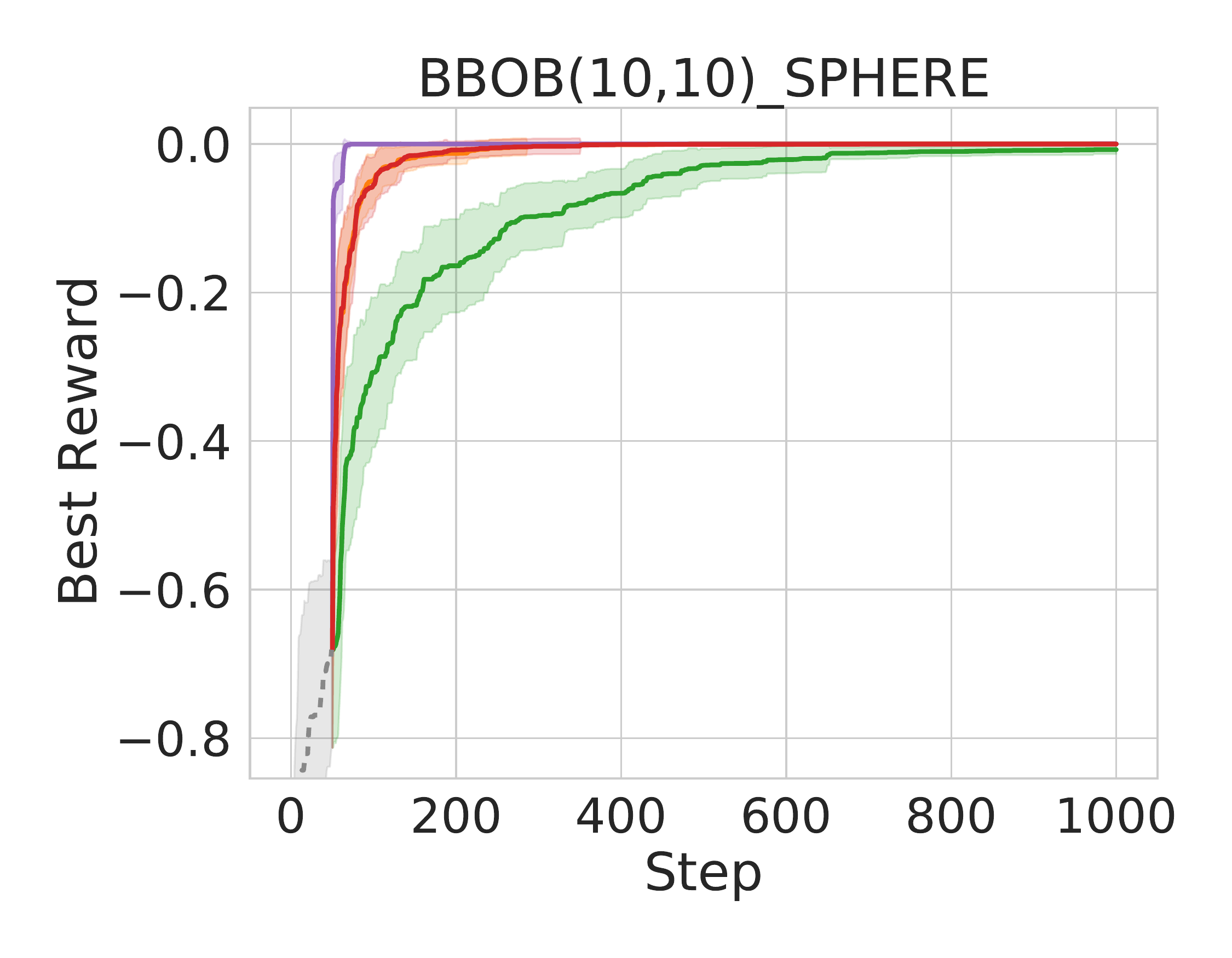}
    \includegraphics[width=\linewidth]{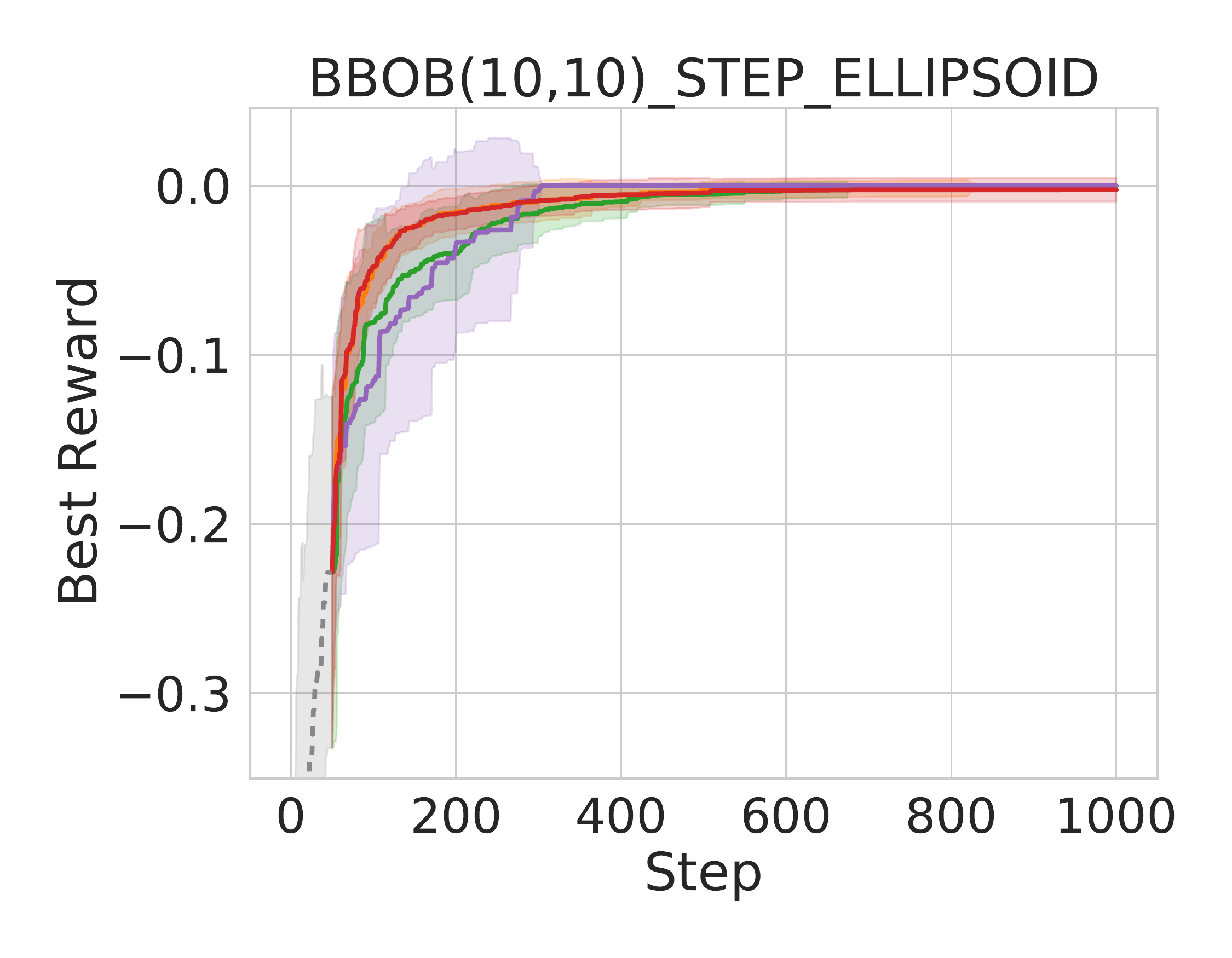}
    \end{multicols}
    \vspace{-30pt}
    \begin{multicols}{3}
    \includegraphics[width=\linewidth]{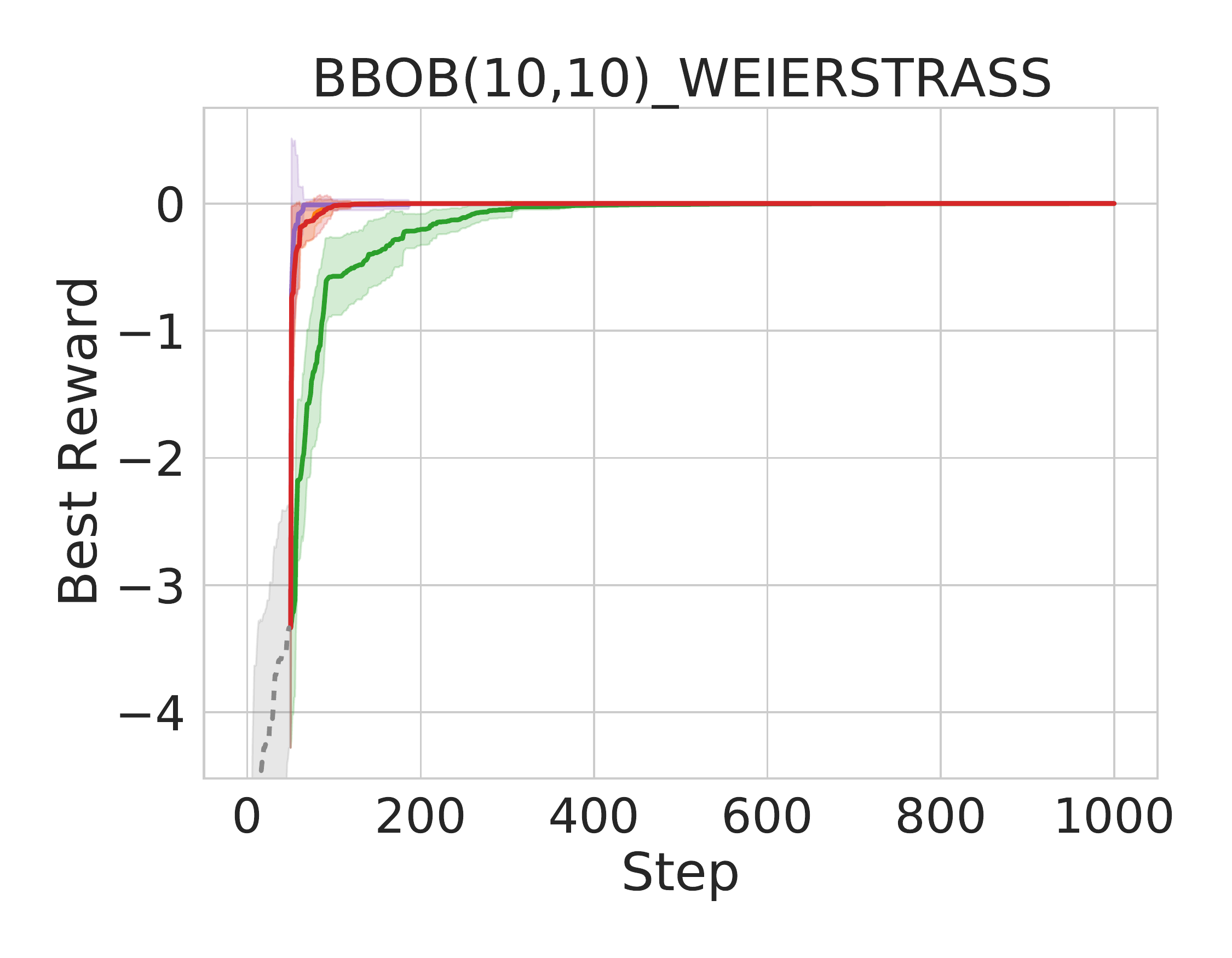}
    \includegraphics[width=\linewidth]{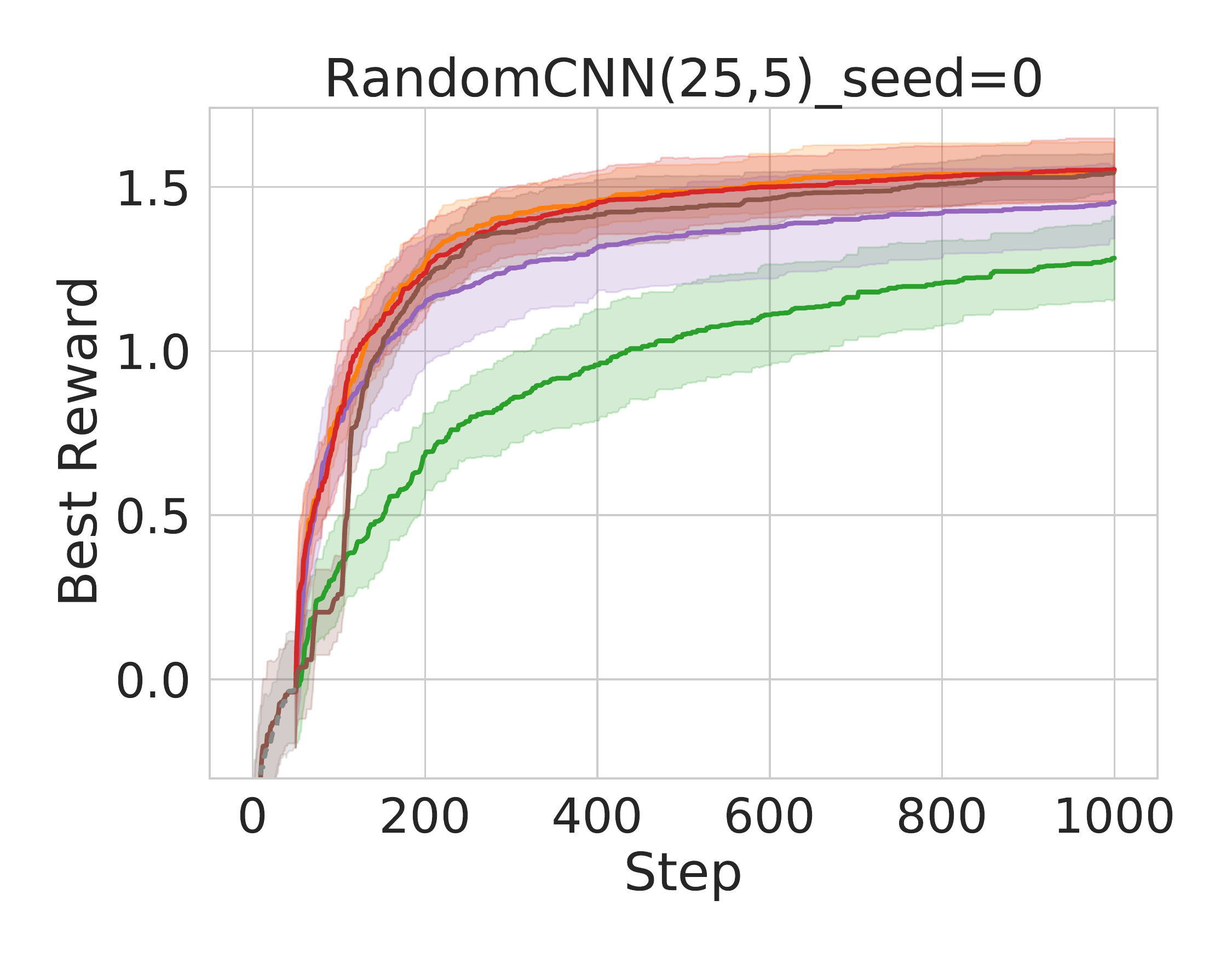}
    \includegraphics[width=\linewidth]{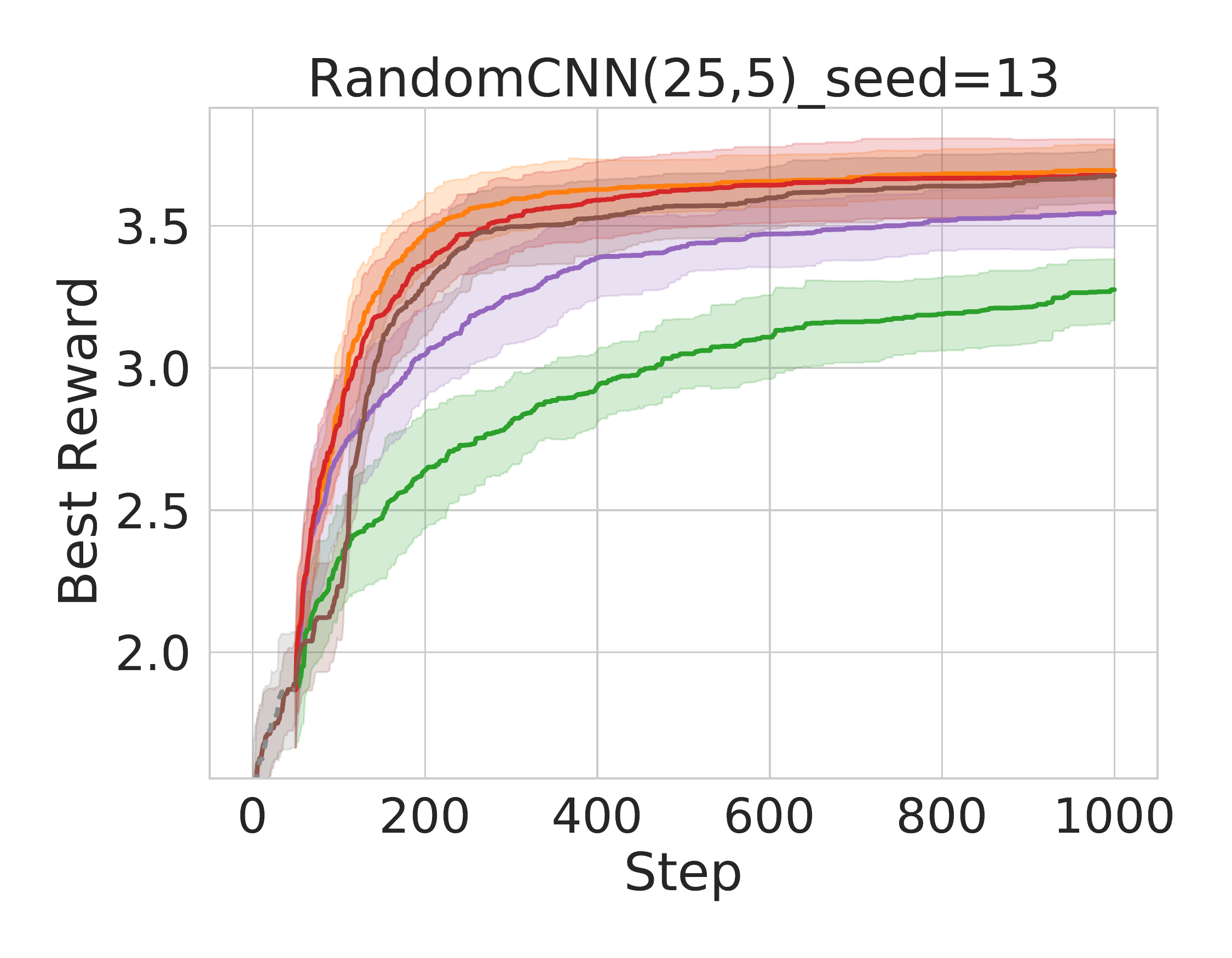}
    \end{multicols}
    \vspace{-30pt}
    \begin{multicols}{3}
    \includegraphics[width=\linewidth]{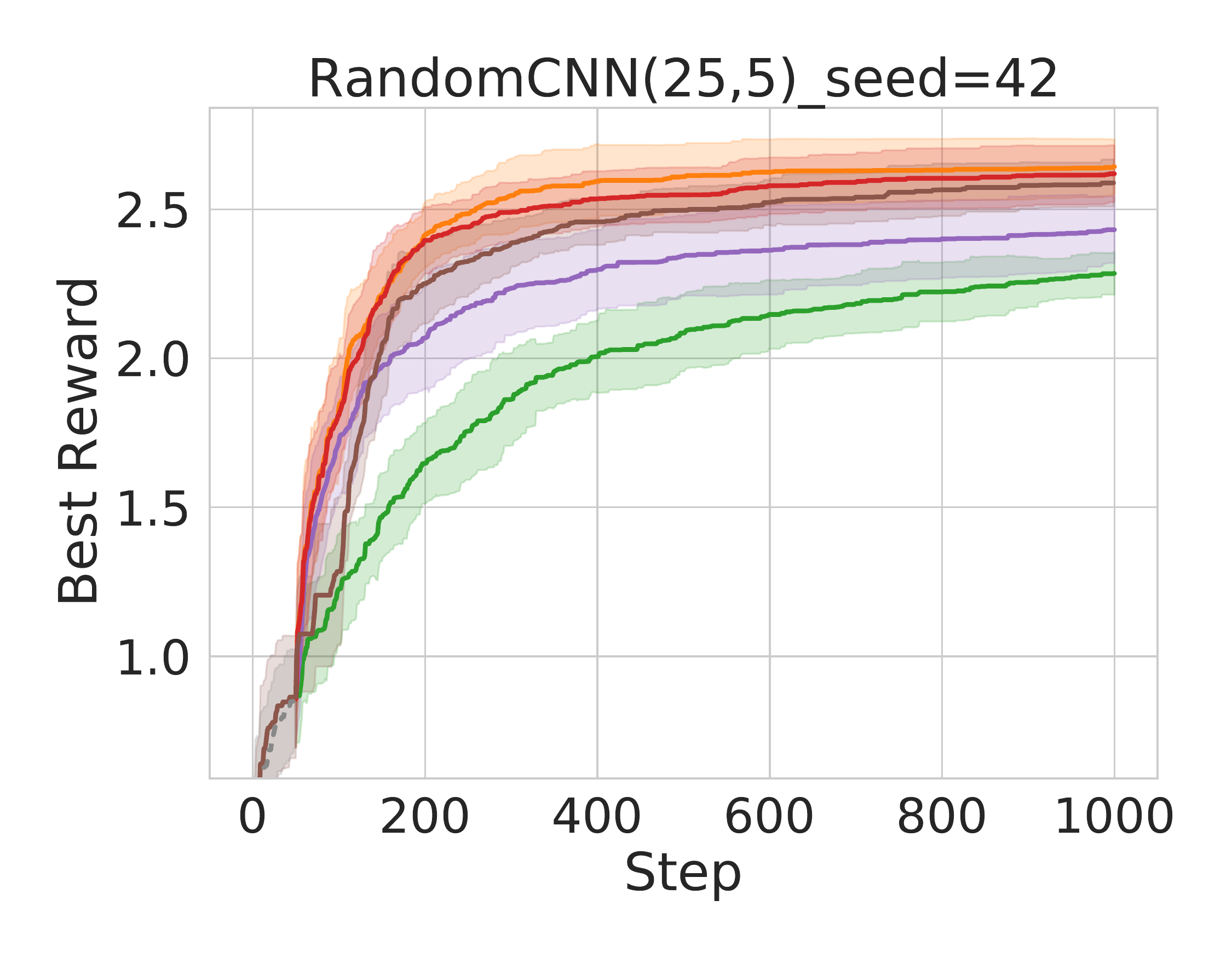}
    \includegraphics[width=\linewidth]{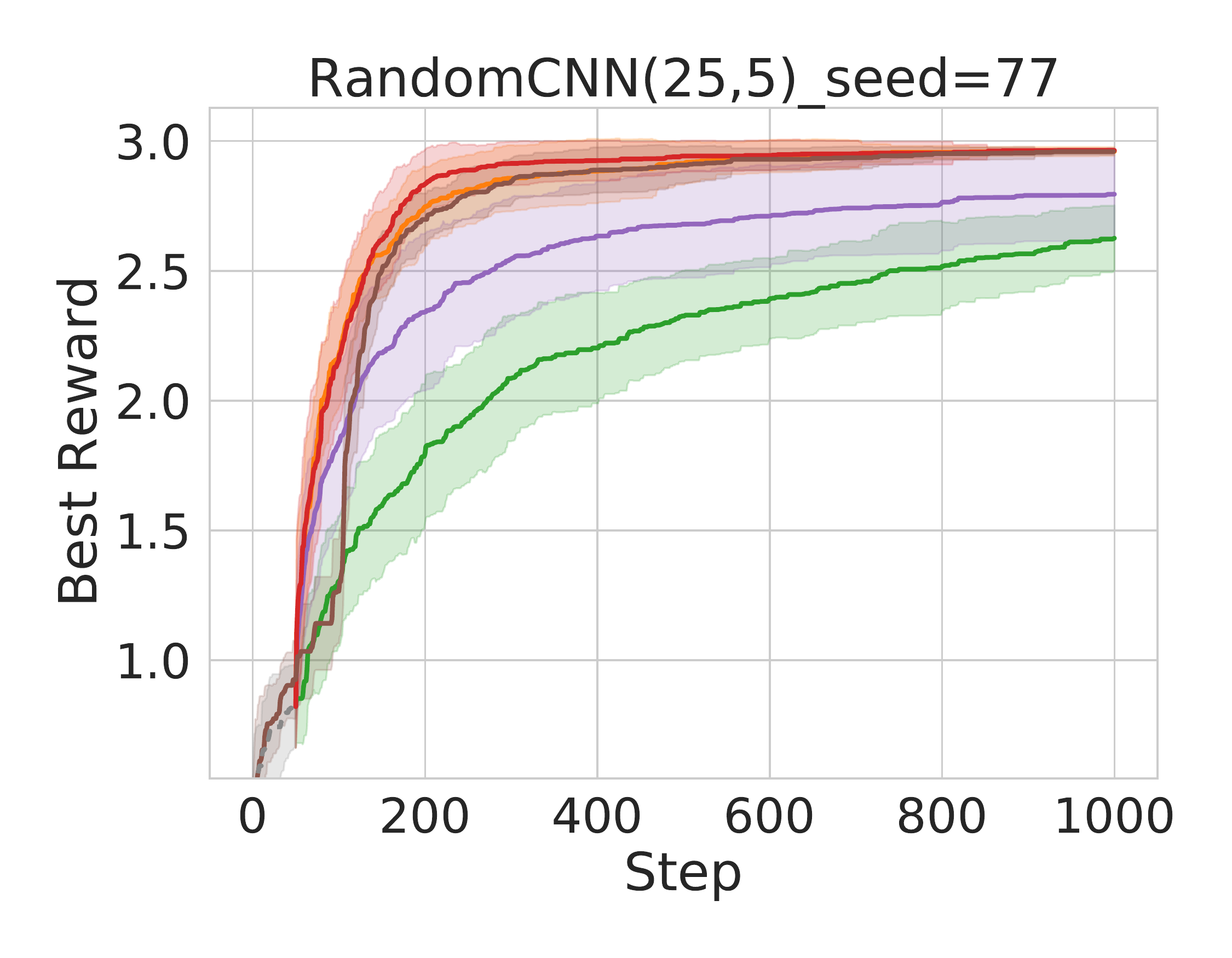}
    \includegraphics[width=\linewidth]{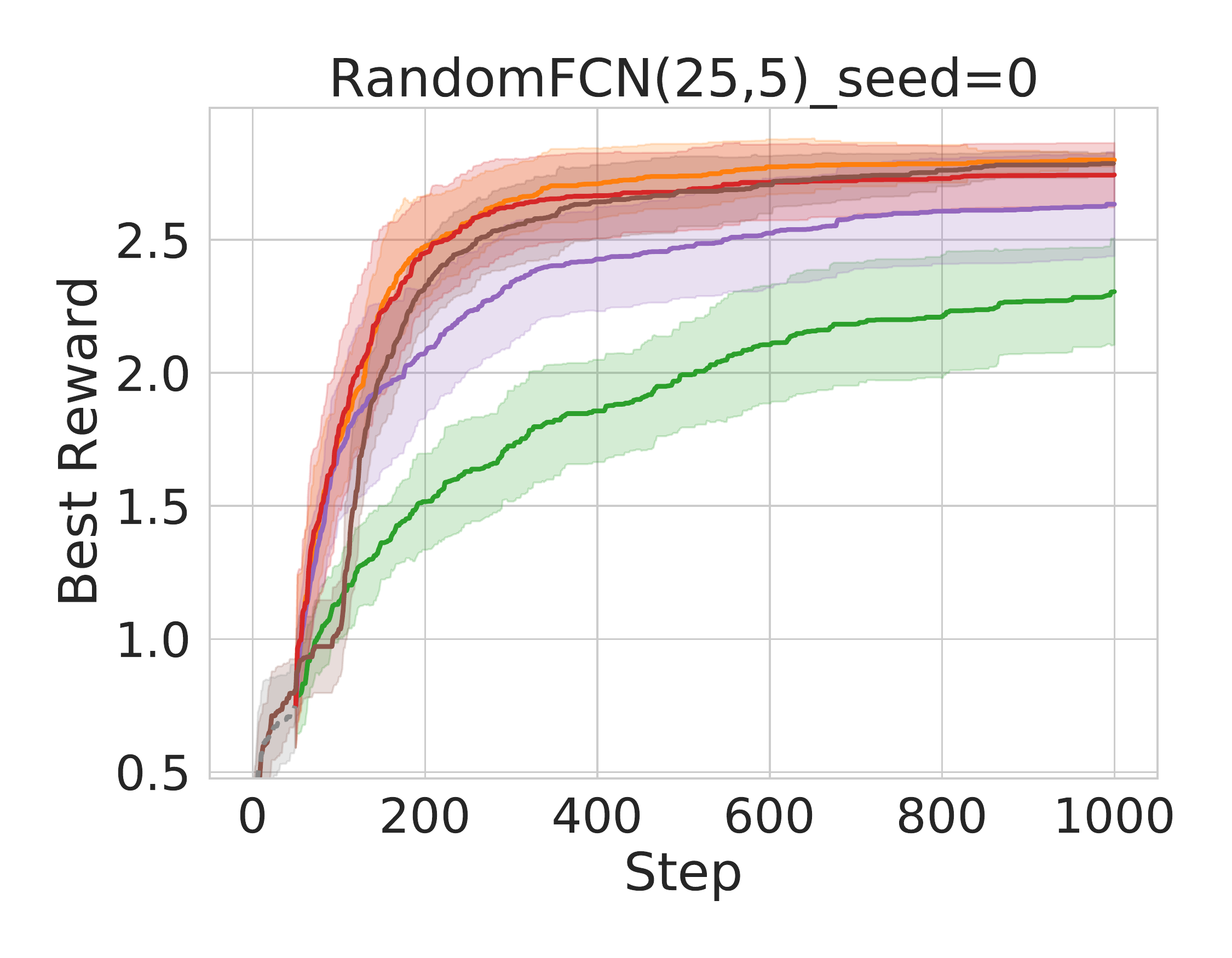}
    \end{multicols}
    \vspace{-20pt}
\centerline{\includegraphics[width=0.9\textwidth]{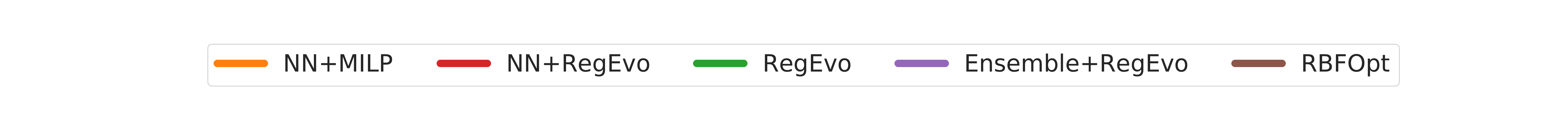}}
\vspace{-15pt}
\caption{Best observed reward as a function of iteration for the first half of all unconstrained problems (Section~\ref{sec:exp_unconstr}), averaged over 20 trials (bands indicate $\pm 1$sd). Dashed grey lines in the first 50 steps indicate the initial randomly sampled dataset, common to all methods except RBFOpt, which performs its own initialization.}
\label{fig:appendix_unconstrained_curves1}

\end{center}
\end{figure*}

\begin{figure*}[ht]
\begin{center}
    \begin{multicols}{3}
    \includegraphics[width=\linewidth]{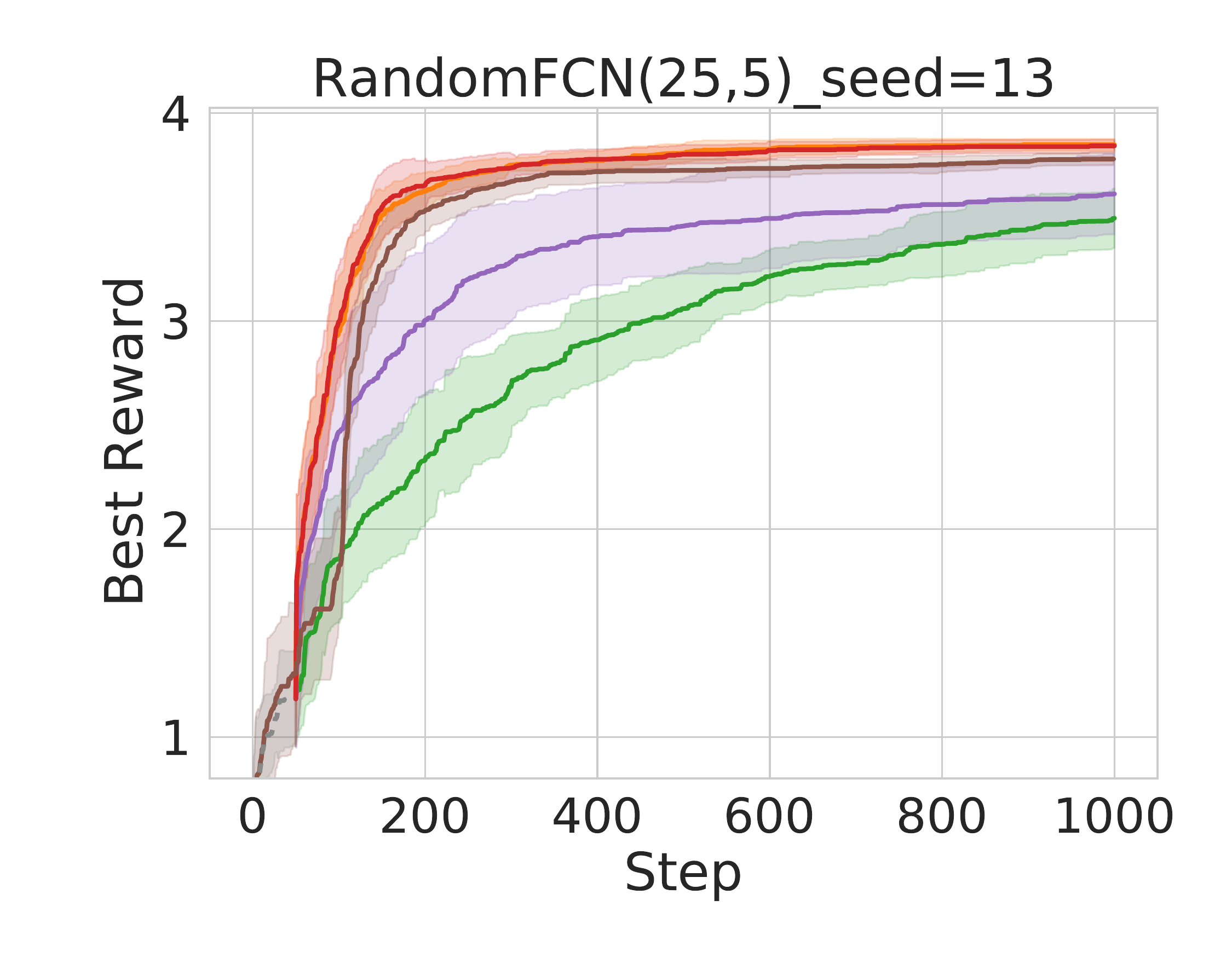}
    \includegraphics[width=\linewidth]{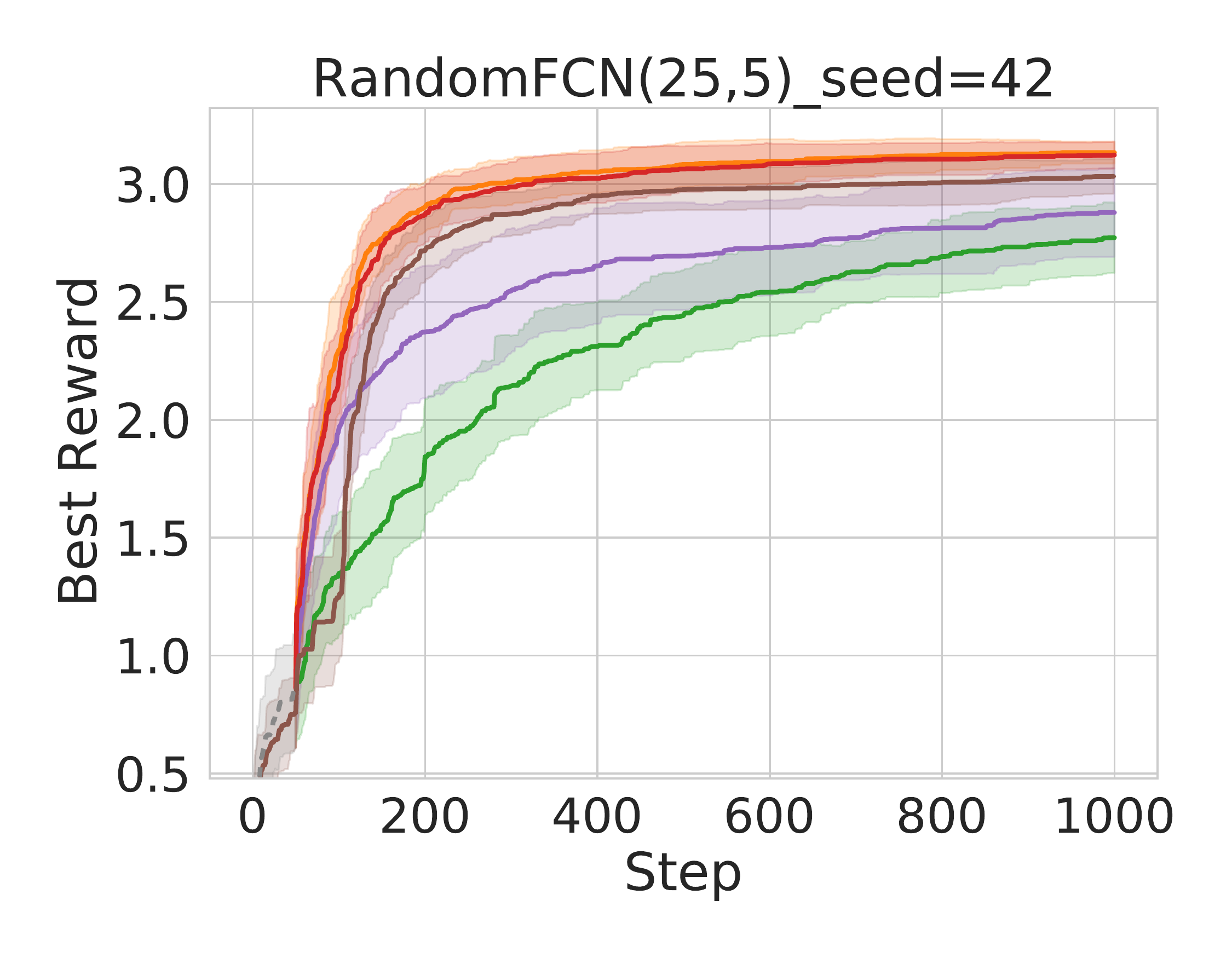}
    \includegraphics[width=\linewidth]{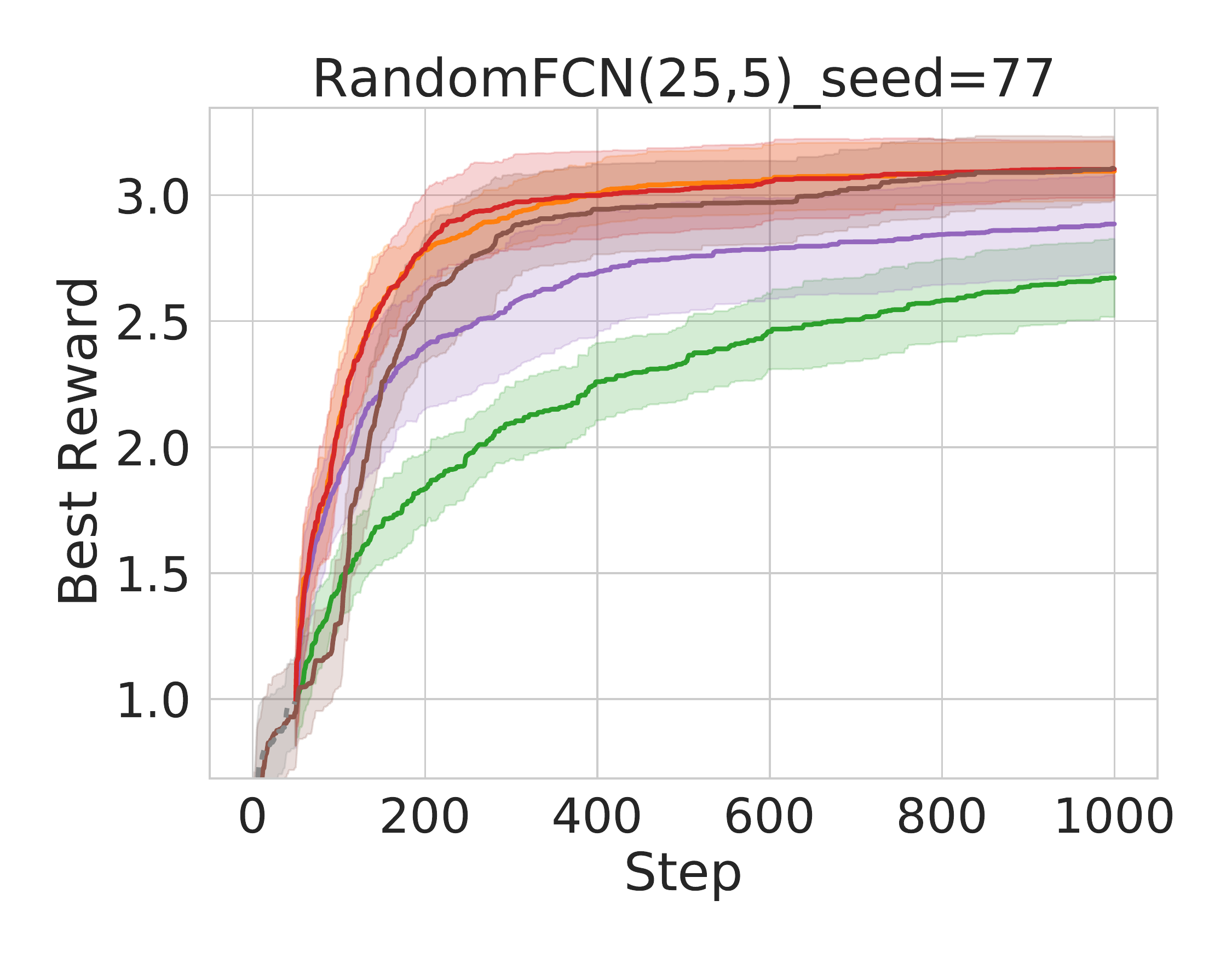}
    \end{multicols}
    \vspace{-30pt}
    \begin{multicols}{3}
    \includegraphics[width=\linewidth]{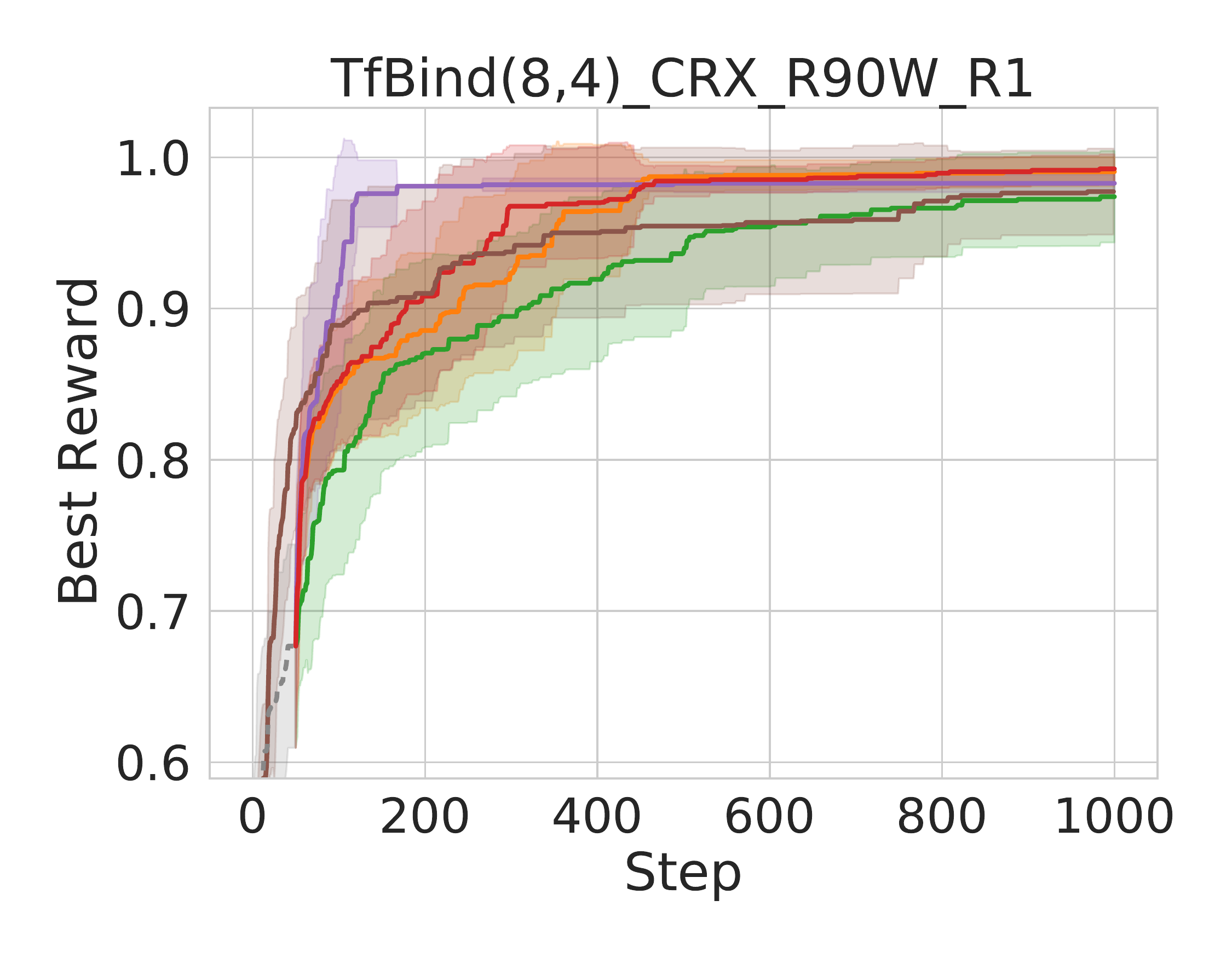}
    \includegraphics[width=\linewidth]{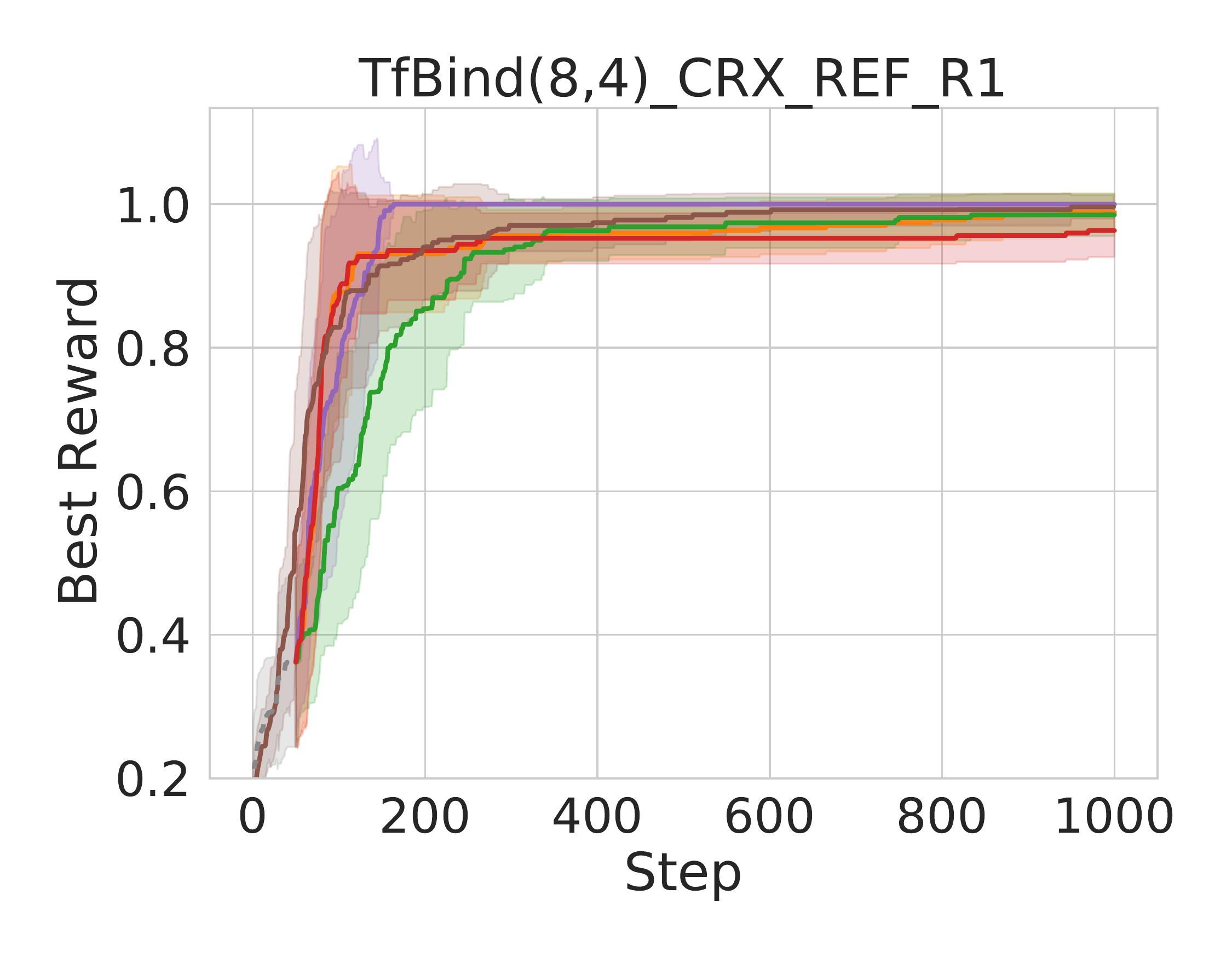}
    \includegraphics[width=\linewidth]{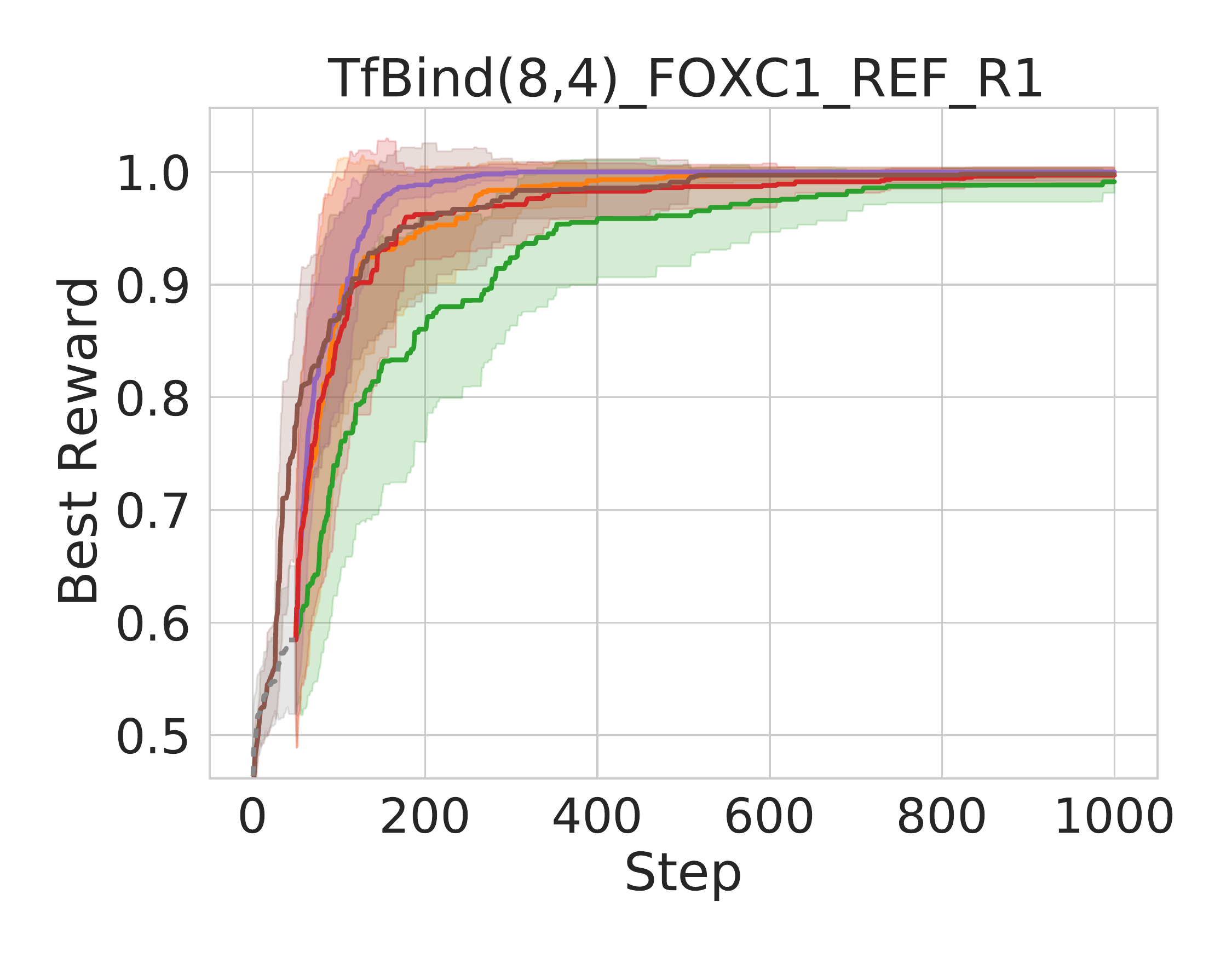}
    \end{multicols}
    \vspace{-30pt}
    \begin{multicols}{3}
    \includegraphics[width=\linewidth]{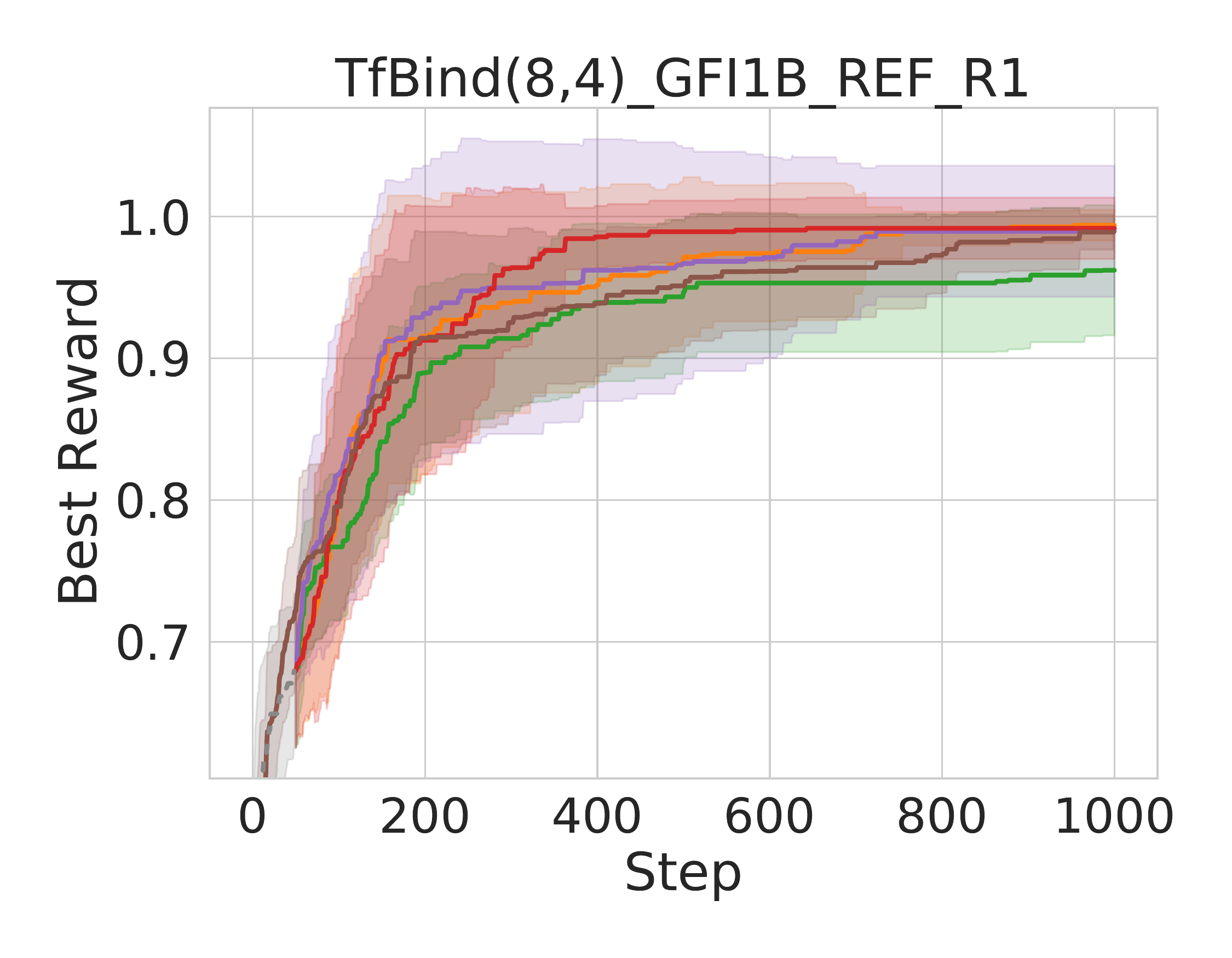}
    \includegraphics[width=\linewidth]{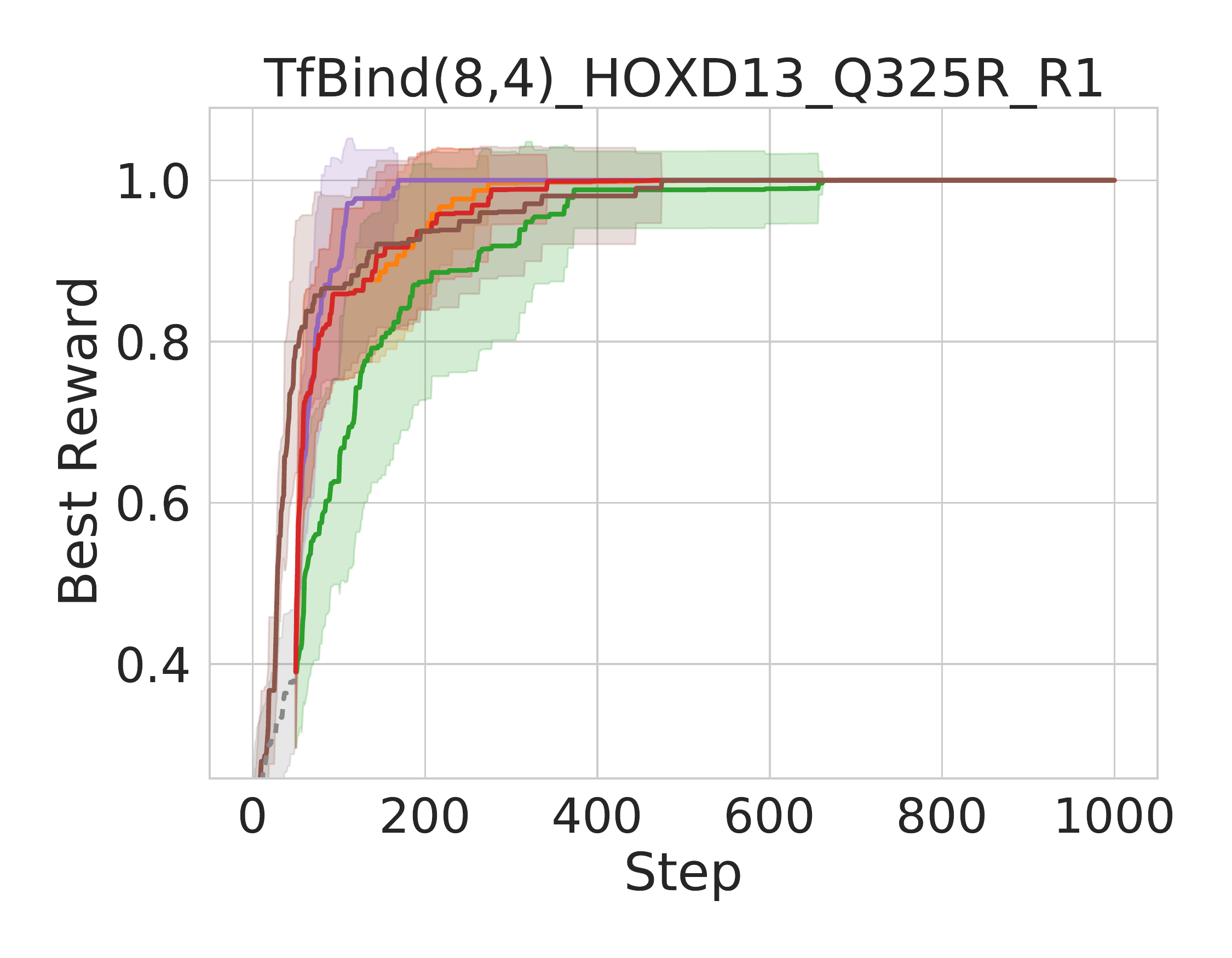}
    \includegraphics[width=\linewidth]{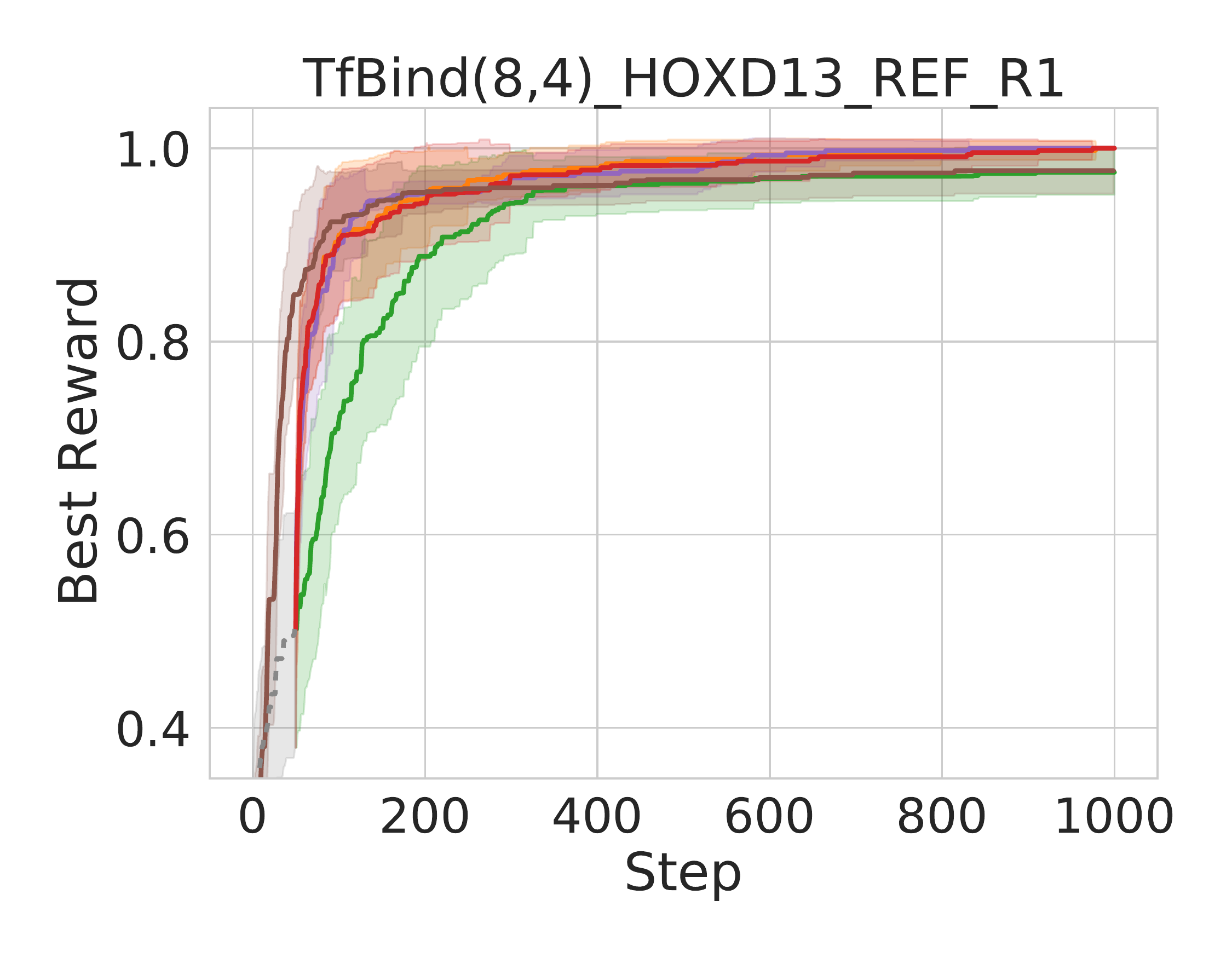}
    \end{multicols}
    \vspace{-30pt}
    \begin{multicols}{3}
    \includegraphics[width=\linewidth]{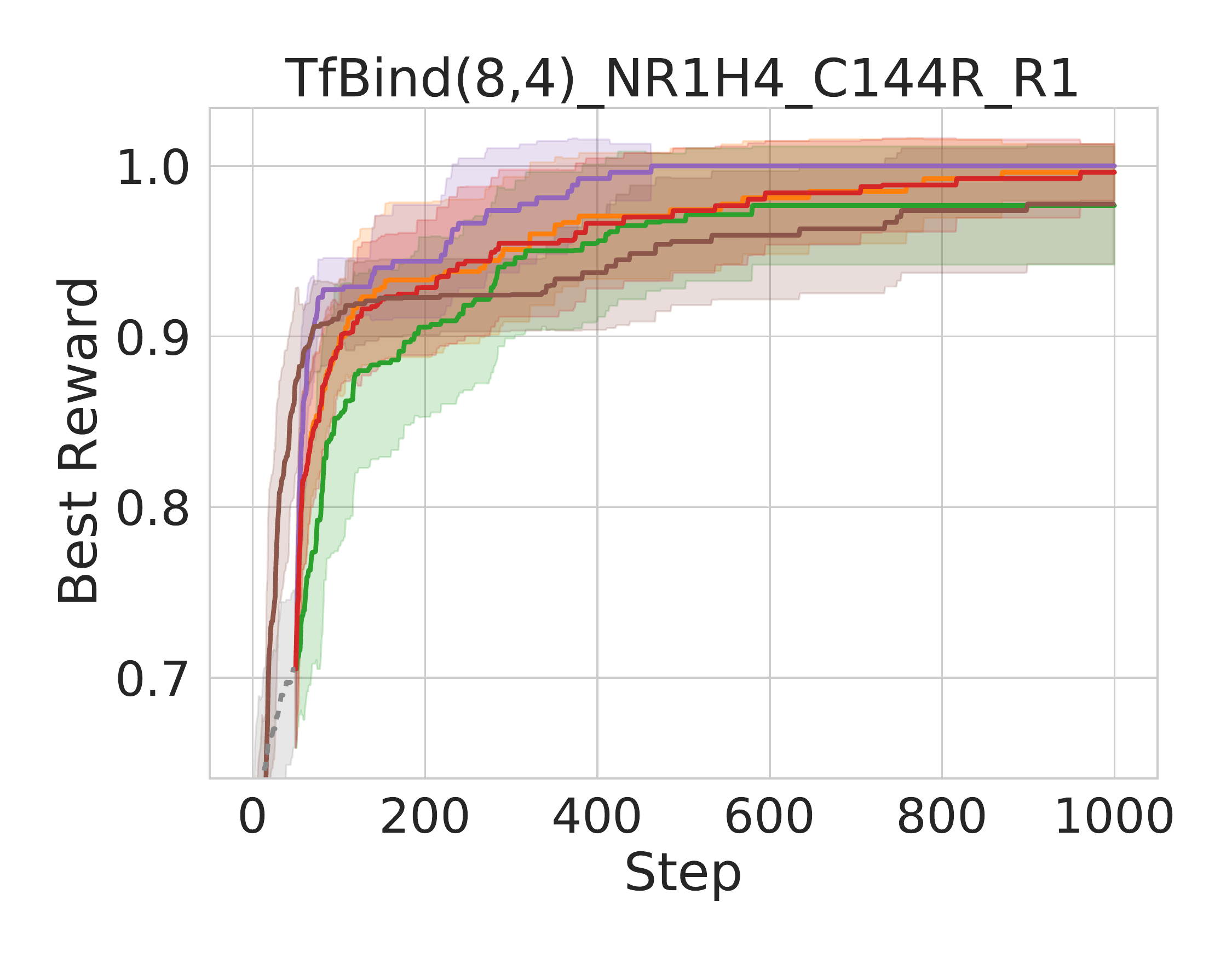}
    \includegraphics[width=\linewidth]{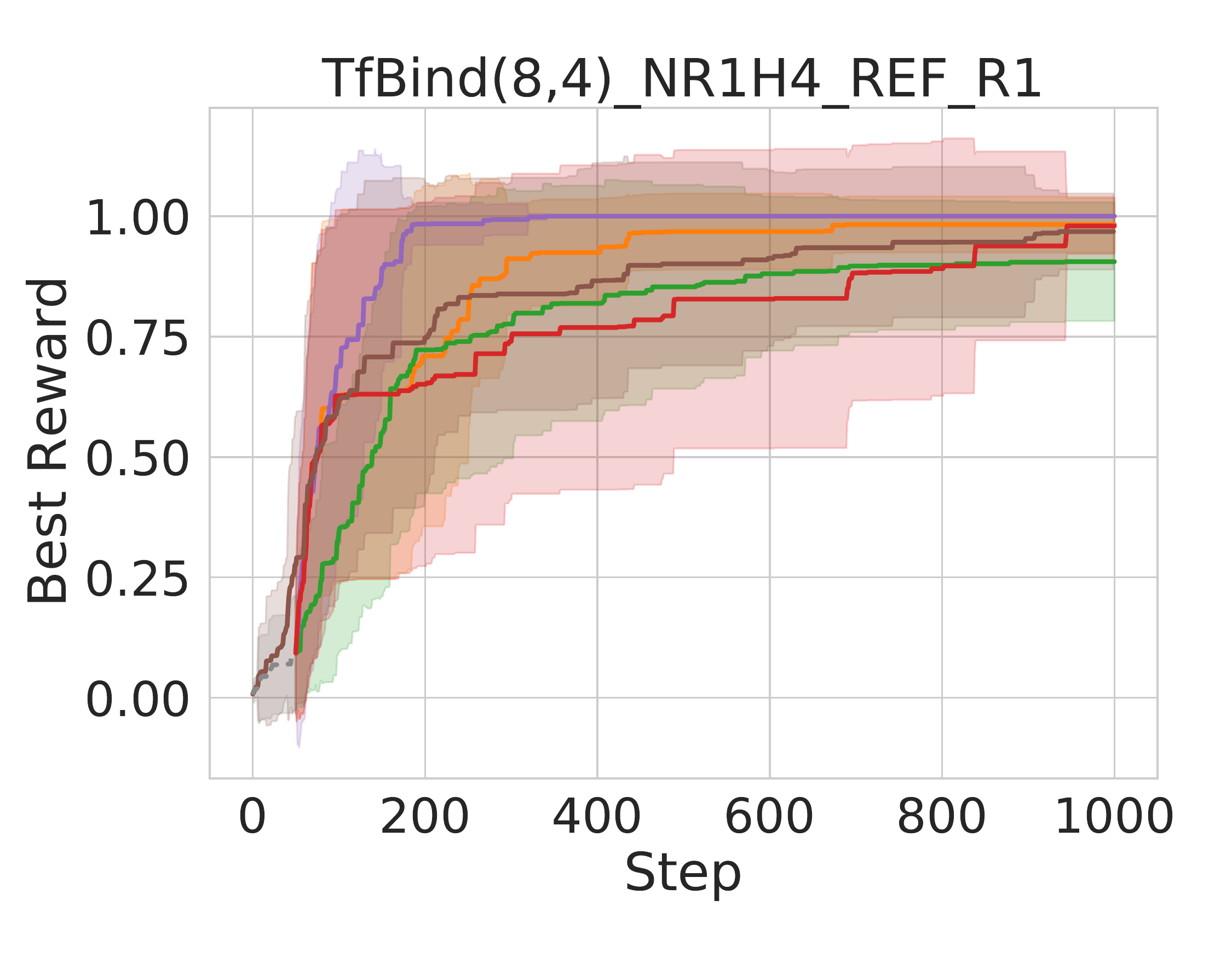}
    \includegraphics[width=\linewidth]{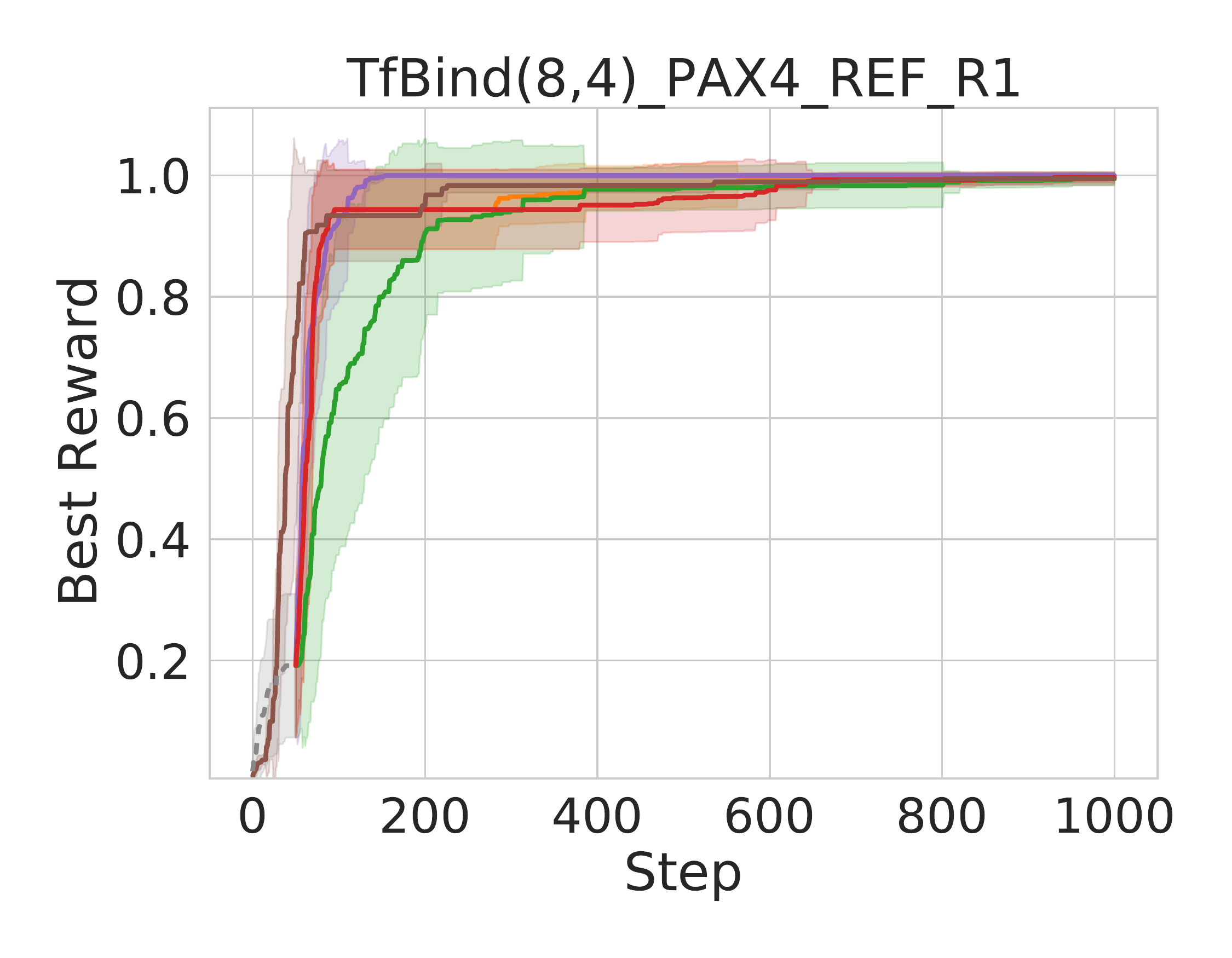}
    \end{multicols}
    \vspace{-30pt}
    \begin{multicols}{3}
    \includegraphics[width=\linewidth]{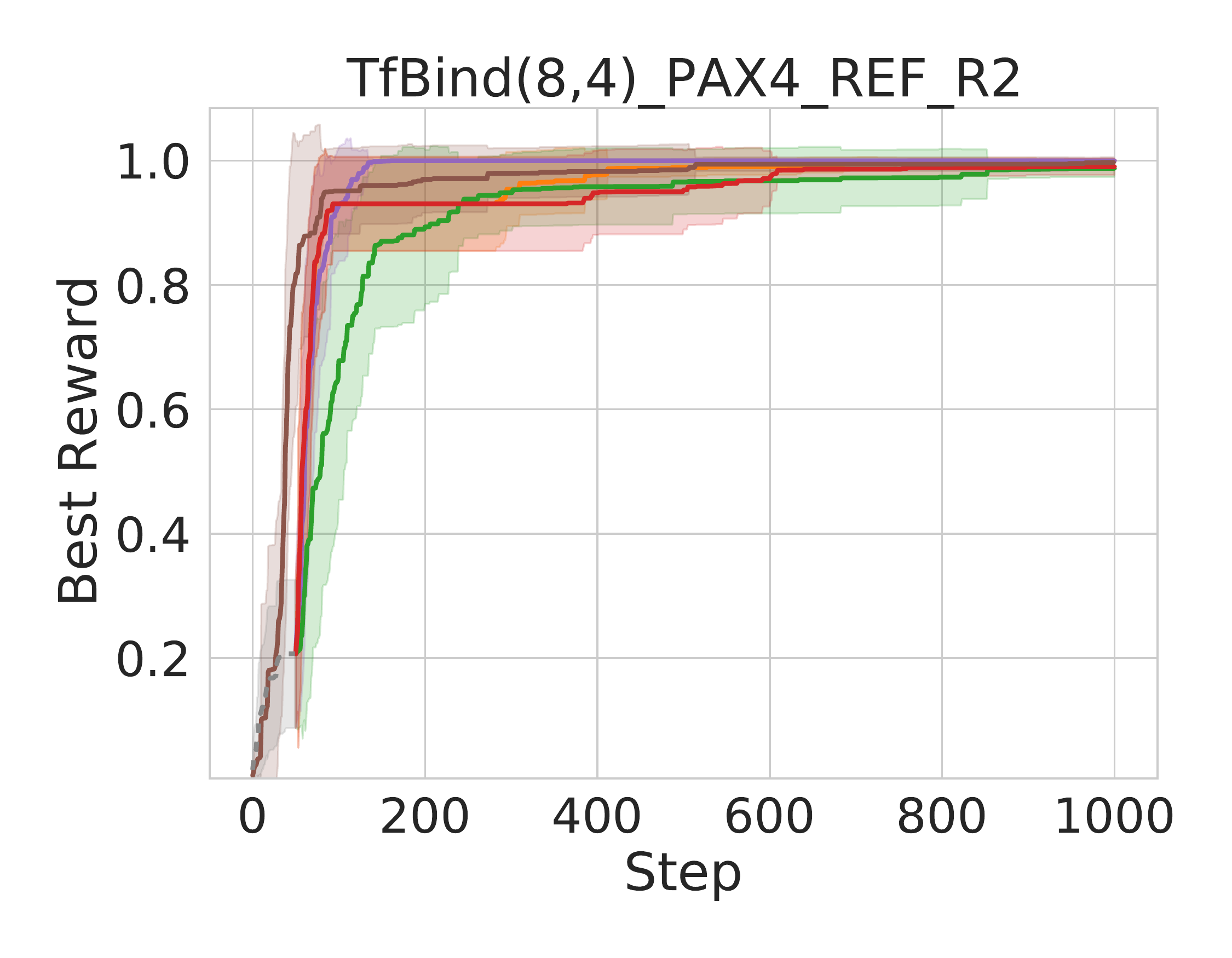}
    \includegraphics[width=\linewidth]{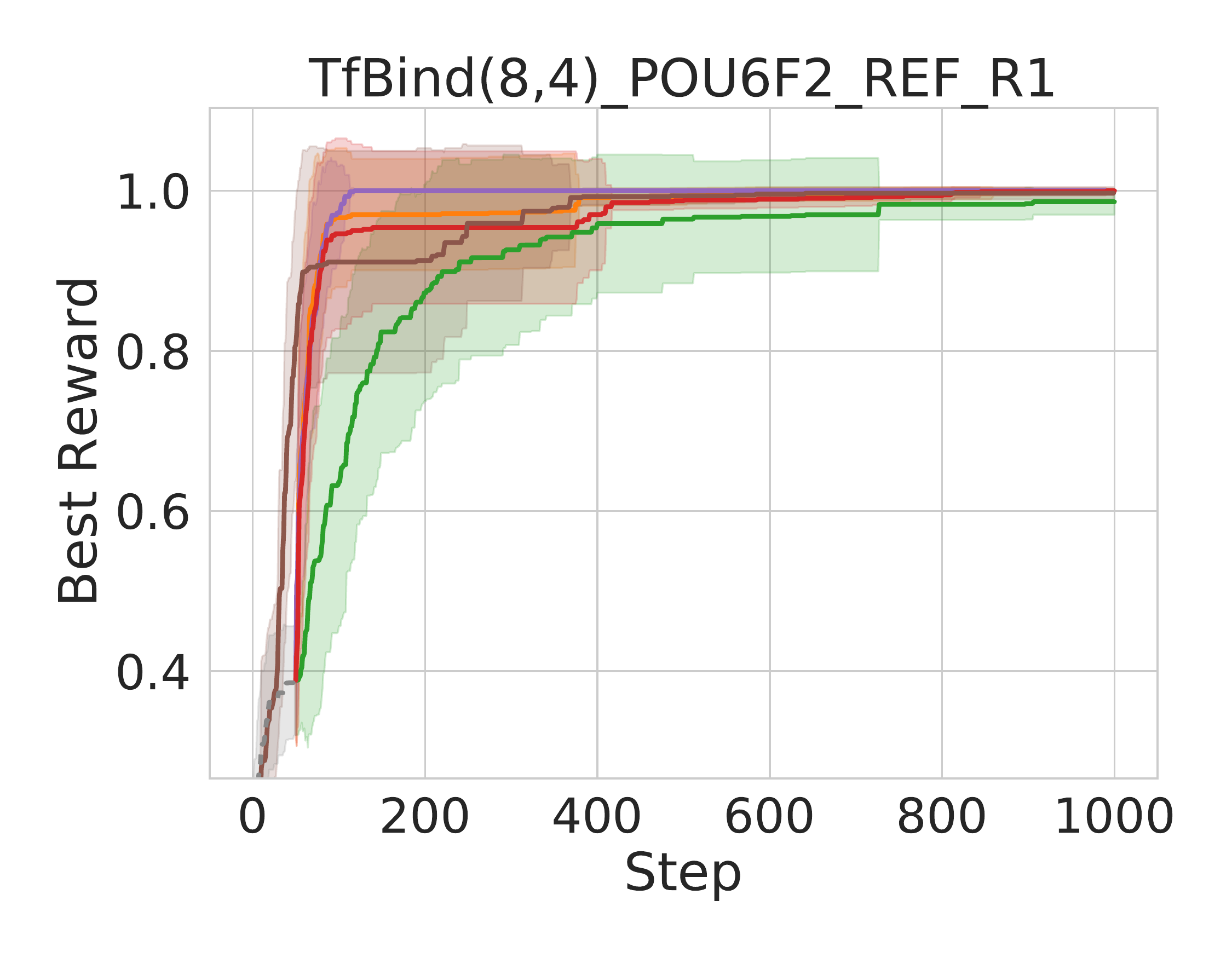}
    \includegraphics[width=\linewidth]{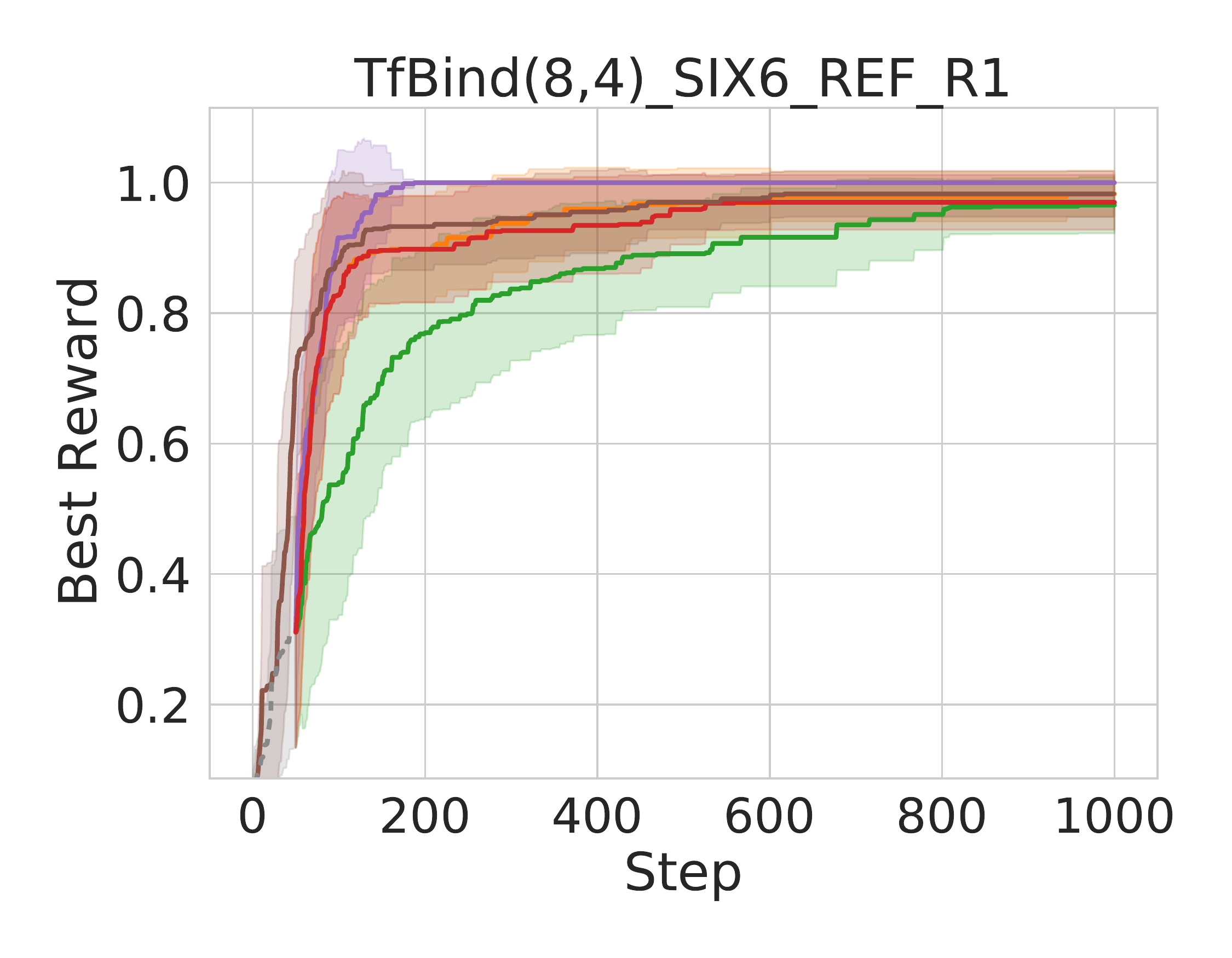}
    \end{multicols}
    \vspace{-20pt}
\centerline{\includegraphics[width=0.9\textwidth]{a1_unconstr_legend_with_rbfopt.pdf}}
\vspace{-15pt}
\caption{Best observed reward as a function of iteration for the second half of all unconstrained problems (Section~\ref{sec:exp_unconstr}), averaged over 20 trials (bands indicate $\pm 1$sd). Dashed grey lines in the first 50 steps indicate the initial randomly sampled dataset, common to all methods except RBFOpt, which performs its own initialization.}
\label{fig:appendix_unconstrained_curves2}
\end{center}
\end{figure*}

Figure~\ref{fig:appendix_unconstrained_scores_auc} plots the distribution of algorithms' normalized AUC scores over all unconstrained problems, split by objective function class. The relative performance of algorithms in terms of this new AUC metric does not differ significantly from what we found for final reward (Section~\ref{sec:exp_unconstr}, Figure~\ref{fig:unconstrained_scores}). Figures~\ref{fig:appendix_unconstrained_curves1} and~\ref{fig:appendix_unconstrained_curves2} plot the individual reward curves as function of outer-loop iteration for each unconstrained problem.

\begin{figure*}[htbp]
\begin{center}
\includegraphics[width=0.33\textwidth]{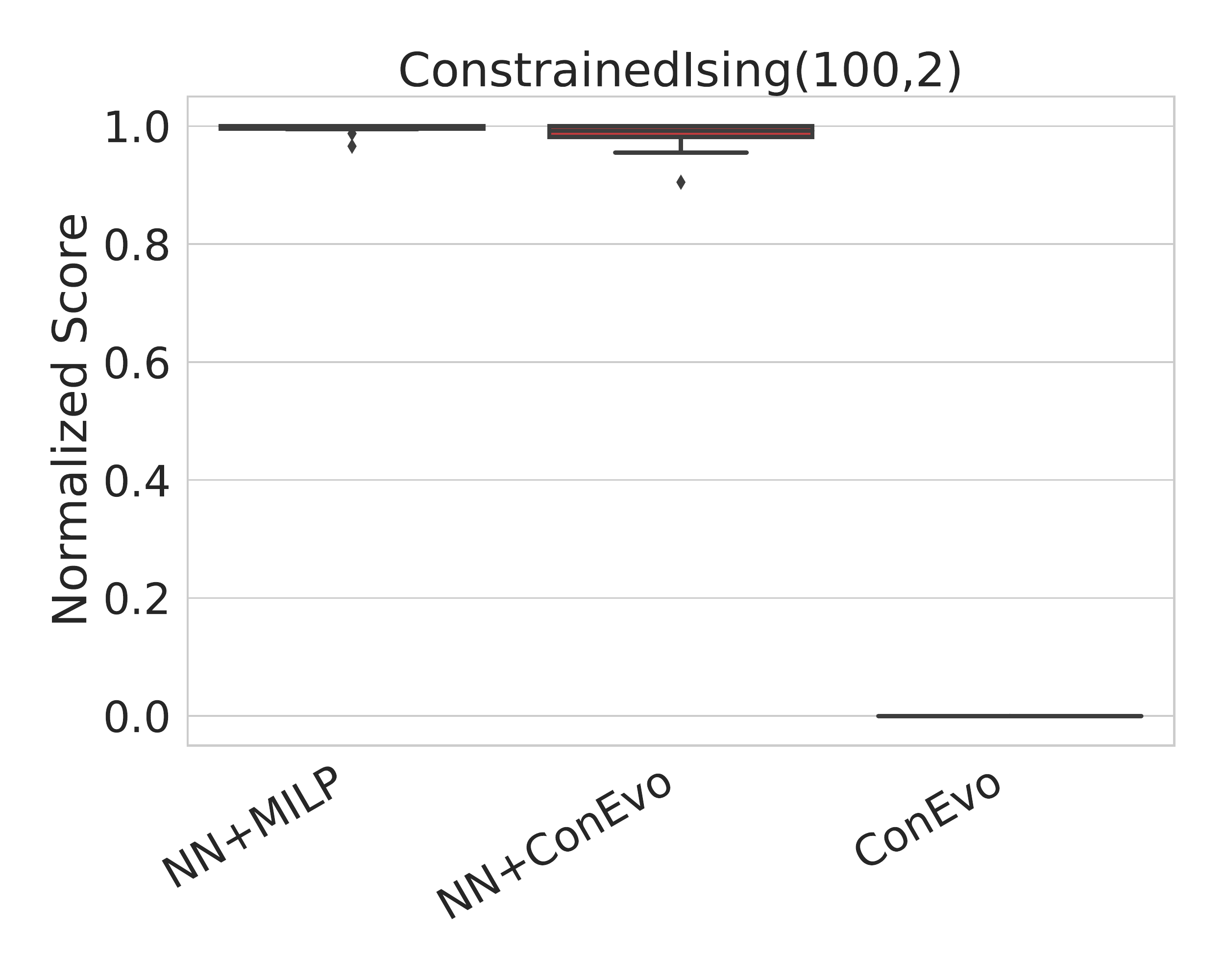}
\includegraphics[width=0.33\textwidth]{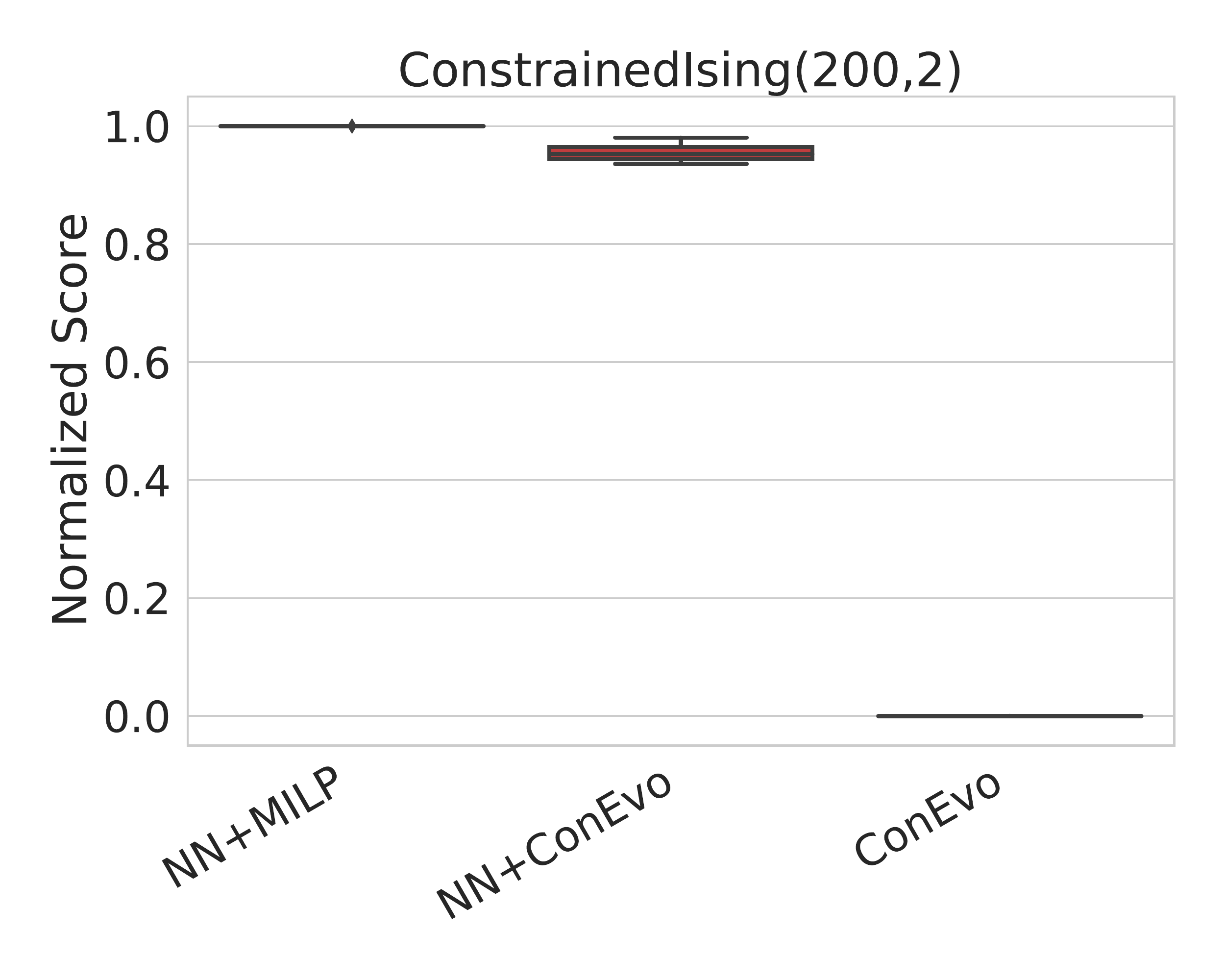}
\includegraphics[width=0.33\textwidth]{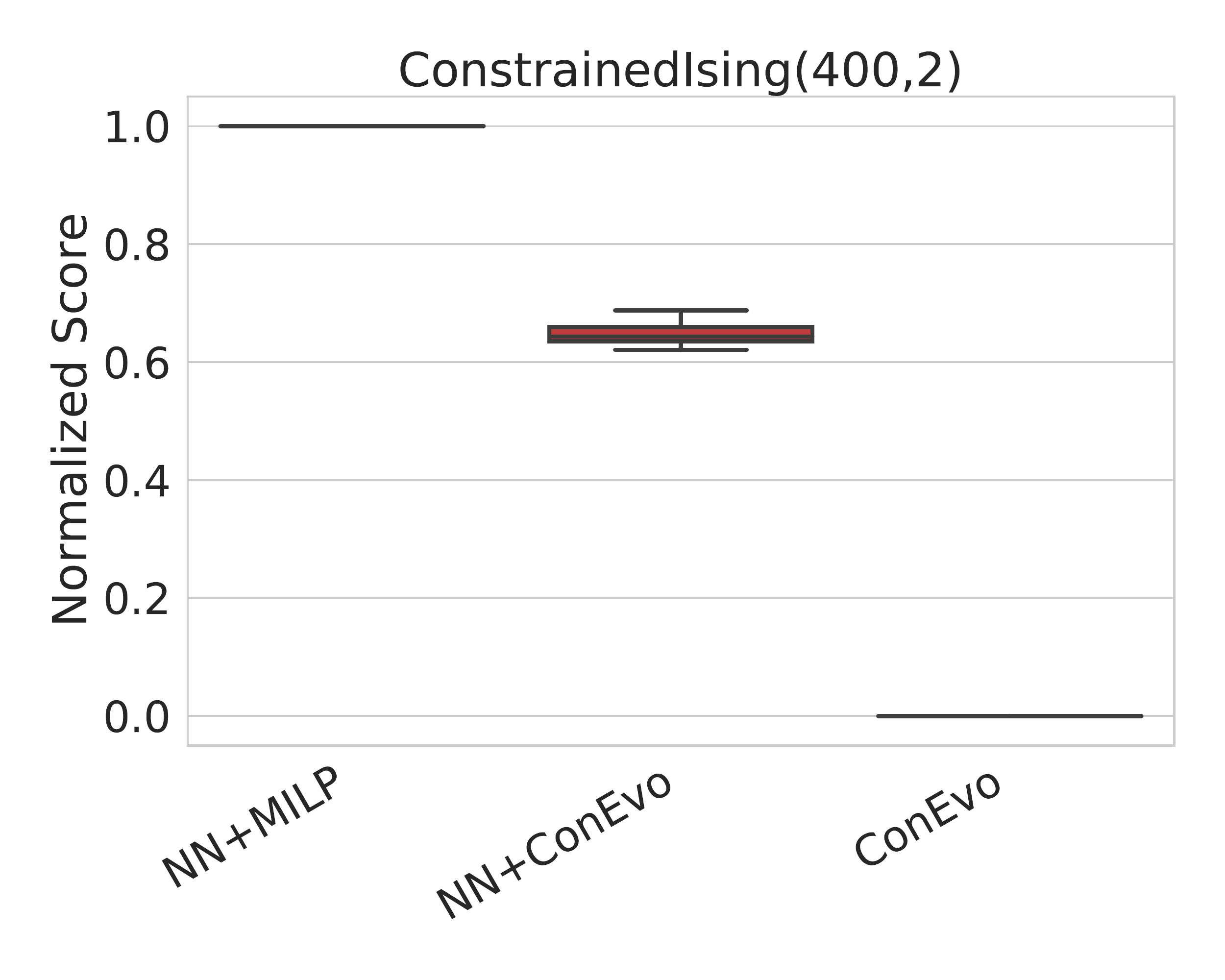}
\caption{Distribution of algorithms' normalized max-reward scores (defined Section~\ref{sec:exp_tasks}) on constrained Ising problems from Section~\ref{sec:exp_constr}, split by problem size (left-to-right: $n=100,200$ and $400$). Higher is better. Exact optimization of the acquisition function during MBBO (\textit{NN+MILP}) provides significant benefits compared to evolutionary heuristics (\textit{NN+ConEvo}) at large scales. Individual reward curves for each problem are given in Figures~\ref{fig:appendix_constrained_curves1}~and~\ref{fig:appendix_constrained_curves2}.}
\label{fig:appendix_constrained_scores}
\end{center}
\end{figure*}

\begin{figure*}[htbp]
    \includegraphics[width=0.48\textwidth]{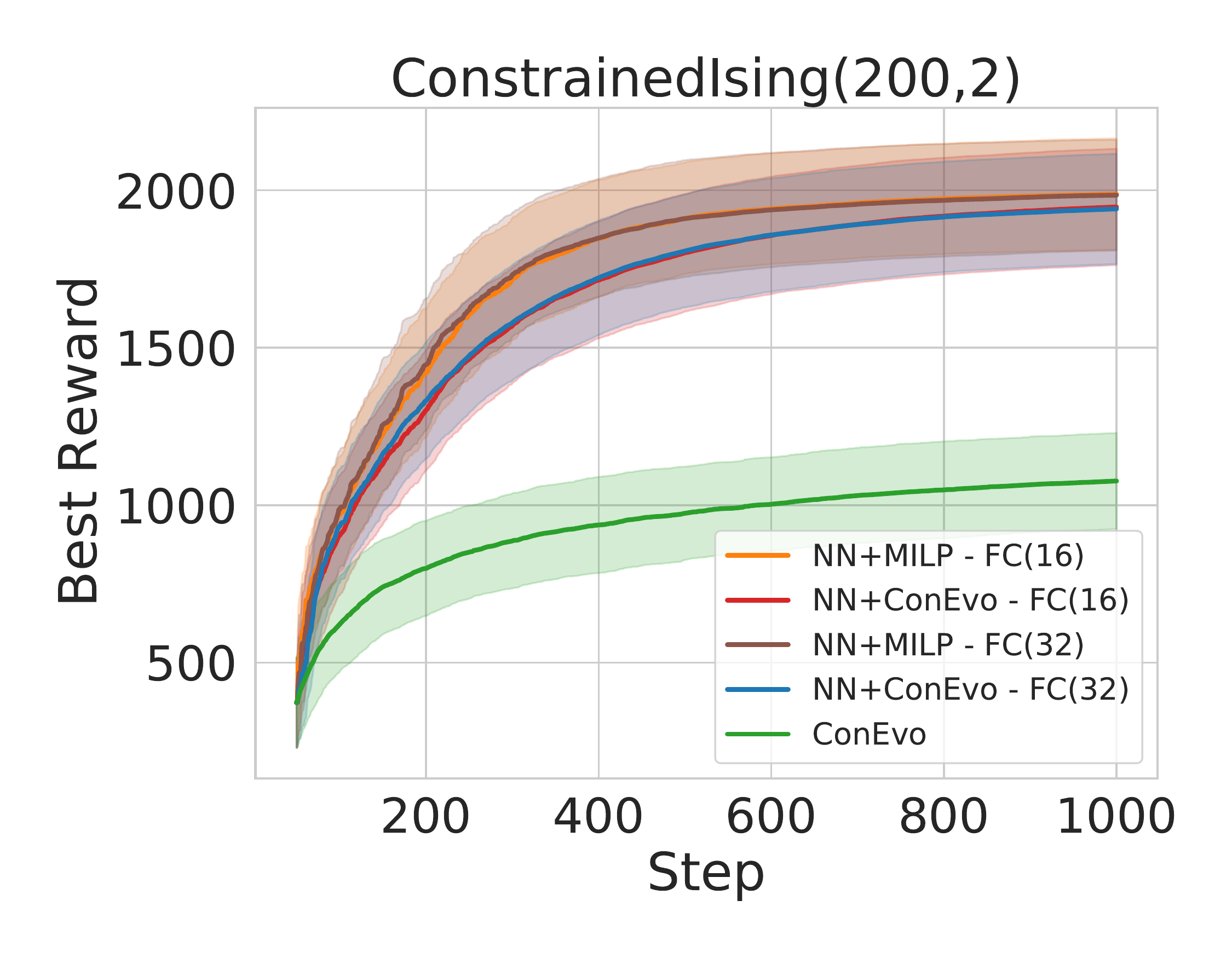}
    \includegraphics[width=0.48\textwidth]{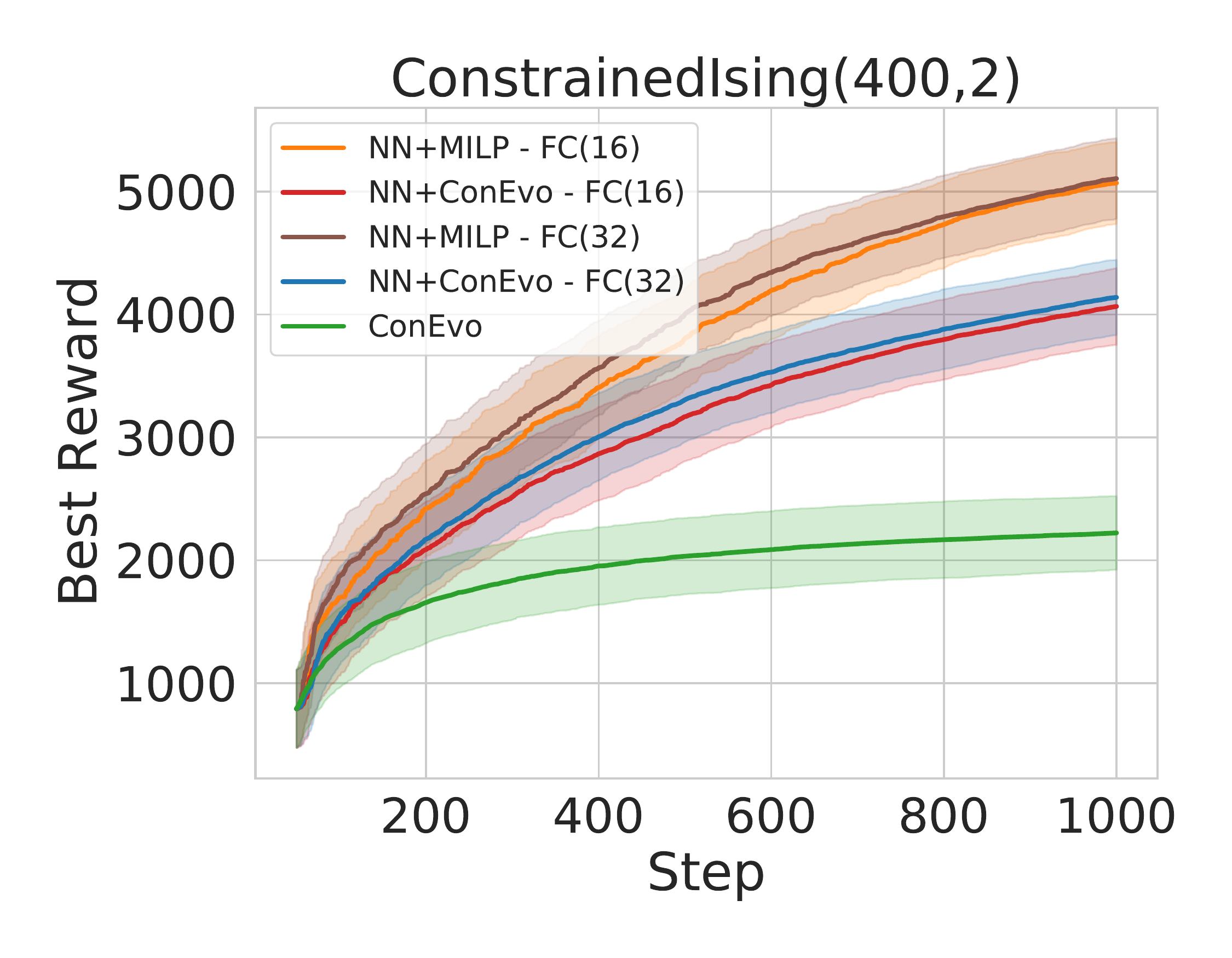}
	\caption{Best observed reward as a function of iteration for constrained Ising models with $n=200$ (left) and $n=400$ (right) variables, as in Section~\ref{sec:exp_constr}. FC(16) and FC(32) represent the runs where the surrogate neural network has a single layer of 16 and 32 ReLUs respectively. Note that the error bands are larger here since we aggregate over all instances within each problem size.}\label{fig:constrained_ising_model_larger_withfc32}
\end{figure*}

\subsection{Constrained Optimization}
\label{sec:appendix_results_constr}

\subsubsection{Normalized Max Reward}

\begin{figure*}[ht]
\begin{center}
    \begin{multicols}{3}
    \includegraphics[width=\linewidth]{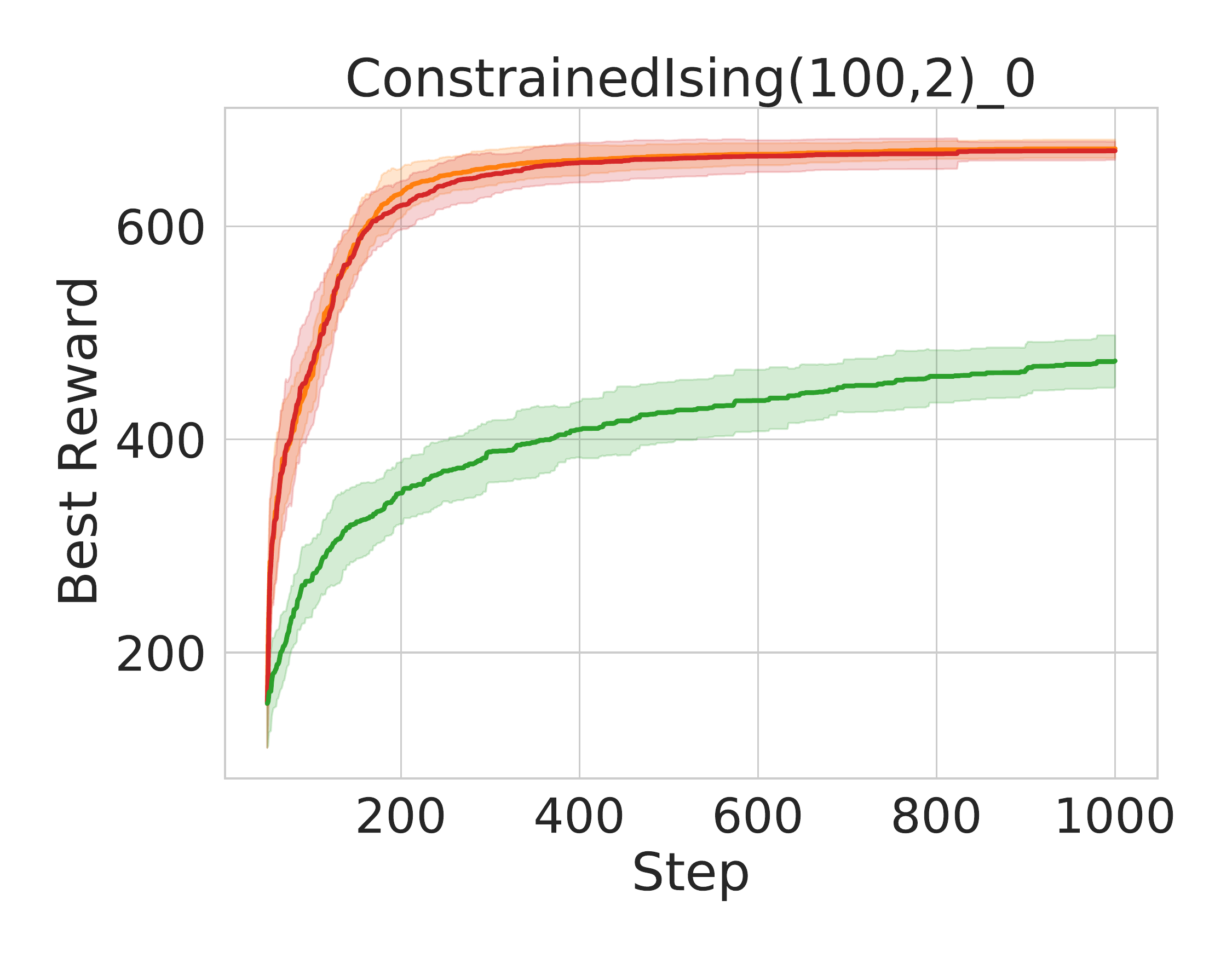}
    \includegraphics[width=\linewidth]{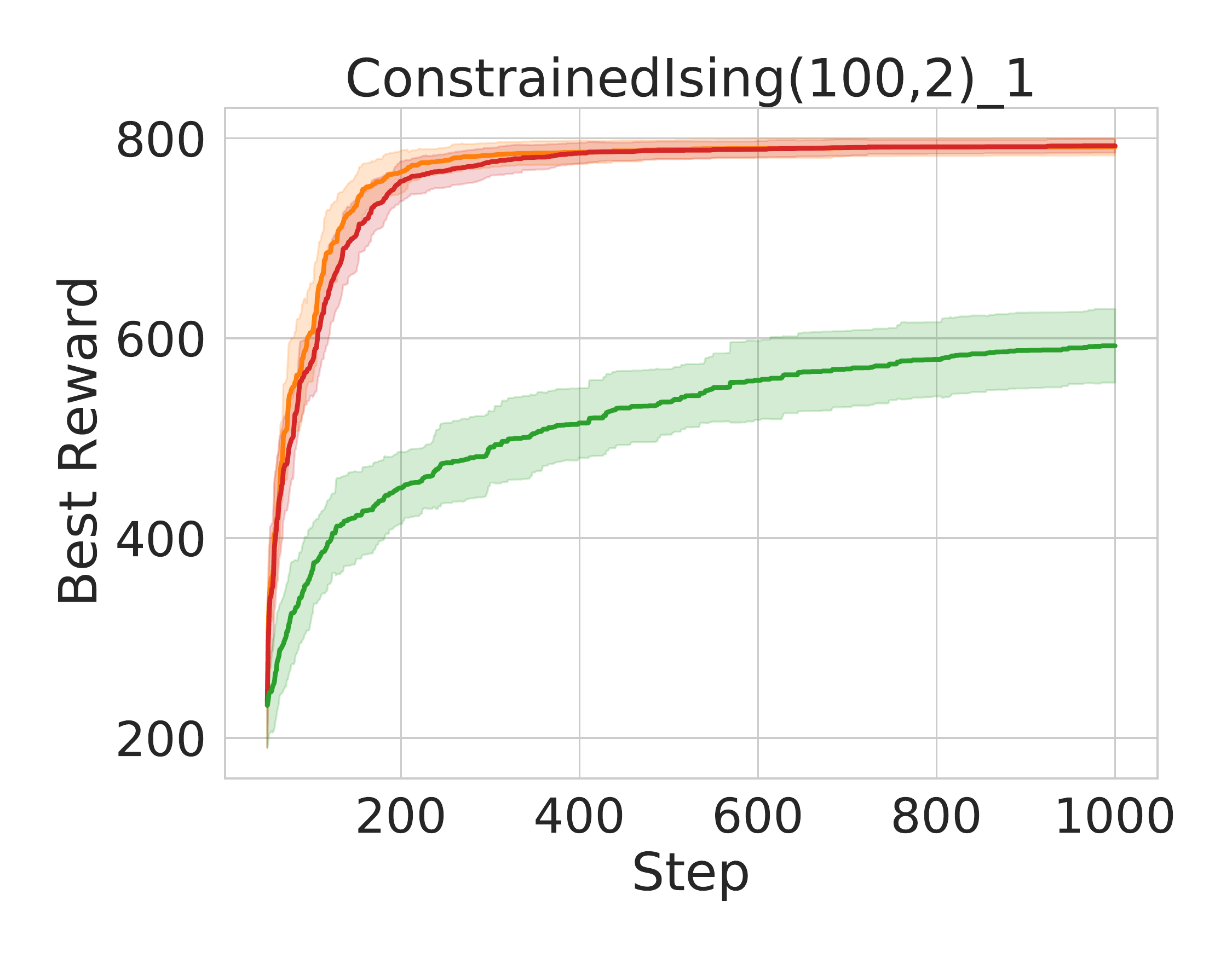}
    \includegraphics[width=\linewidth]{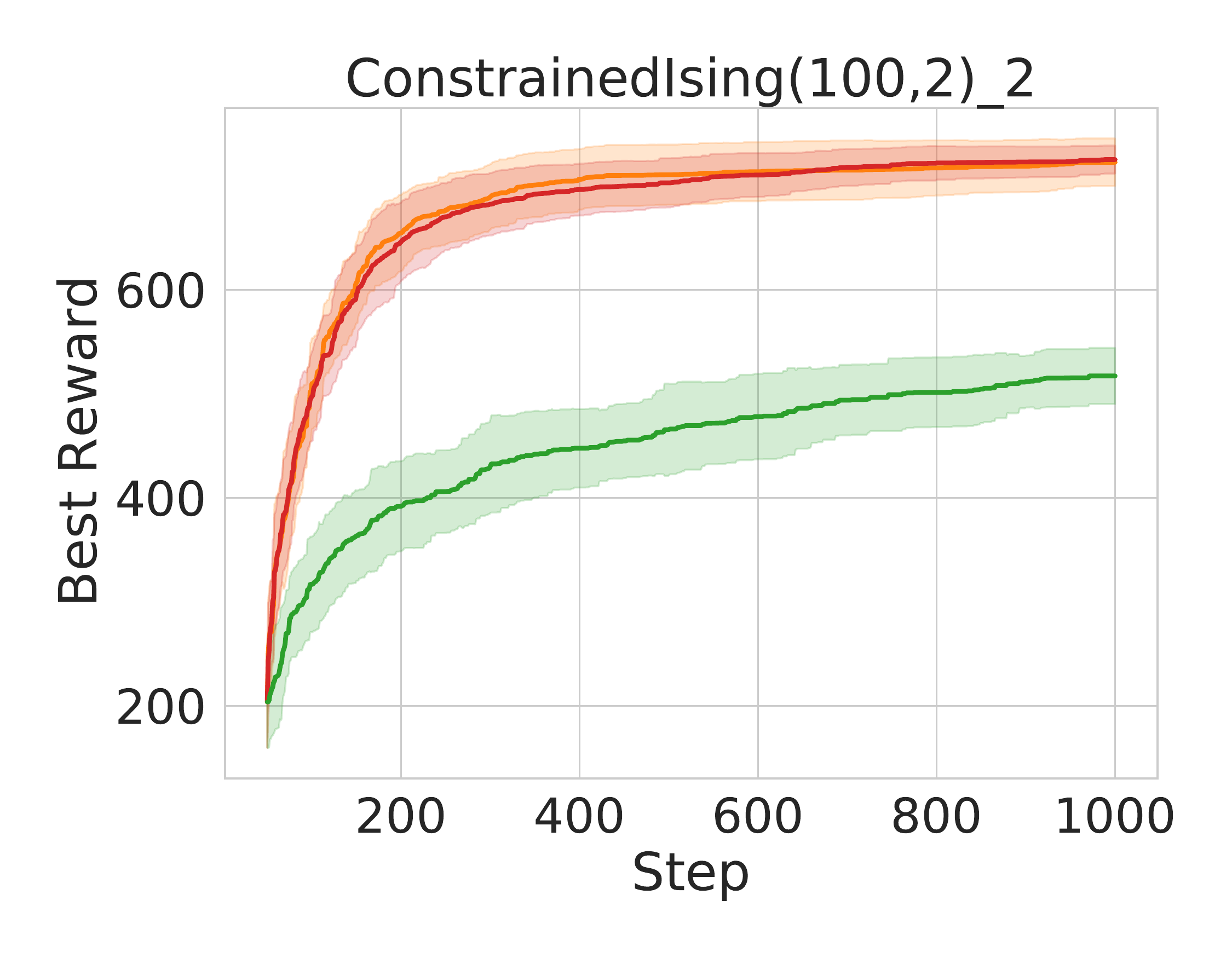}
    \end{multicols}
    \vspace{-30pt}
    \begin{multicols}{3}
    \includegraphics[width=\linewidth]{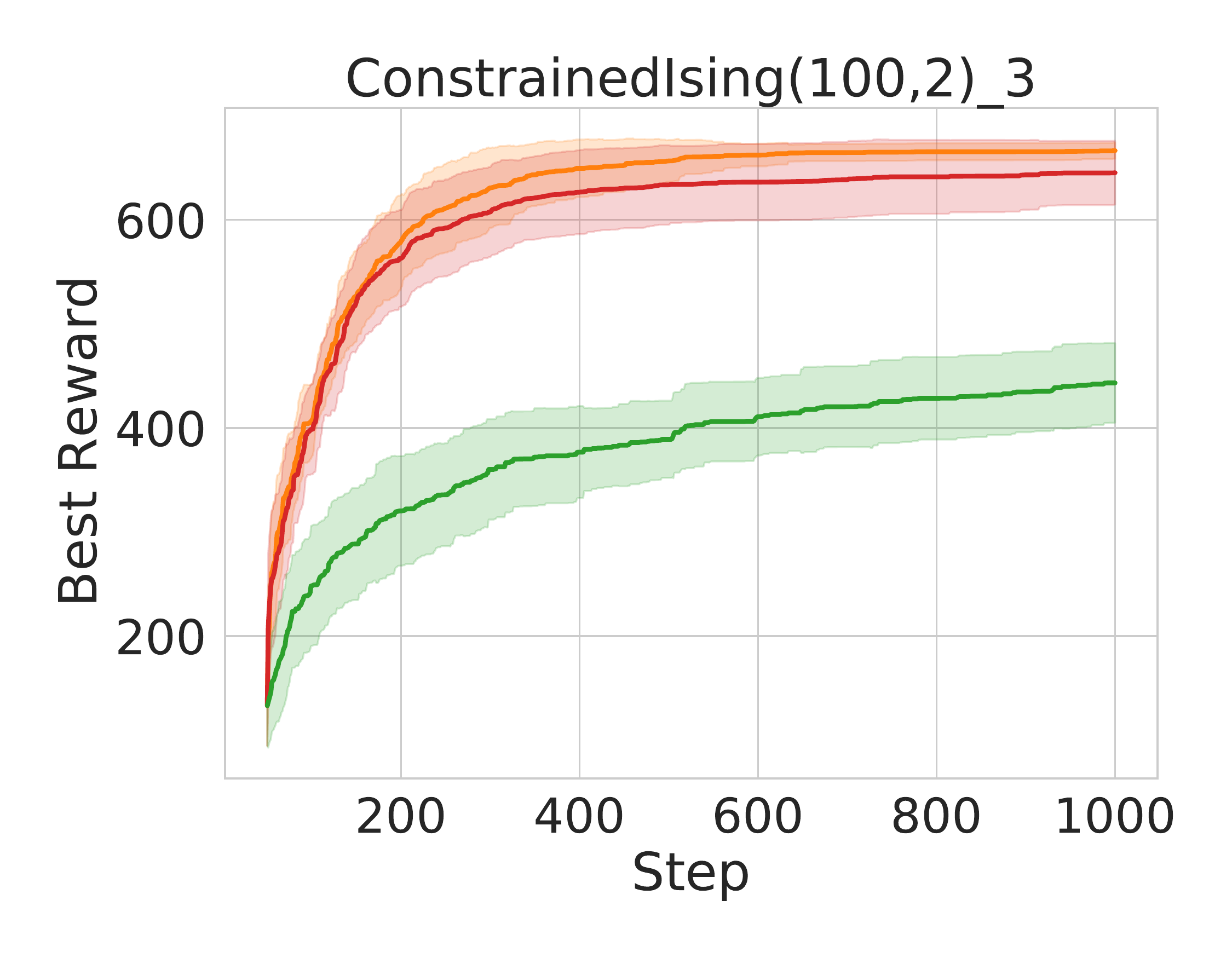}
    \includegraphics[width=\linewidth]{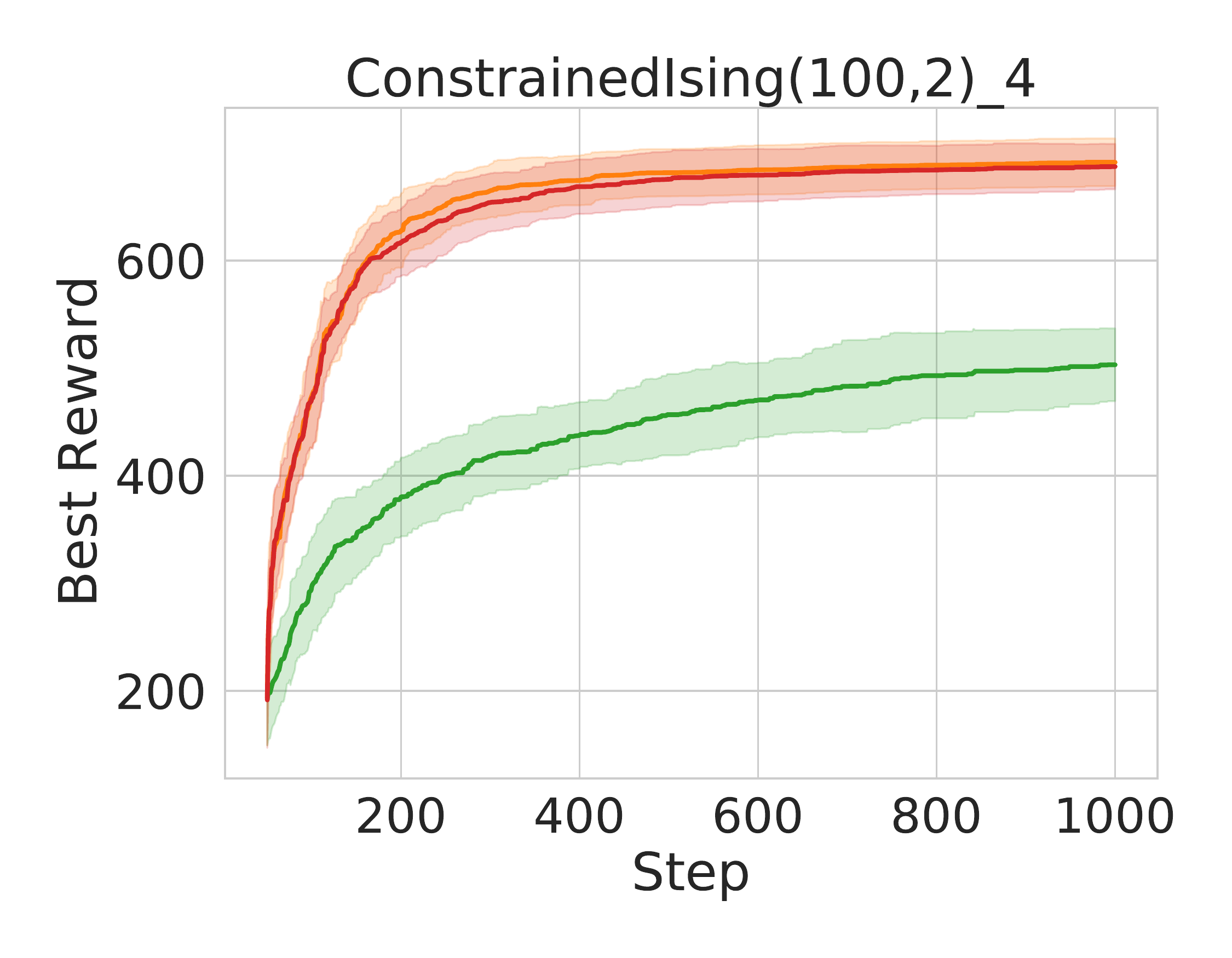}
    \includegraphics[width=\linewidth]{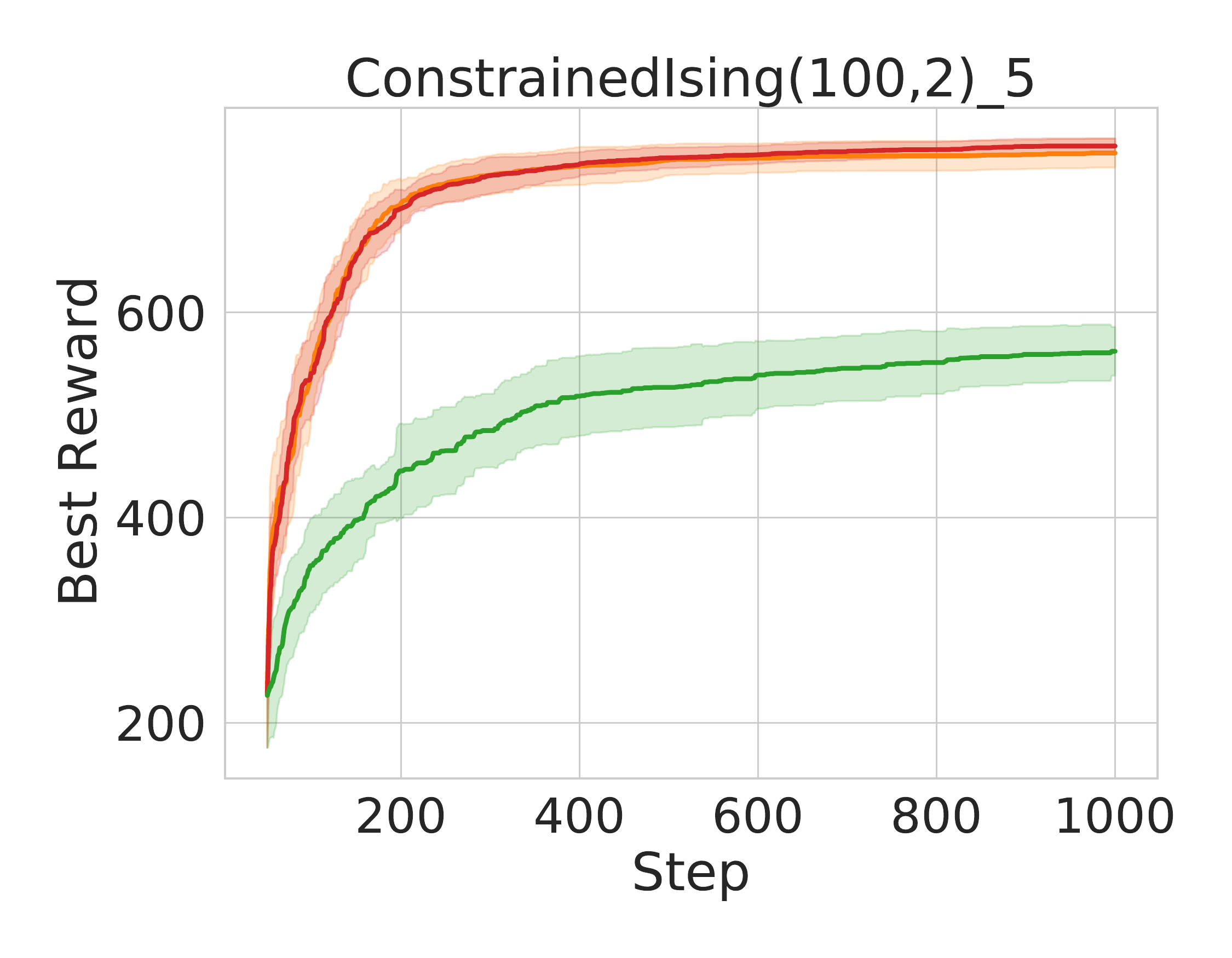}
    \end{multicols}
    \vspace{-30pt}
    \begin{multicols}{3}
    \includegraphics[width=\linewidth]{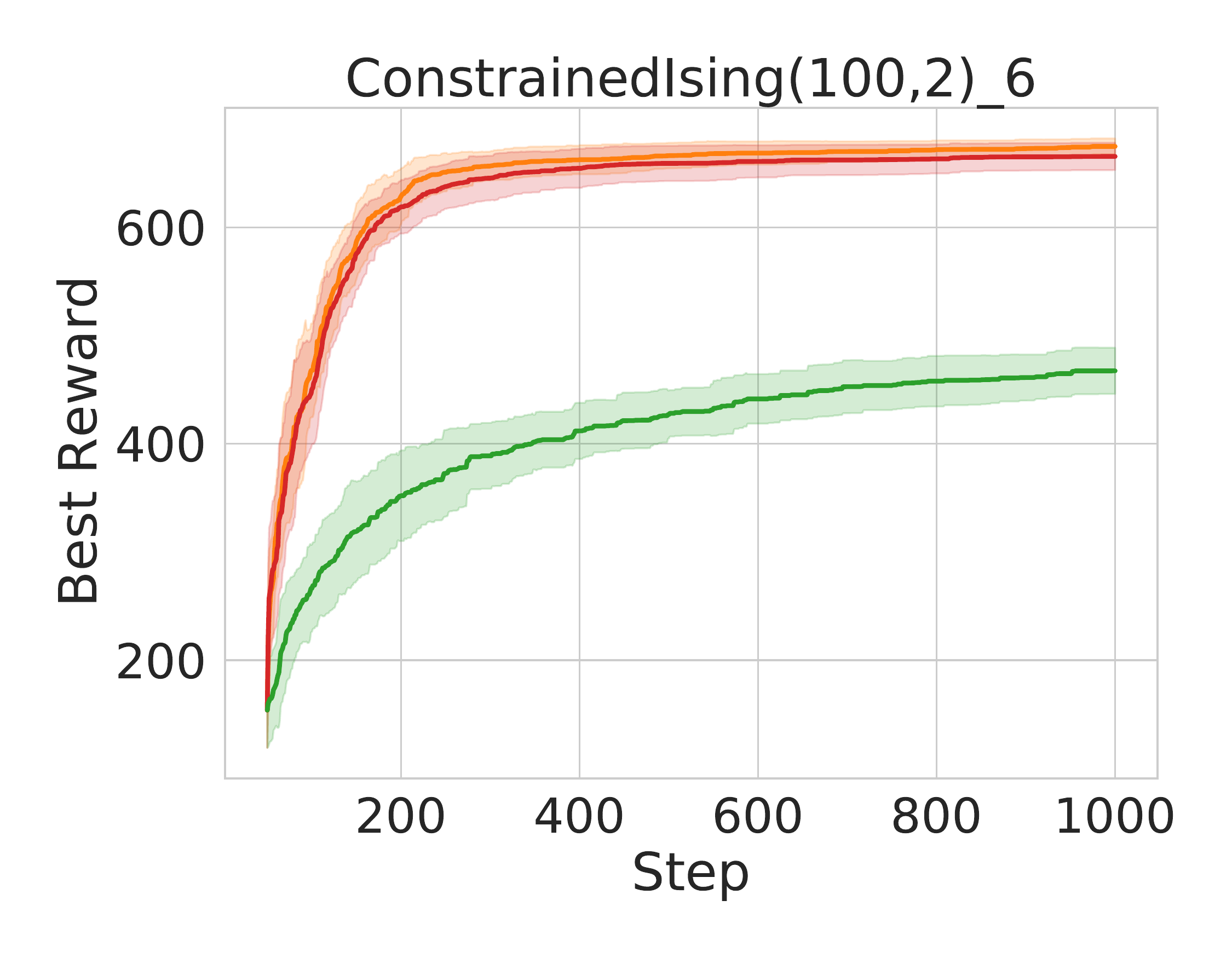}
    \includegraphics[width=\linewidth]{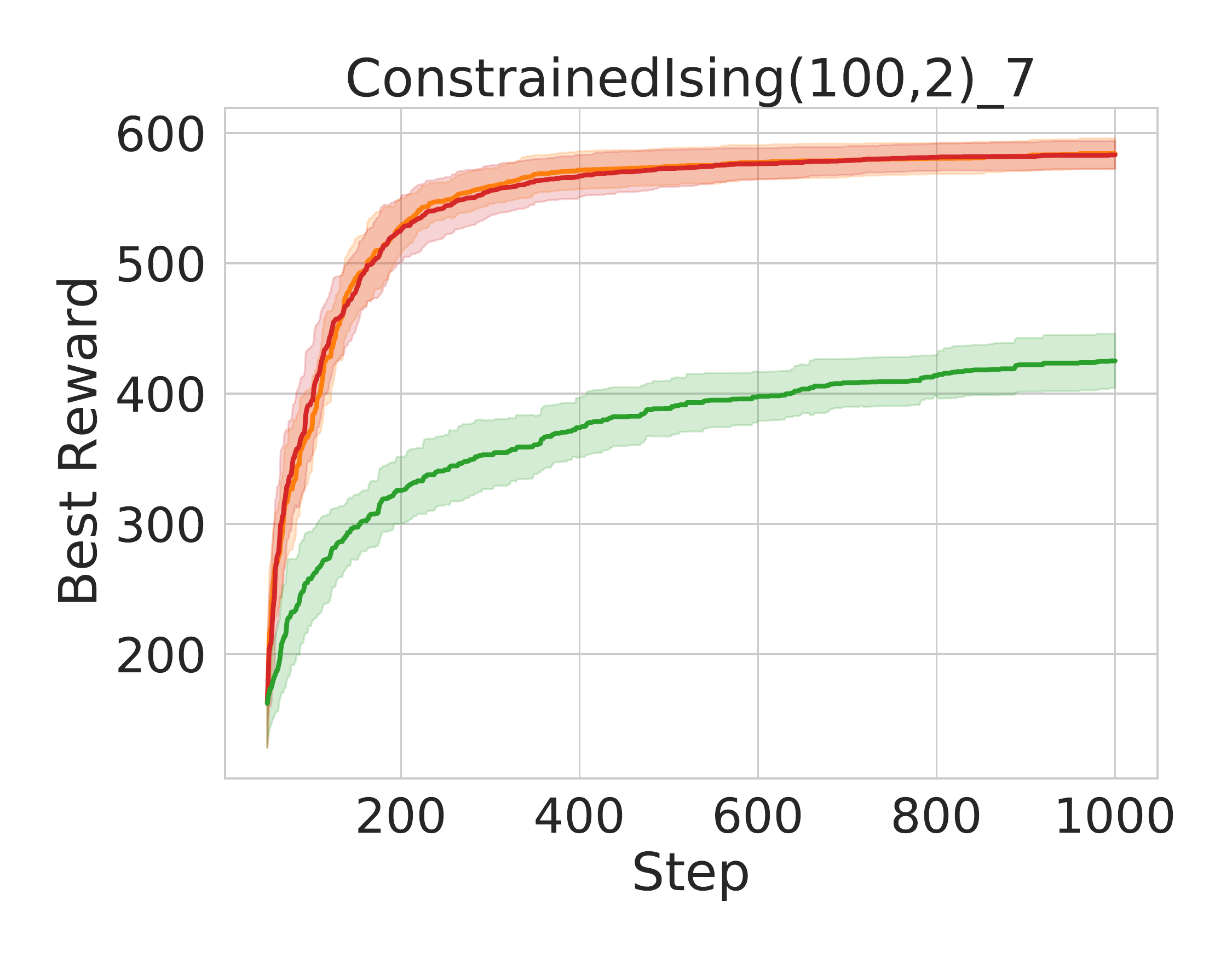}
    \includegraphics[width=\linewidth]{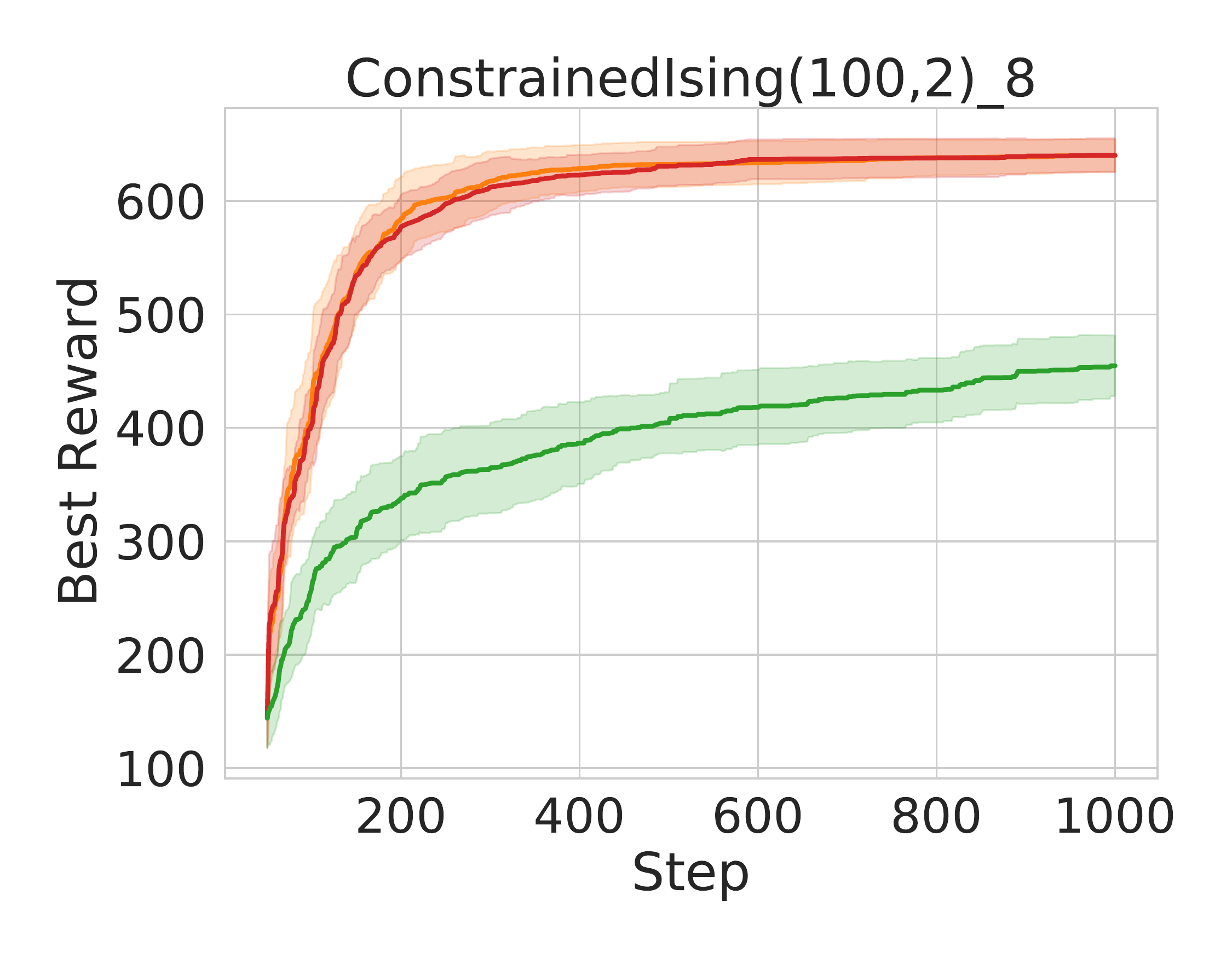}
    \end{multicols}
    \vspace{-30pt}
    \begin{multicols}{3}
    \includegraphics[width=\linewidth]{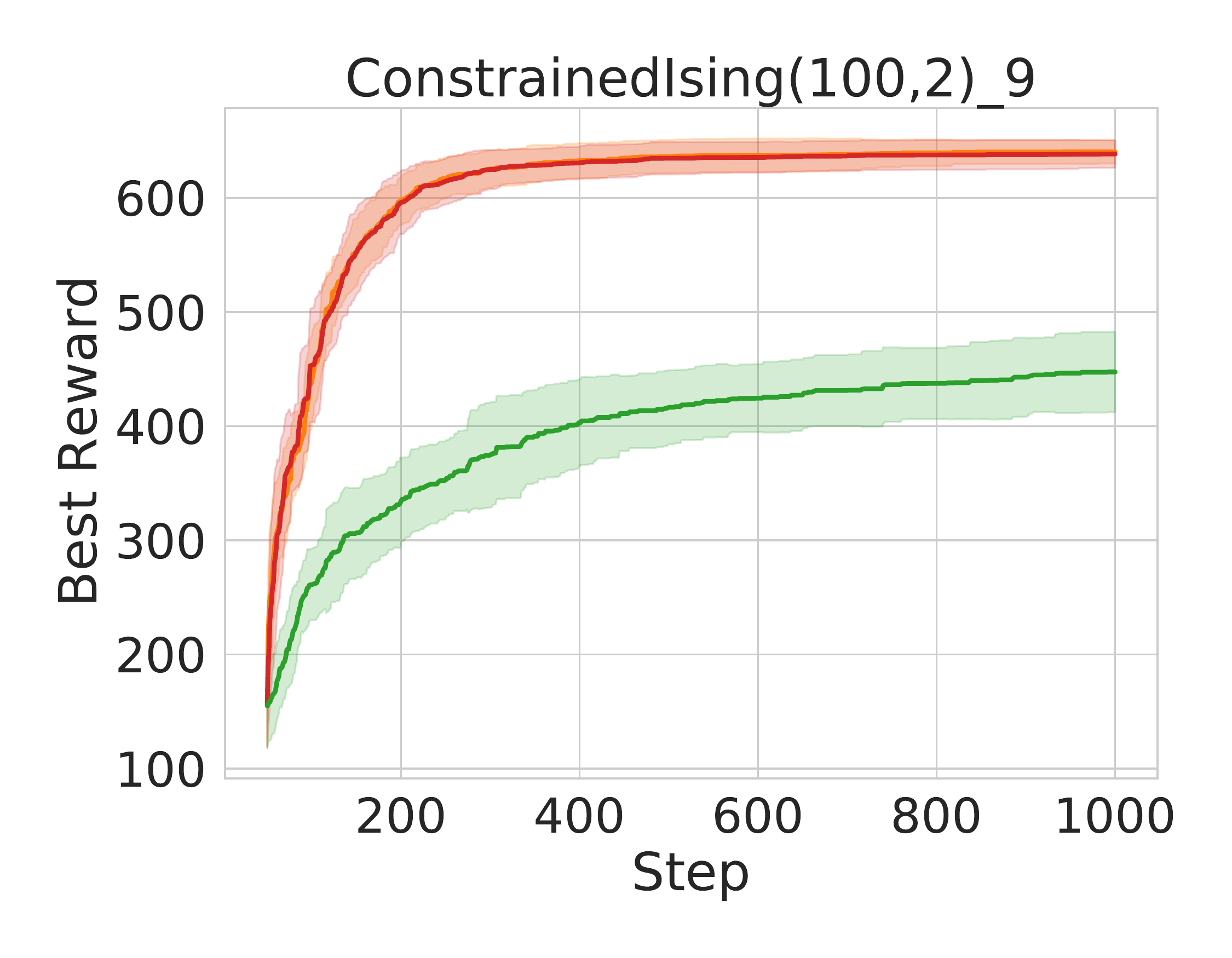}
    \includegraphics[width=\linewidth]{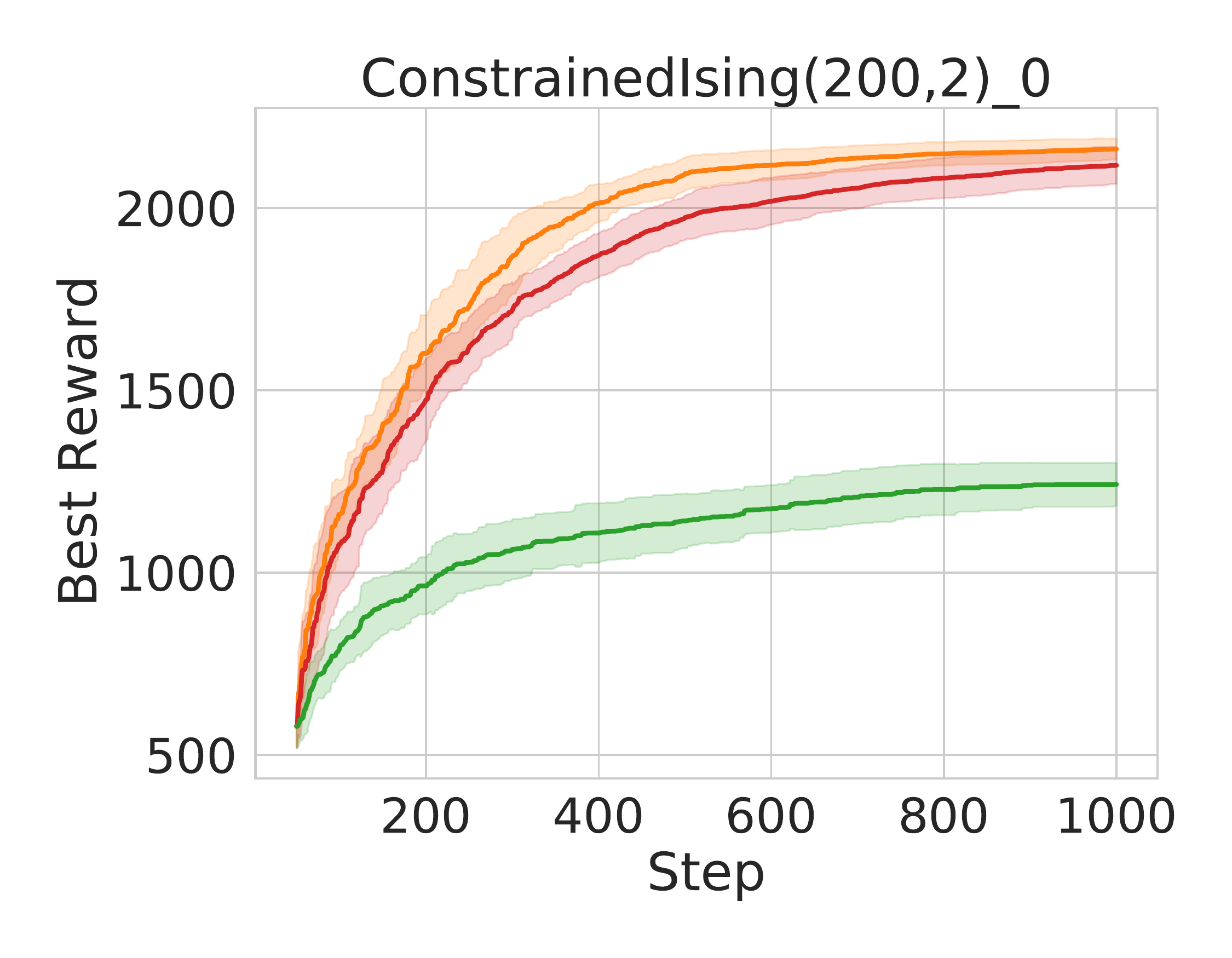}
    \includegraphics[width=\linewidth]{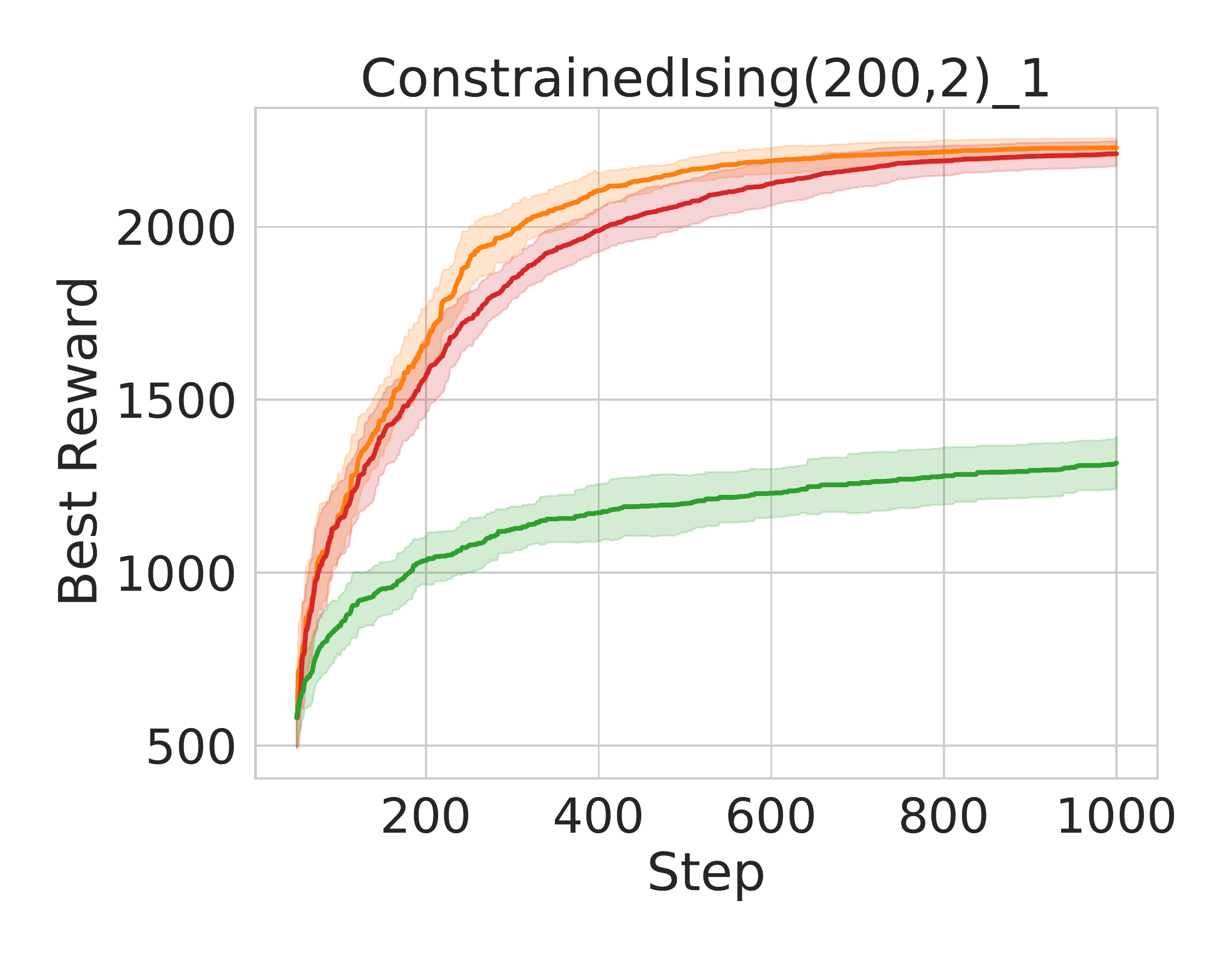}
    \end{multicols}
    \vspace{-30pt}
    \begin{multicols}{3}
    \includegraphics[width=\linewidth]{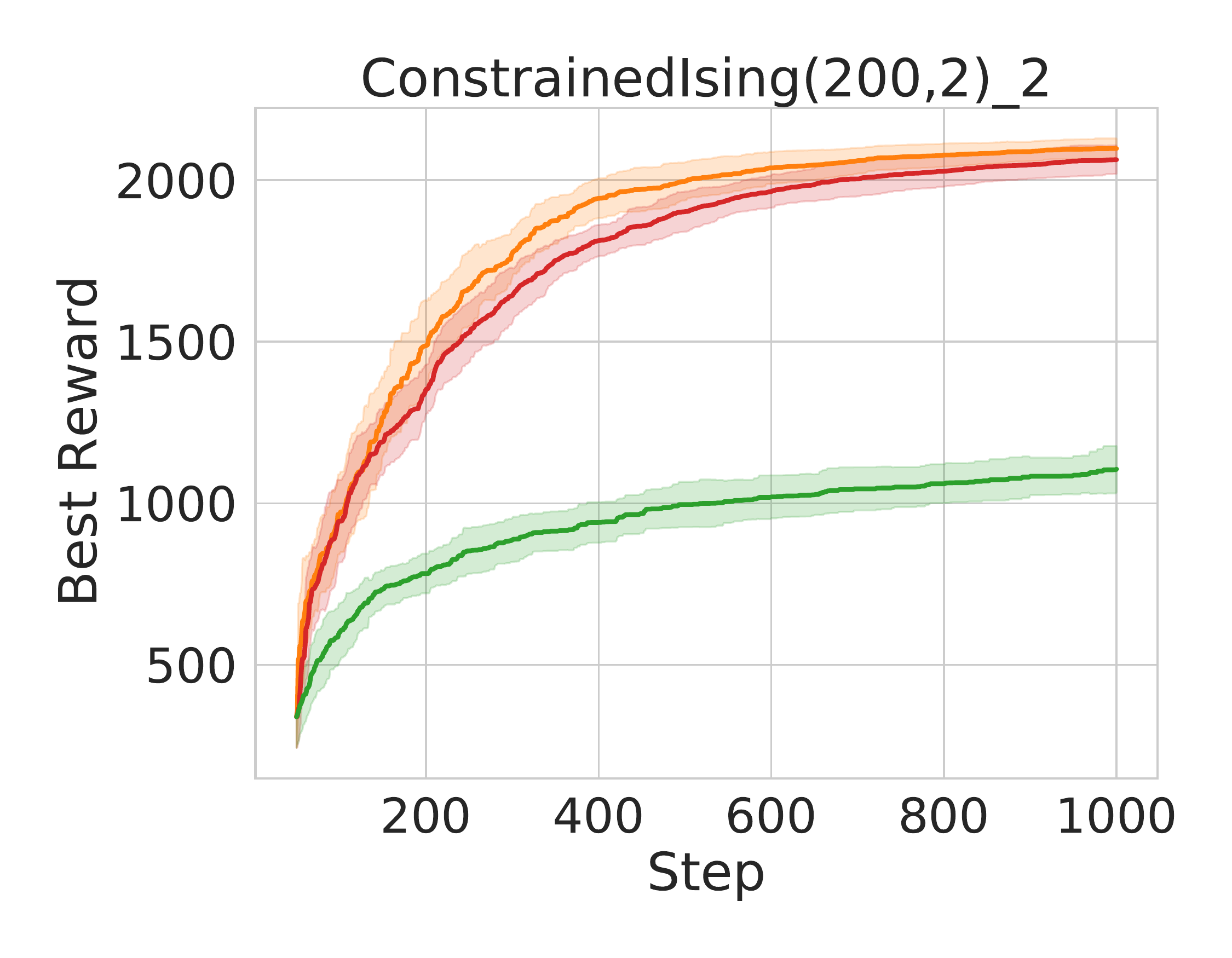}
    \includegraphics[width=\linewidth]{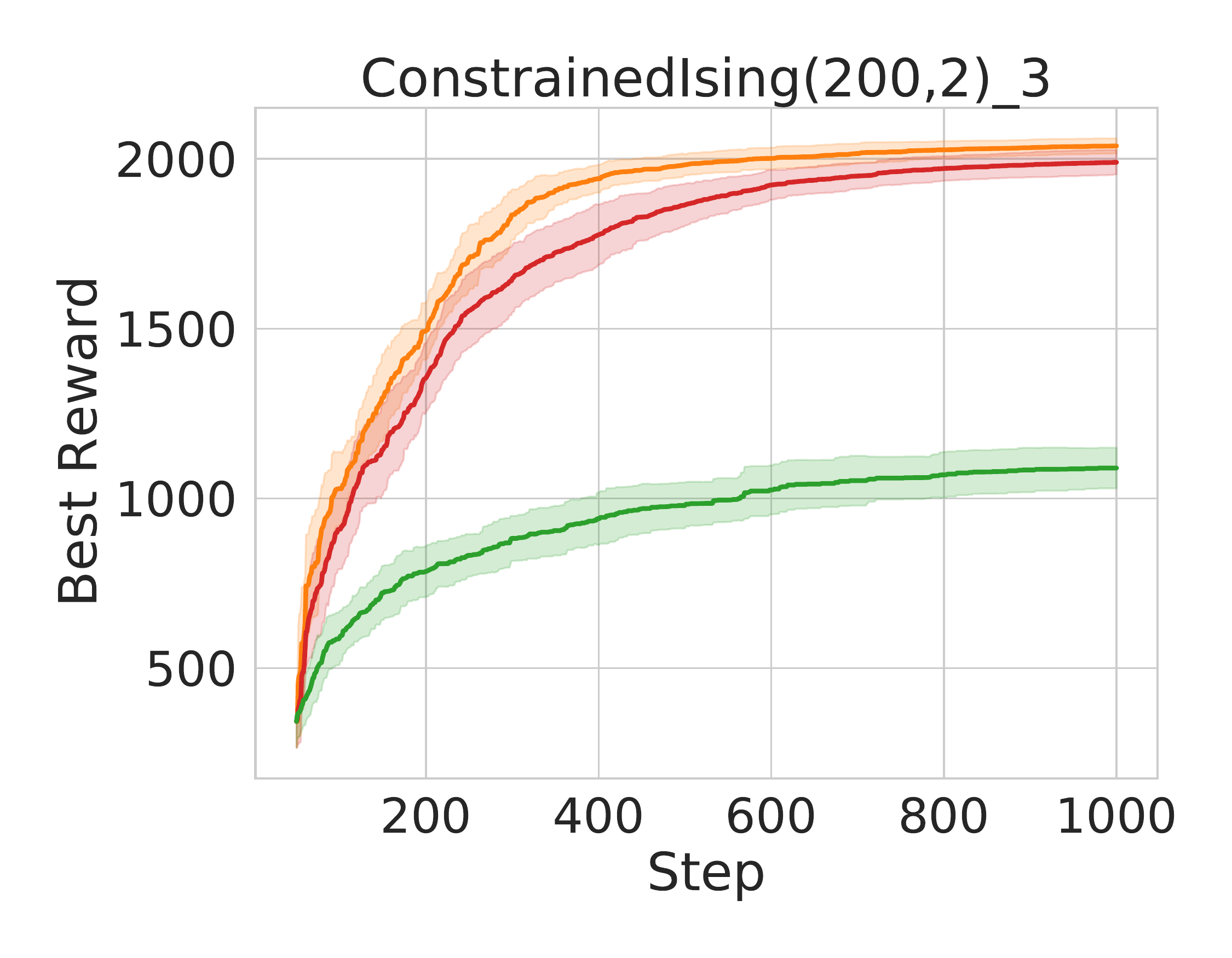}
    \includegraphics[width=\linewidth]{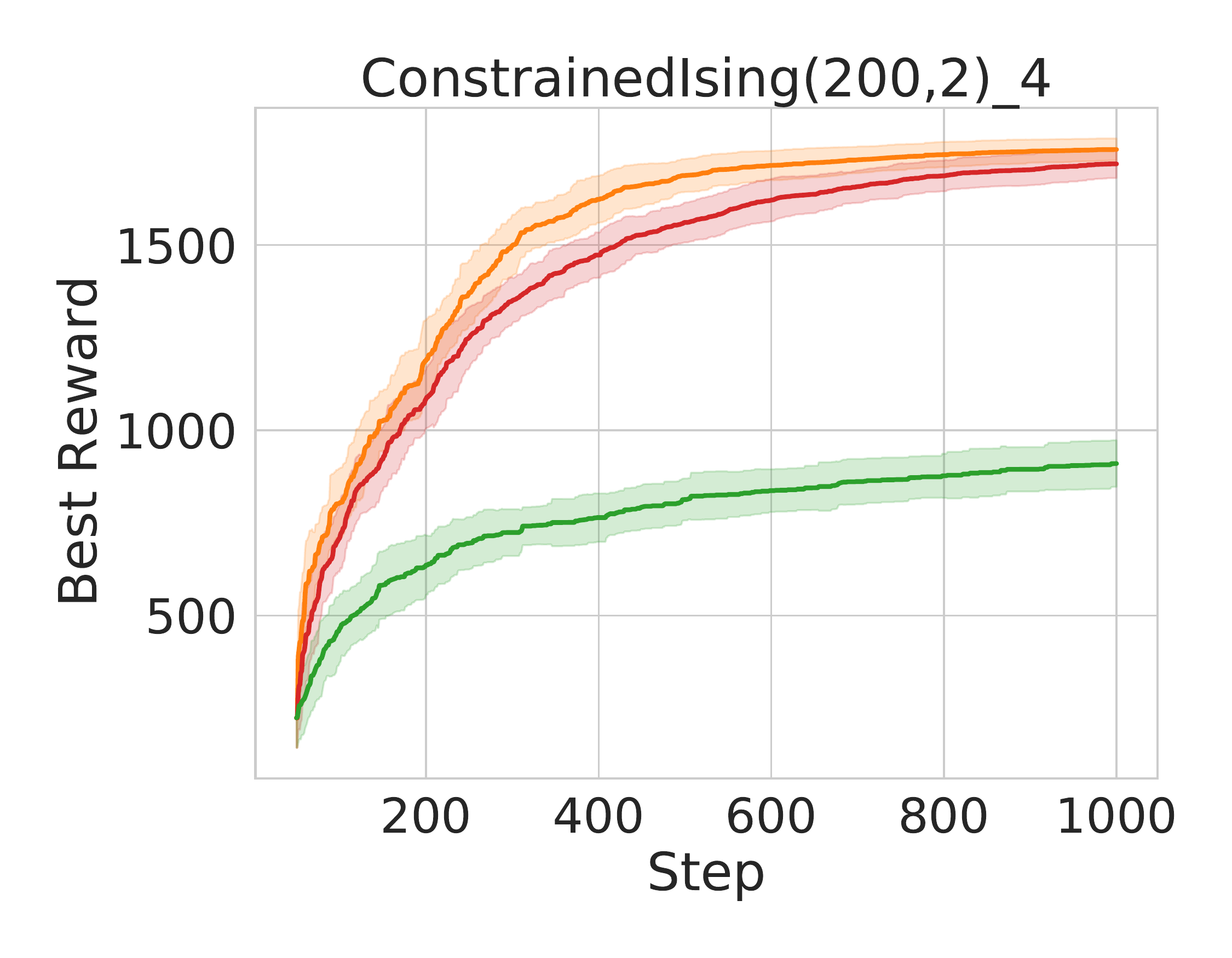}
    \end{multicols}
    \vspace{-20pt}
\centerline{\includegraphics[width=0.9\textwidth]{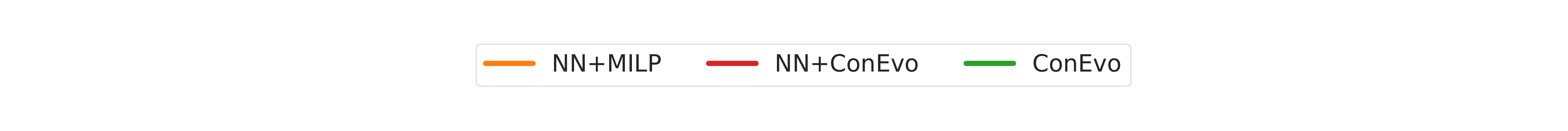}}
\vspace{-15pt}
\caption{Best observed reward as a function of iteration for the first half of all constrained Ising model problems (Section~\ref{sec:exp_constr}), averaged over 20 or 10 trials for $n=100$ or $n \geq 200$ respectively (bands indicate $\pm 1$sd). Initial randomly sampled set of 50 points is omitted.}
\label{fig:appendix_constrained_curves1}
\end{center}
\end{figure*}

\begin{figure*}[ht]
\begin{center}
    \begin{multicols}{3}
    \includegraphics[width=\linewidth]{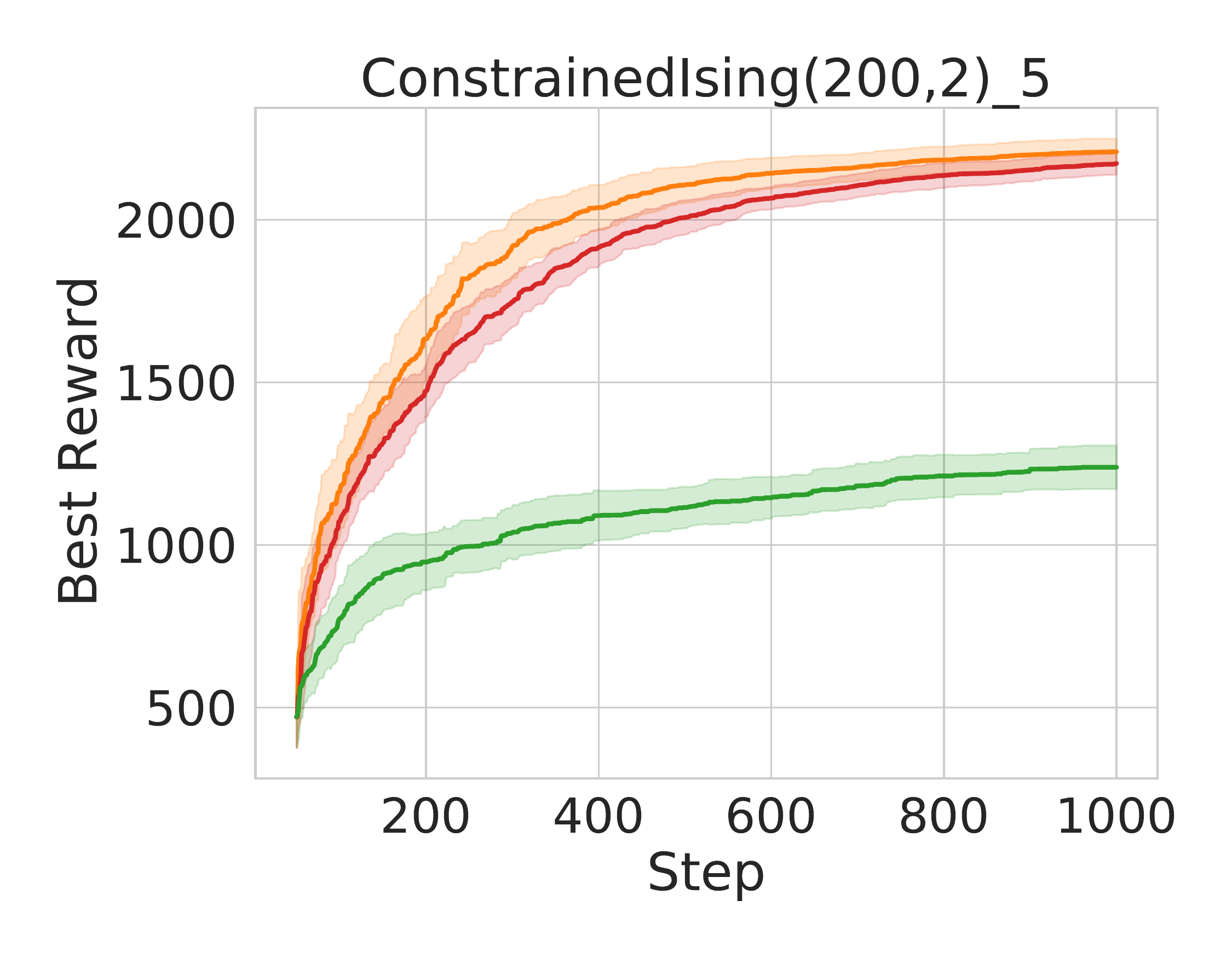}
    \includegraphics[width=\linewidth]{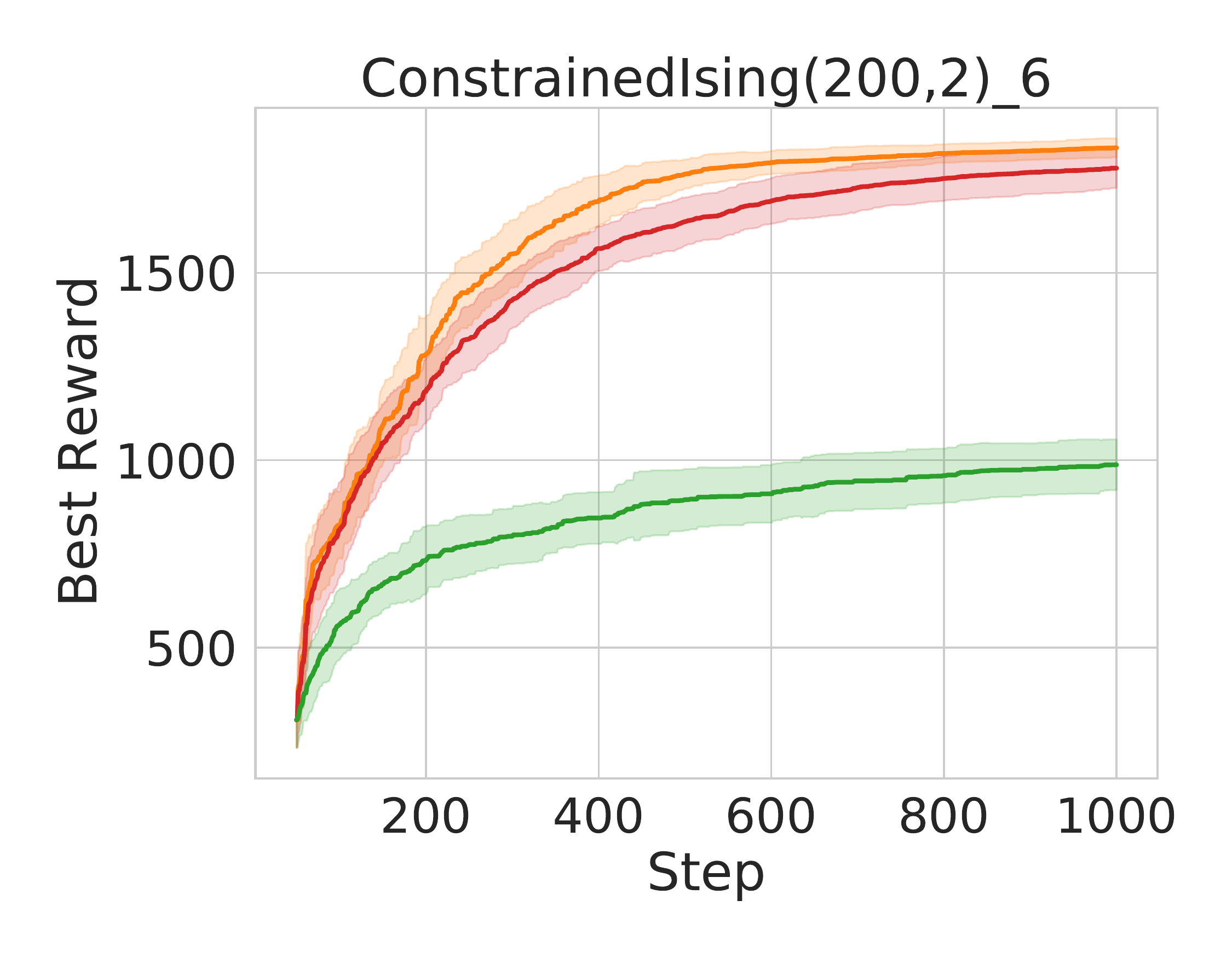}
    \includegraphics[width=\linewidth]{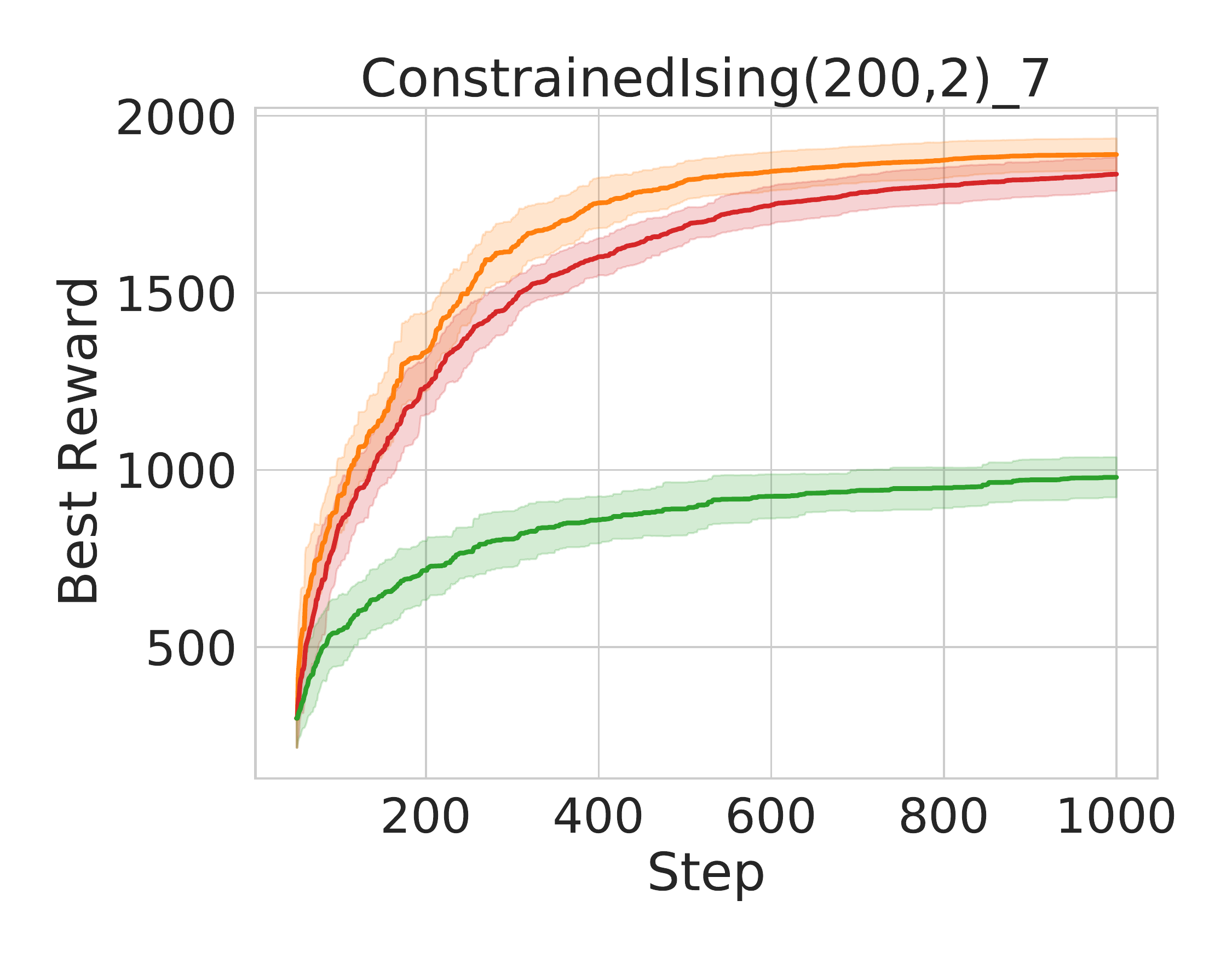}
    \end{multicols}
    \vspace{-30pt}
    \begin{multicols}{3}
    \includegraphics[width=\linewidth]{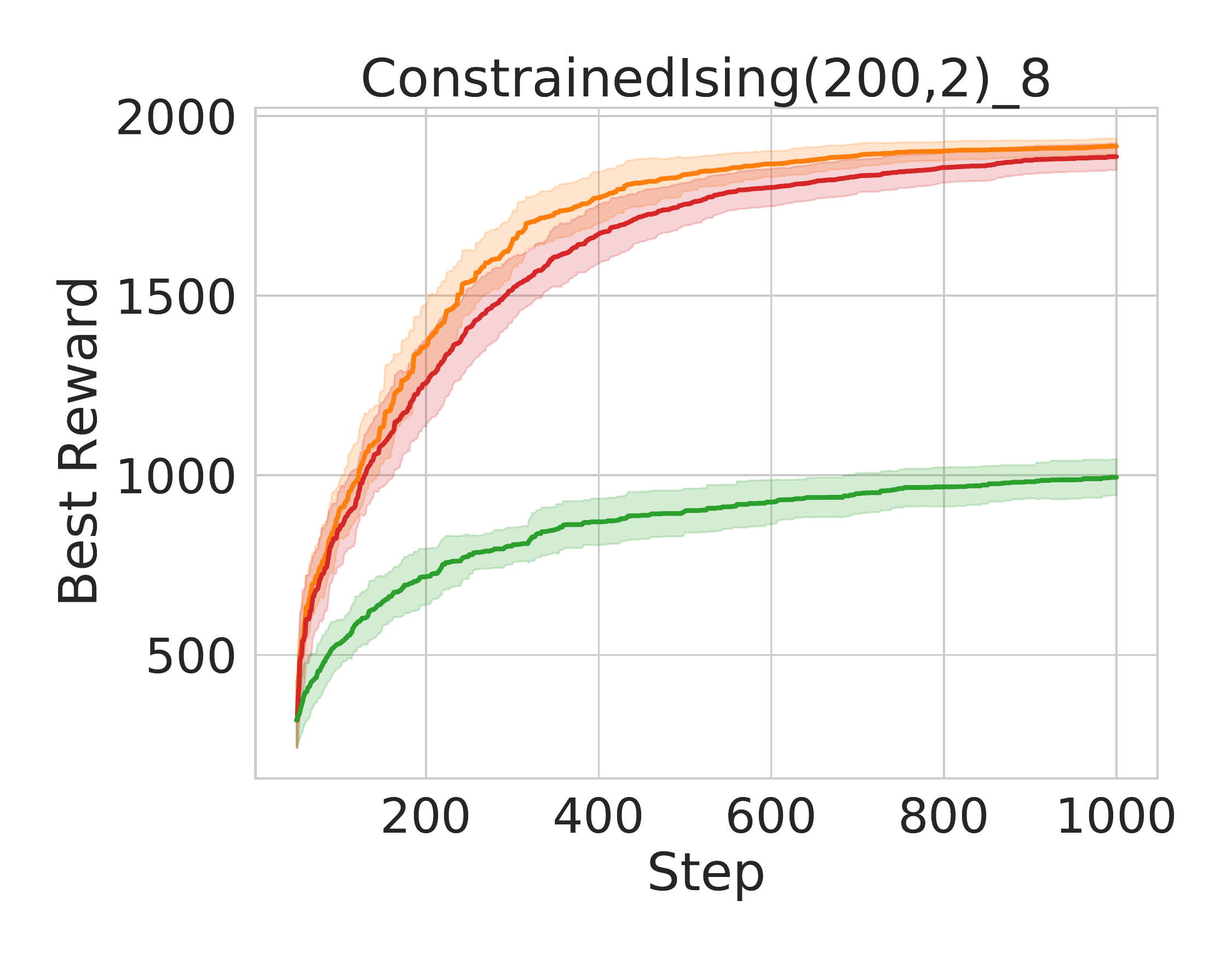}
    \includegraphics[width=\linewidth]{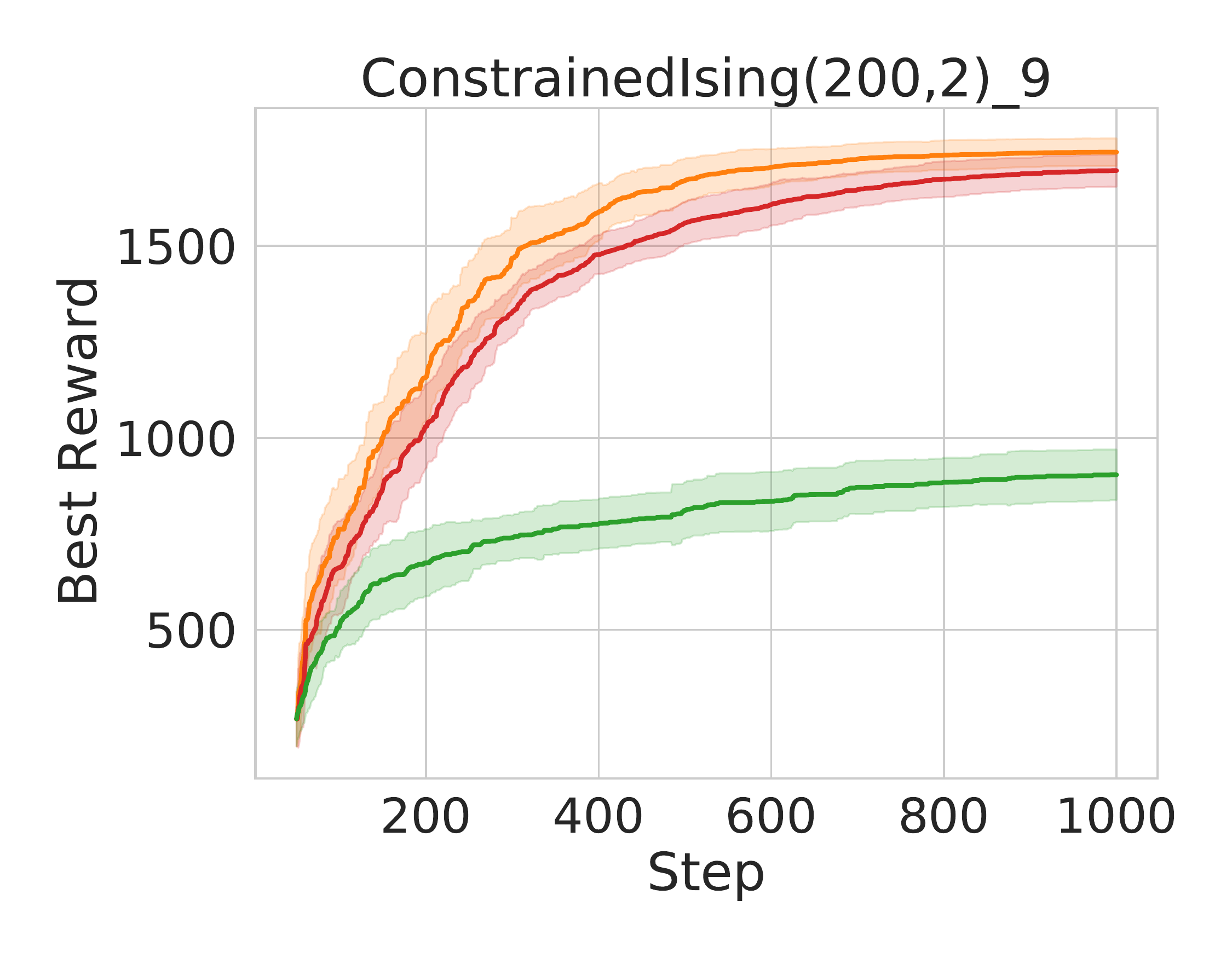}
    \includegraphics[width=\linewidth]{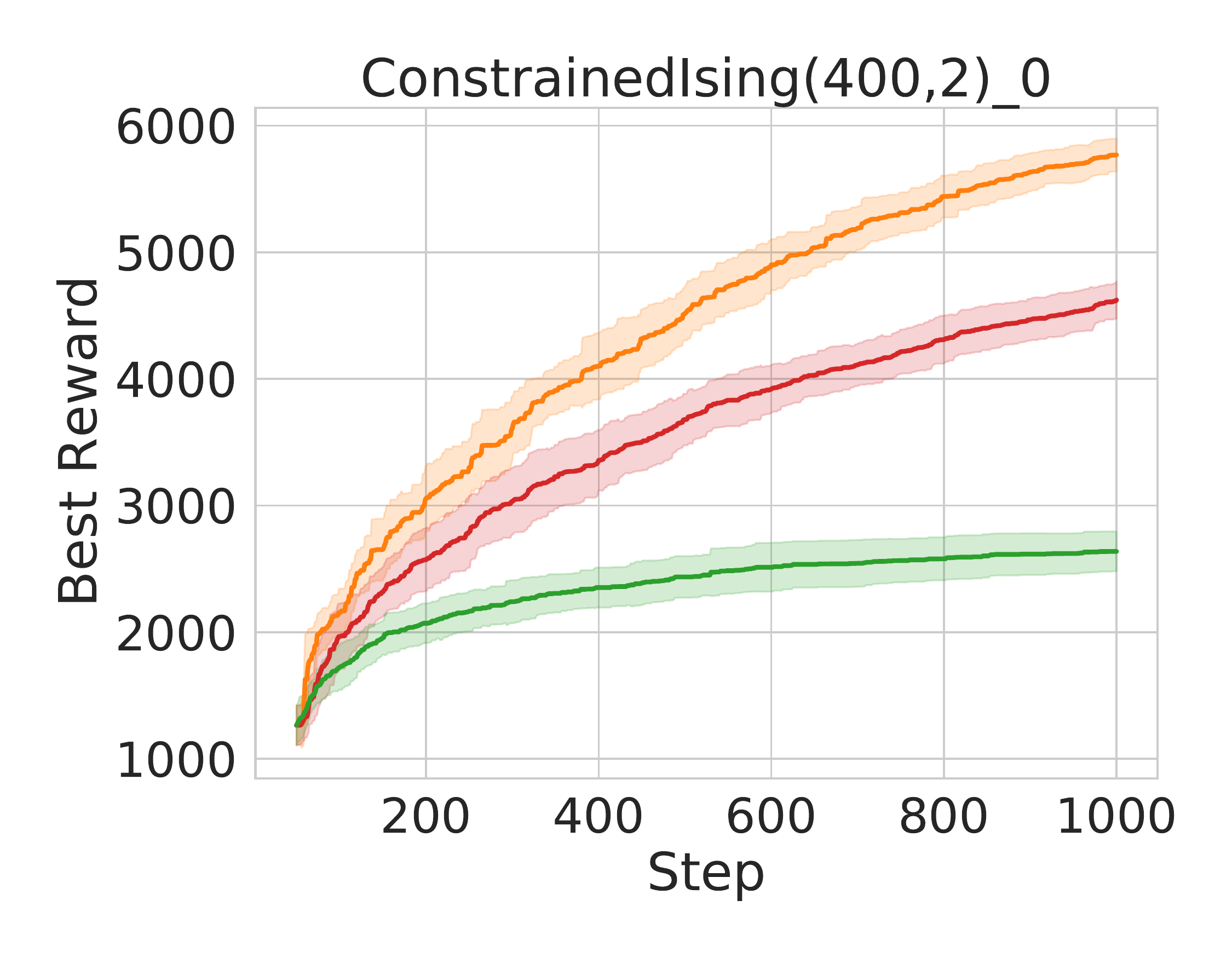}
    \end{multicols}
    \vspace{-30pt}
    \begin{multicols}{3}
    \includegraphics[width=\linewidth]{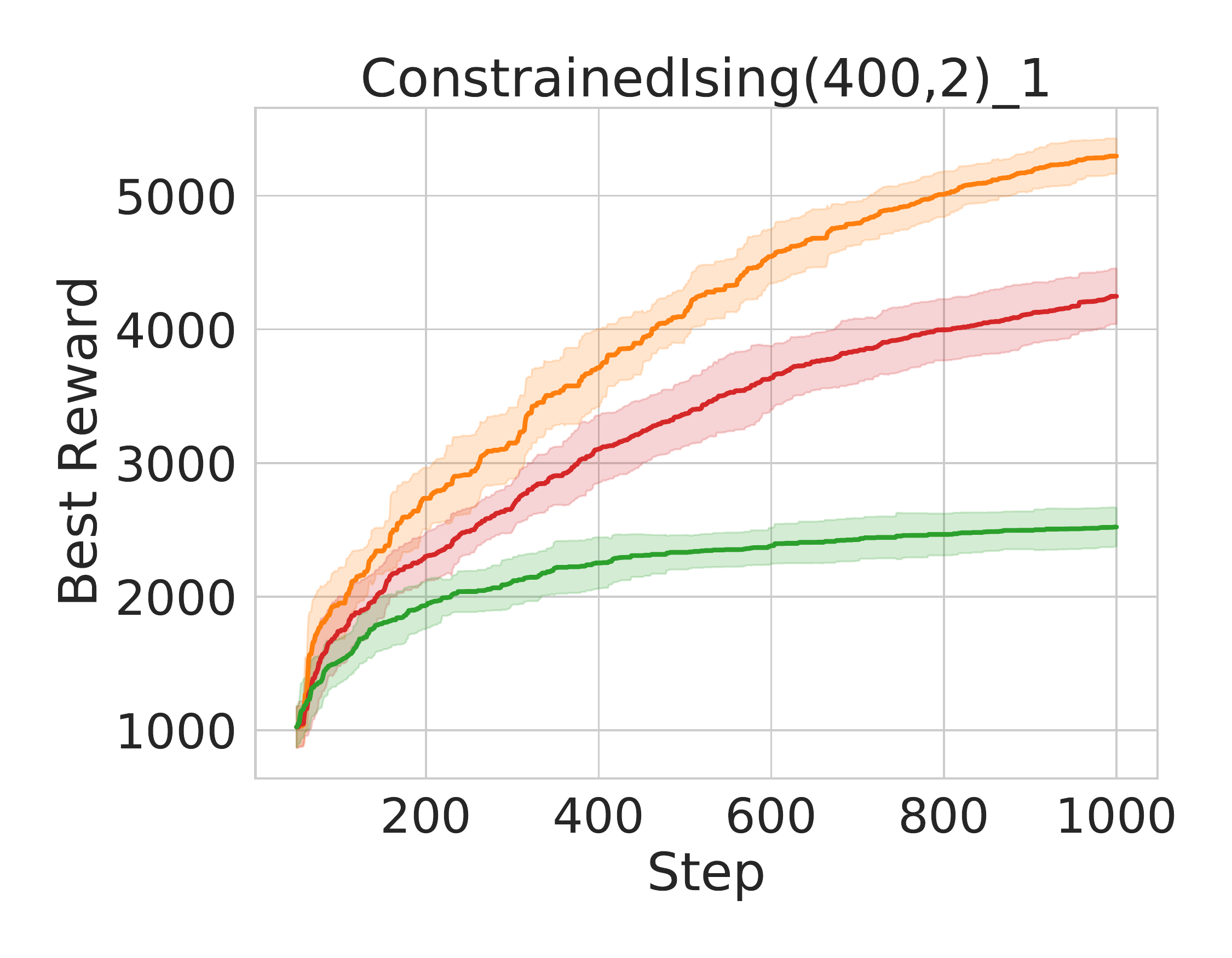}
    \includegraphics[width=\linewidth]{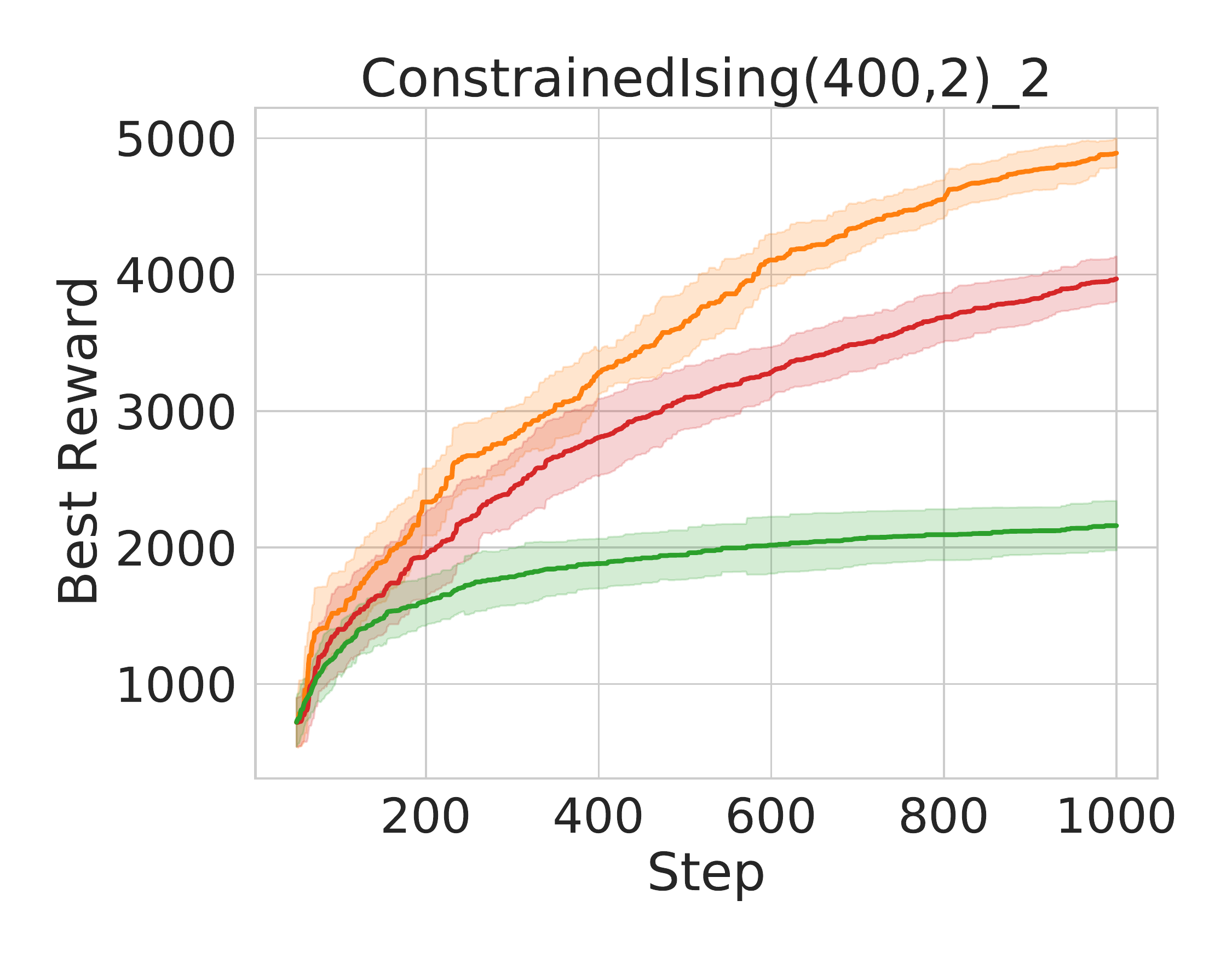}
    \includegraphics[width=\linewidth]{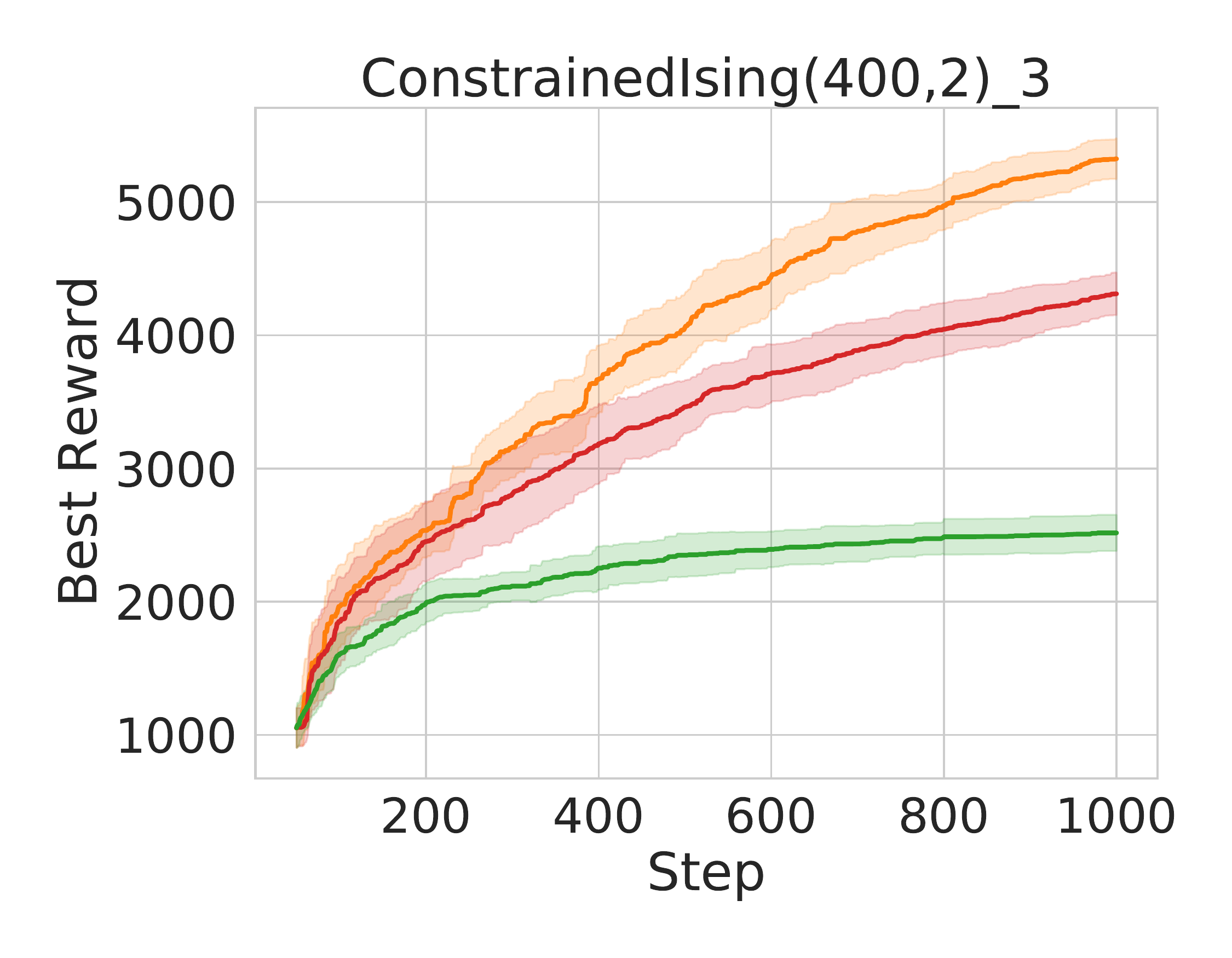}
    \end{multicols}
    \vspace{-30pt}
    \begin{multicols}{3}
    \includegraphics[width=\linewidth]{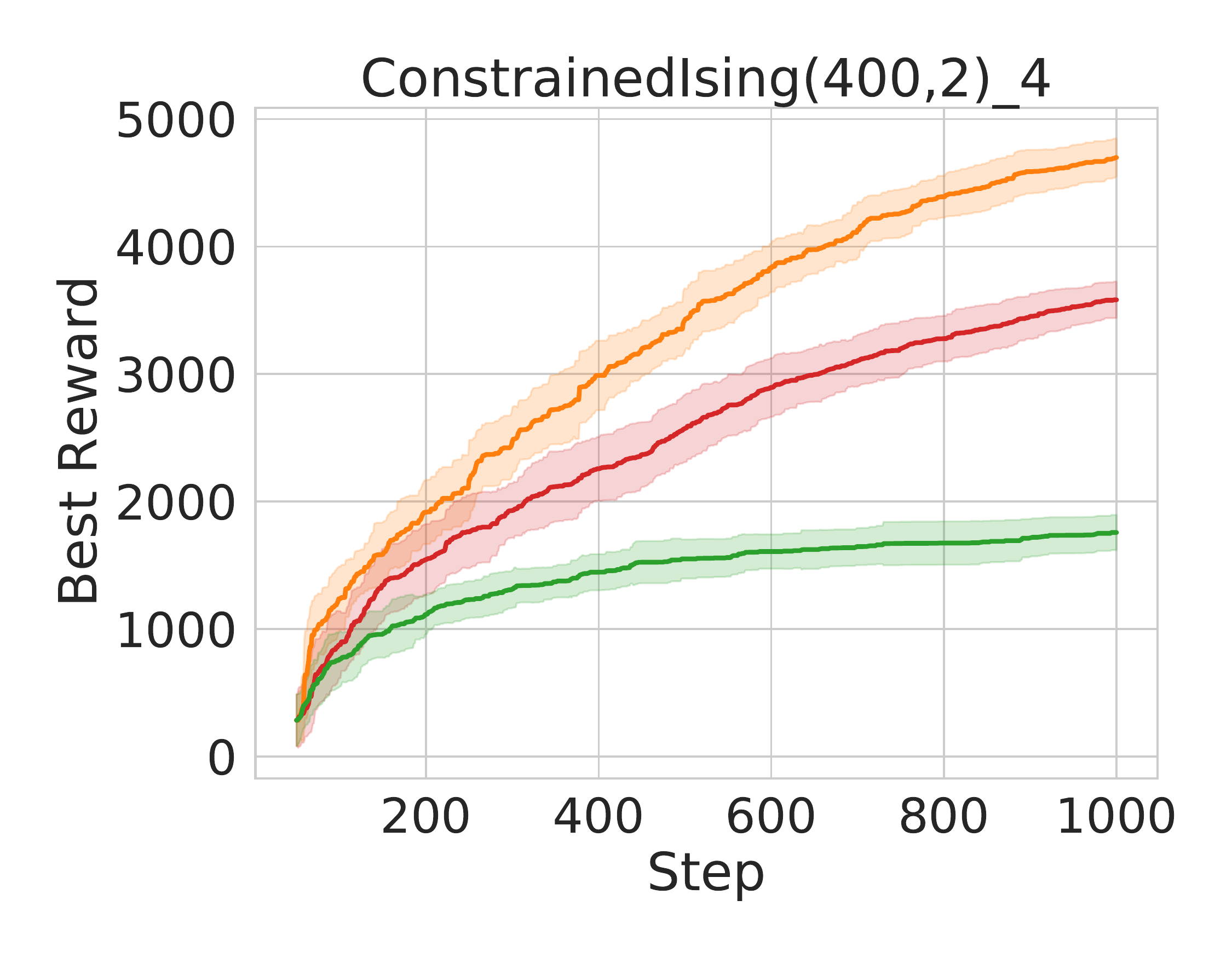}
    \includegraphics[width=\linewidth]{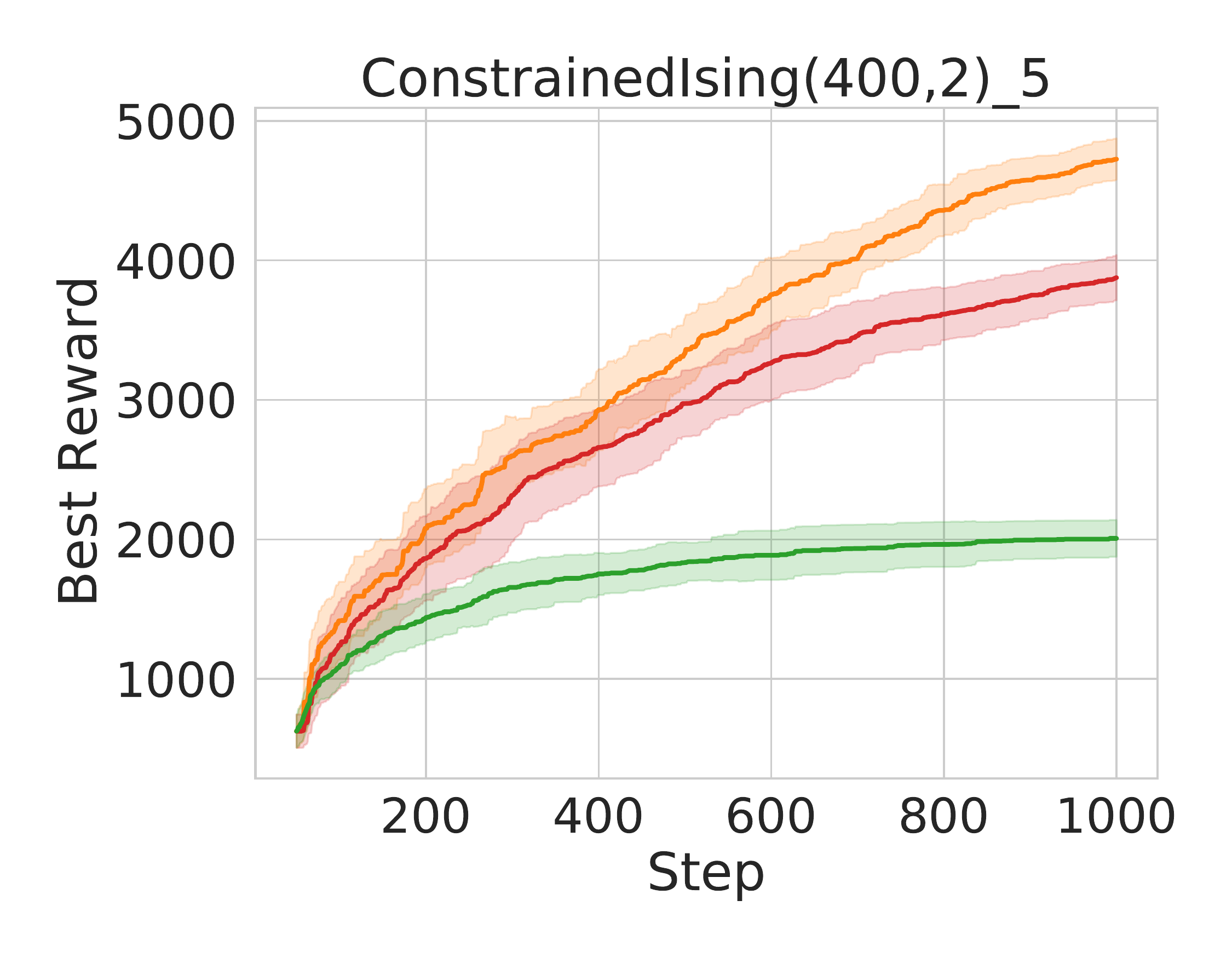}
    \includegraphics[width=\linewidth]{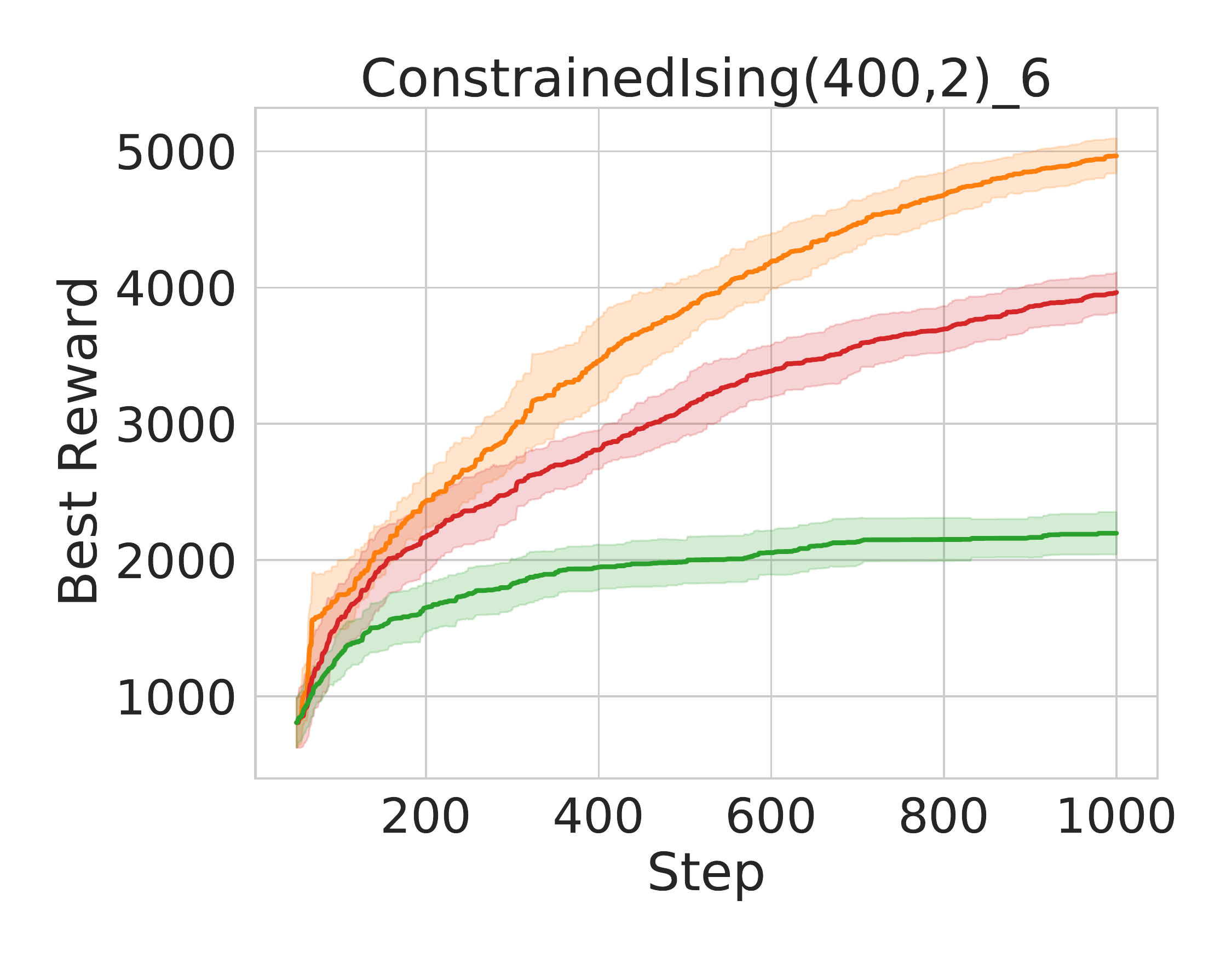}
    \end{multicols}
    \vspace{-30pt}
    \begin{multicols}{3}
    \includegraphics[width=\linewidth]{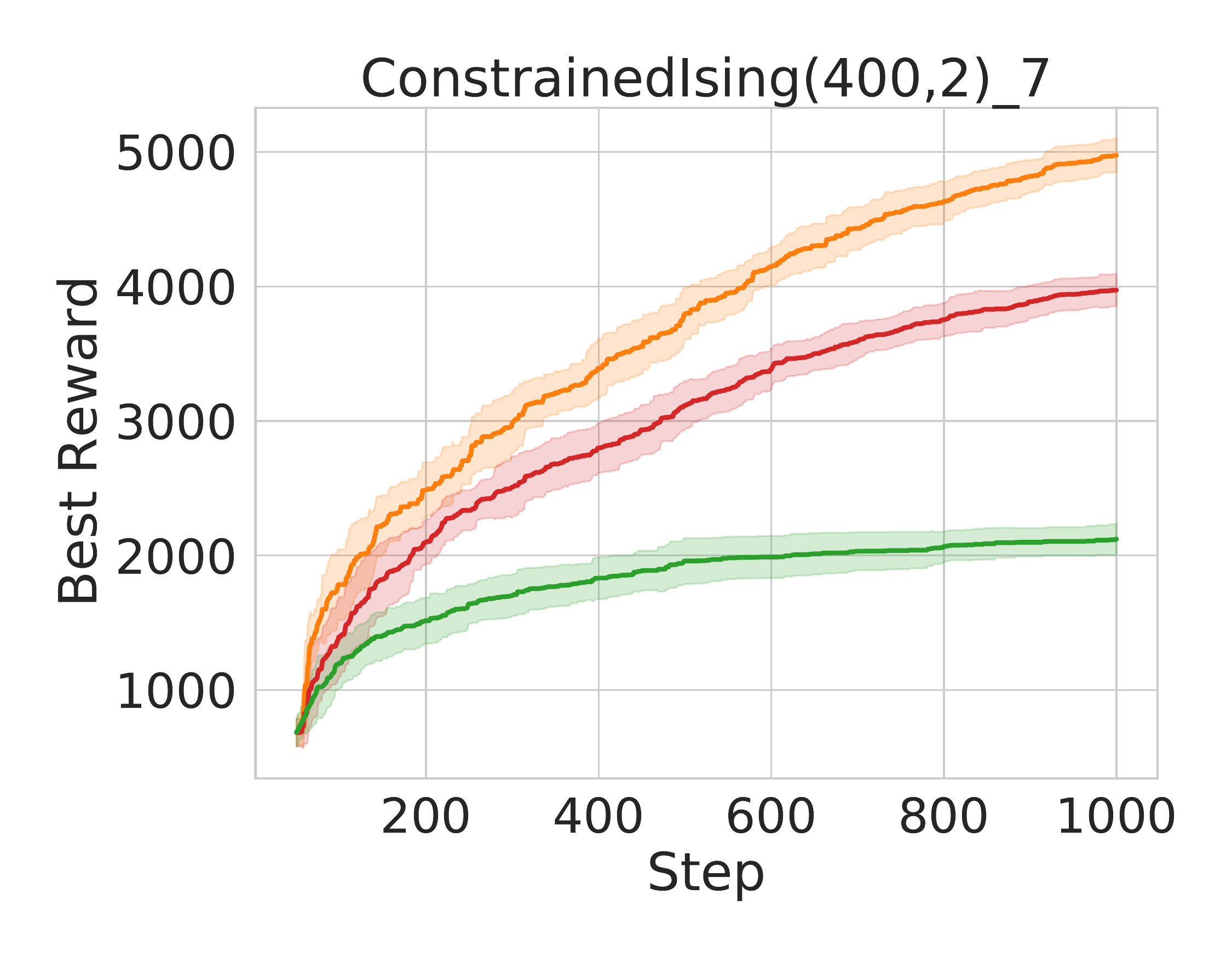}
    \includegraphics[width=\linewidth]{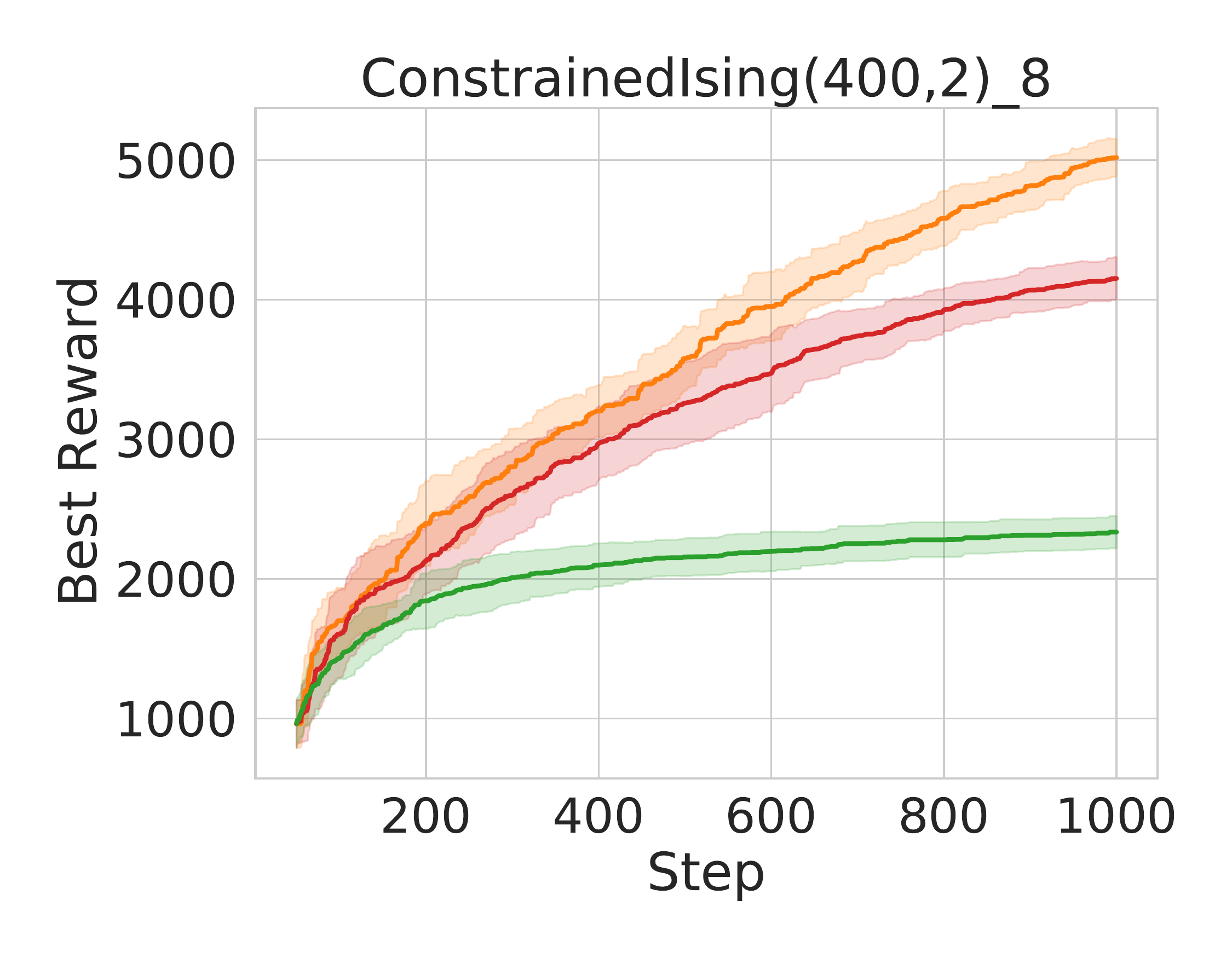}
    \includegraphics[width=\linewidth]{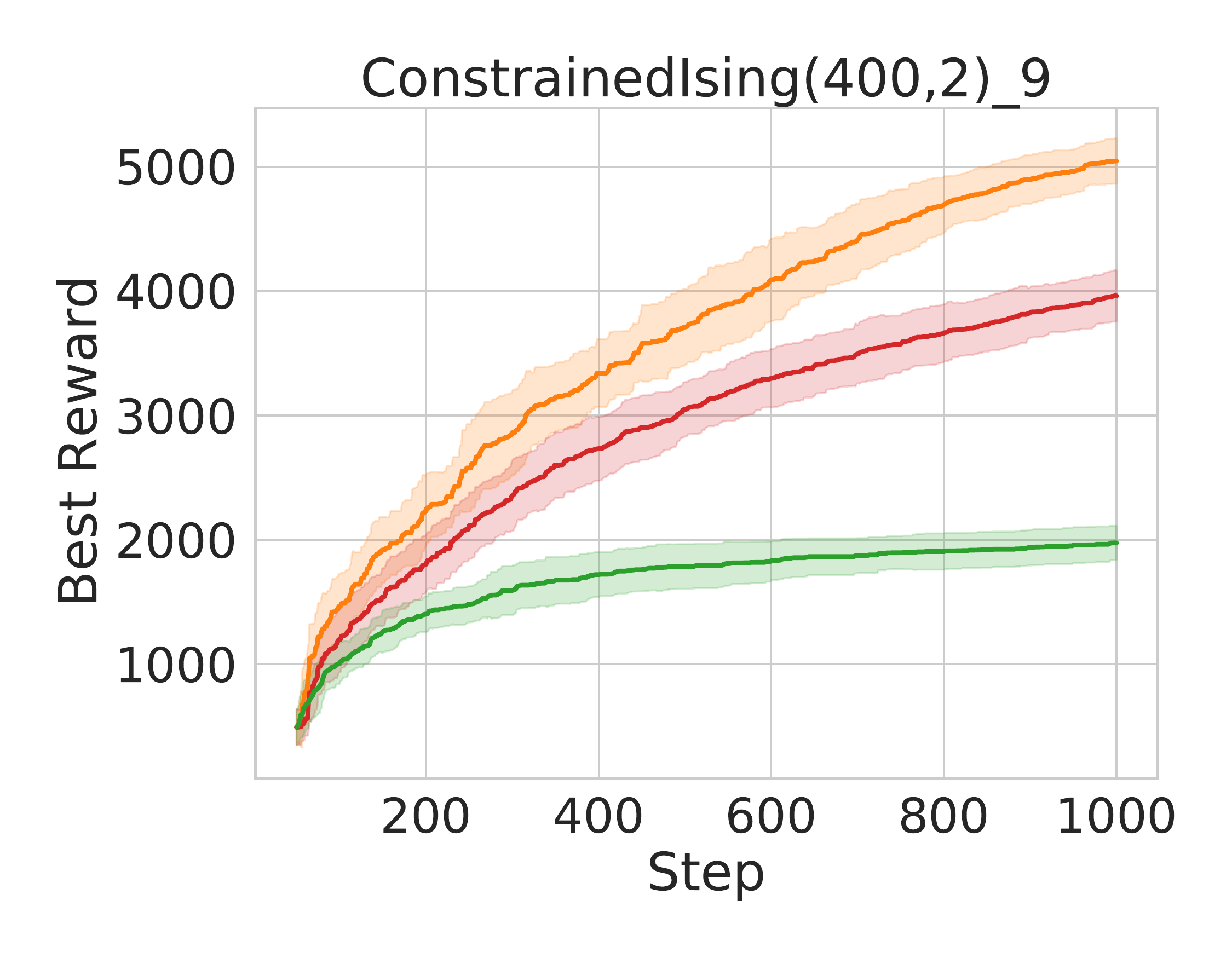}
    \end{multicols}
    \vspace{-20pt}
\centerline{\includegraphics[width=0.9\textwidth]{constr_ising_legend.pdf}}
\vspace{-15pt}
\caption{Best observed reward as a function of iteration for the second half of all constrained Ising model problems (Section~\ref{sec:exp_constr}), averaged over 10 trials (bands indicate $\pm 1$sd). Initial randomly sampled set of 50 points is omitted.}
\label{fig:appendix_constrained_curves2}
\end{center}
\end{figure*}

For the sake of completeness, we include in Figure~\ref{fig:appendix_constrained_scores} the distributions of algorithms' normalized max-reward scores over all constrained Ising problems (Section~\ref{sec:exp_constr}), paralleling Figure~\ref{fig:unconstrained_scores} for the unconstrained problems. As noted, while \textit{NN+MILP} and \textit{NN+ConEvo} perform similarly in the smaller instances, the former considerably improves over the latter as the problem size increases. The small variance in normalized scores reflects the fact that algorithms' relative reward progression was qualitatively similar in all problems within a given size, as can be seen in the individual reward curves in Figures~\ref{fig:appendix_constrained_curves1} and~\ref{fig:appendix_constrained_curves2}.

\subsubsection{Surrogate Model Capacity}
We next perform an experiment to evaluate the impact of the surrogate model's capacity on the quality of solutions as problem size increases. We include ablations of both \textit{NN+MILP} and \textit{NN+ConEvo} where the surrogate model has 32 neurons in the hidden layer, instead of the 16 used for the experiments in the main paper. Figure~\ref{fig:constrained_ising_model_larger_withfc32} plots each algorithms' best observed reward over time, averaged across all trials of all problems. We observe that neither NN+MILP nor NN+ConEvo show substantial improvements in performance when using a larger surrogate network in even the largest instances with 400 binary variables. This suggests that, in this case at least, the relatively small number of training points is a more significant bottleneck for objective approximation than the capacity of the surrogate.

\subsubsection{Random Optimizer}

One of our primary goals in this work is to contrast declarative vs. procedural approaches to handling constraints, as exemplified by the comparison of \textit{NN+MILP} and \textit{NN+ConEvo}. To this end, we include here a third baseline algorithm for the constrained experiments of Section~\ref{sec:exp_constr}. \textit{NN+RejSample} is an ablation of \textit{NN+MILP} where the inner-loop solver samples 10k feasible points uniformly-at-random from the domain and proposes the one with the highest acquisition function value. The configuration is otherwise identical to \textit{NN+MILP}. It is still an MBO algorithm, in the sense that it uses a surrogate to model the black-box objective, but uses naive random search for the inner-loop optimization.

Crucially, if we were to use rejection sampling in the inner loop, the solver could leverage the same exact declarative definition of constraints as in \textit{NN+MILP}. Unfortunately, the size of the feasible set for the subset-equality constraints is prohibitive for true rejection sampling: the chance of finding a feasible point is $<10^{-6}$ when $n=100, k=10$ and smaller for the other problem sizes. We therefore implement a custom sampling algorithm for this class of constraints that generates samples uniformly-at-random from the domain, and is thus equivalent to (though more computationally efficient than) rejection sampling. We note that, much like the custom mutator used by \textit{ConEvo}, this sampling procedure strongly relies on the special disjoint structure of the subset-equality constraints. In general, custom samplers might be much harder to design if constraints interact.

\begin{figure}[ht]
    \centering
    \includegraphics[width=\linewidth]{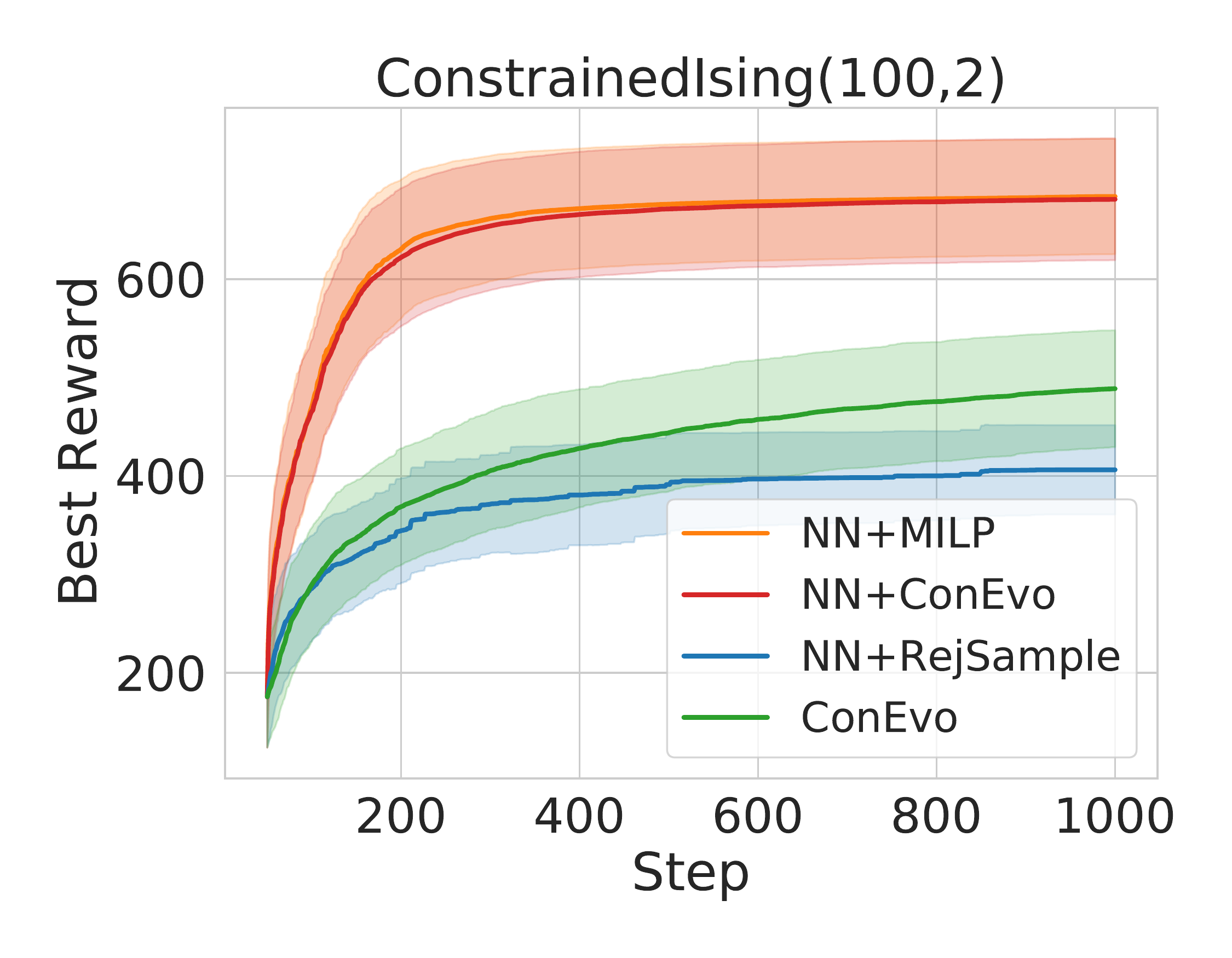}
    \caption{Best observed reward as a function of iteration for constrained Ising problems with $n=100$ variables, including the \textit{NN+RejSample} baseline. Note that the error bands are larger here as we aggregate over all trials of all problem instances, since they exhibit qualitatively similar reward trajectories.}
    \label{fig:appendix_constrained_with_nnrejsample}
\end{figure}

Figure~\ref{fig:appendix_constrained_with_nnrejsample} shows algorithms' best observed reward as a function of iteration averaged over all constrained Ising problems with $n=100$, now including \textit{NN+RejSample}. We do not run the algorithm on $n=200$ and $n=400$, as we expect random search to perform even worse in those larger domains. The poor performance of \textit{NN+RejSample}, even compared to \textit{ConEvo}, highlights the importance of using a high-quality optimizer in the inner-loop. This suggests that rejection sampling-based approaches, though easy to implement given a declarative definition of the constraints, are not likely to be effective when the domain is highly constrained.

\subsection{Practicality of MILP}
\label{sec:appendix_results_pract}

\subsubsection{Additional Timing Results}

We wish first to emphasize that all experiments in this paper were parallelized on a cluster of machines with variable hardware (all of them standard CPU machines with $\sim$1G RAM and $\leq 10$ cores however). As such, we intend our results in Section~\ref{sec:exp_pract} and here as an illustration of the practical computational overhead of our approach, and not as rigorous timing experiments. 

\begin{figure*}[ht]
\begin{center}
\includegraphics[width=\textwidth]{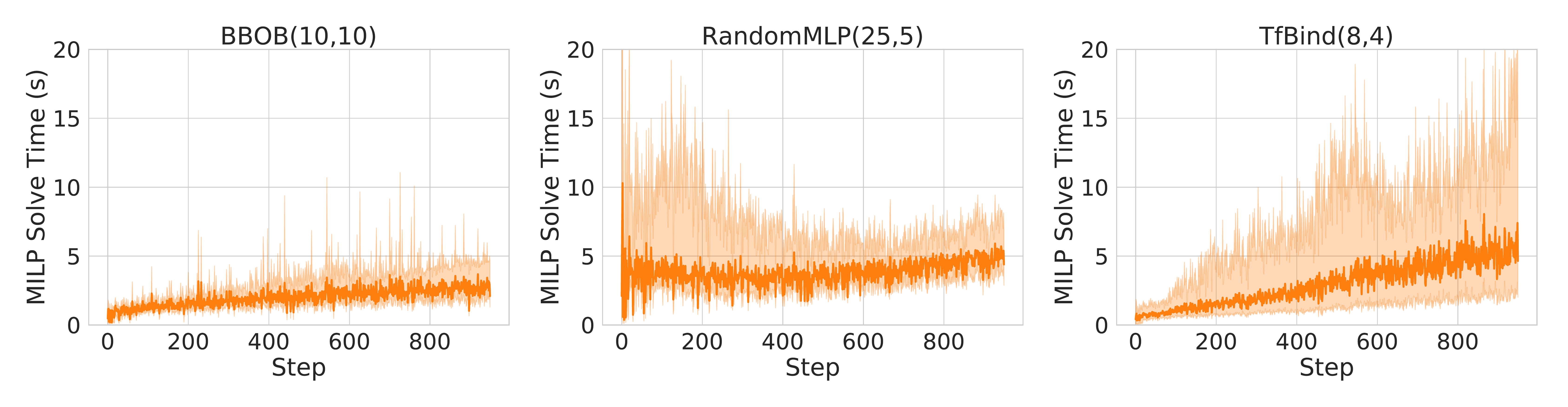}
\includegraphics[width=0.33\textwidth]{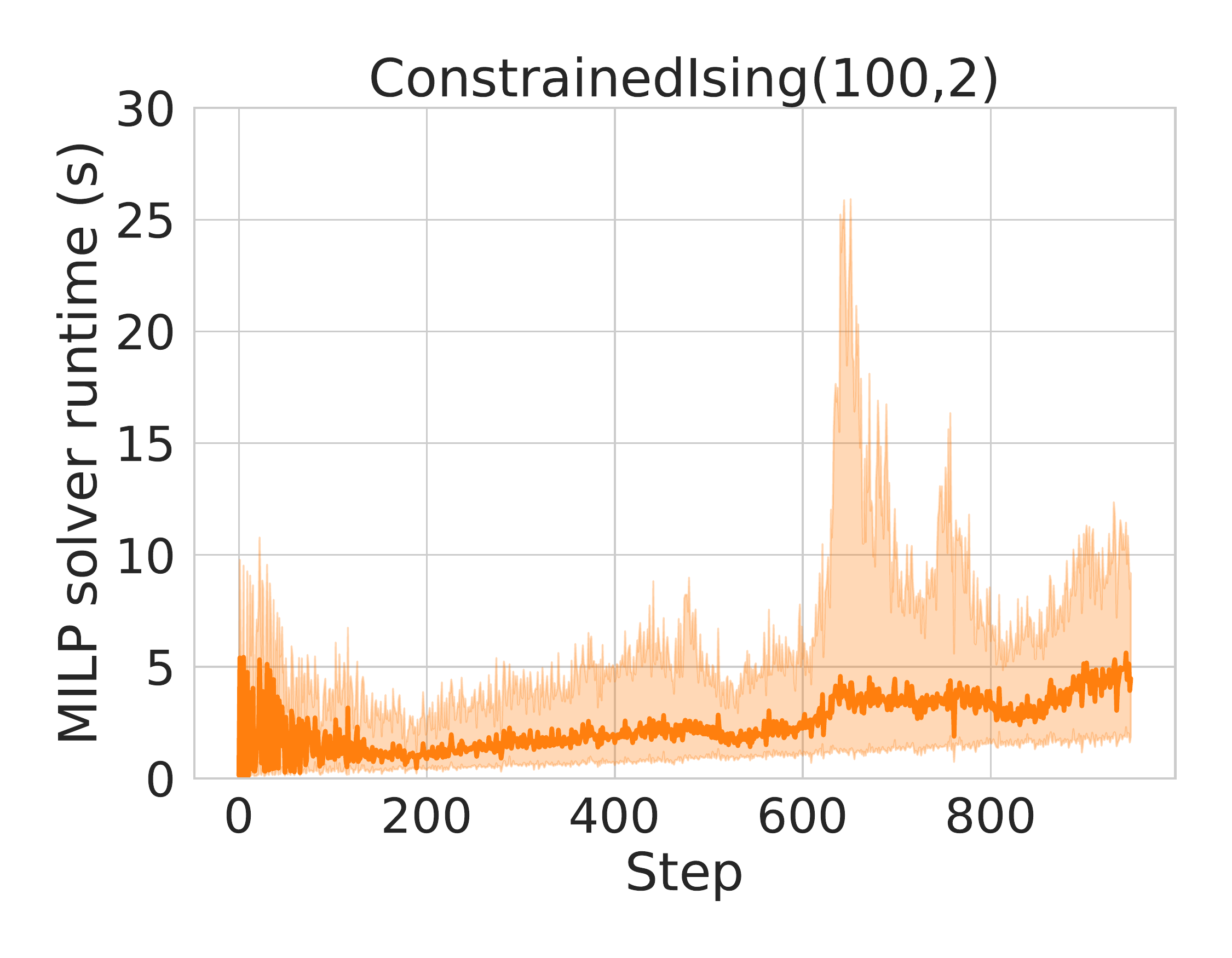}
\includegraphics[width=0.33\textwidth]{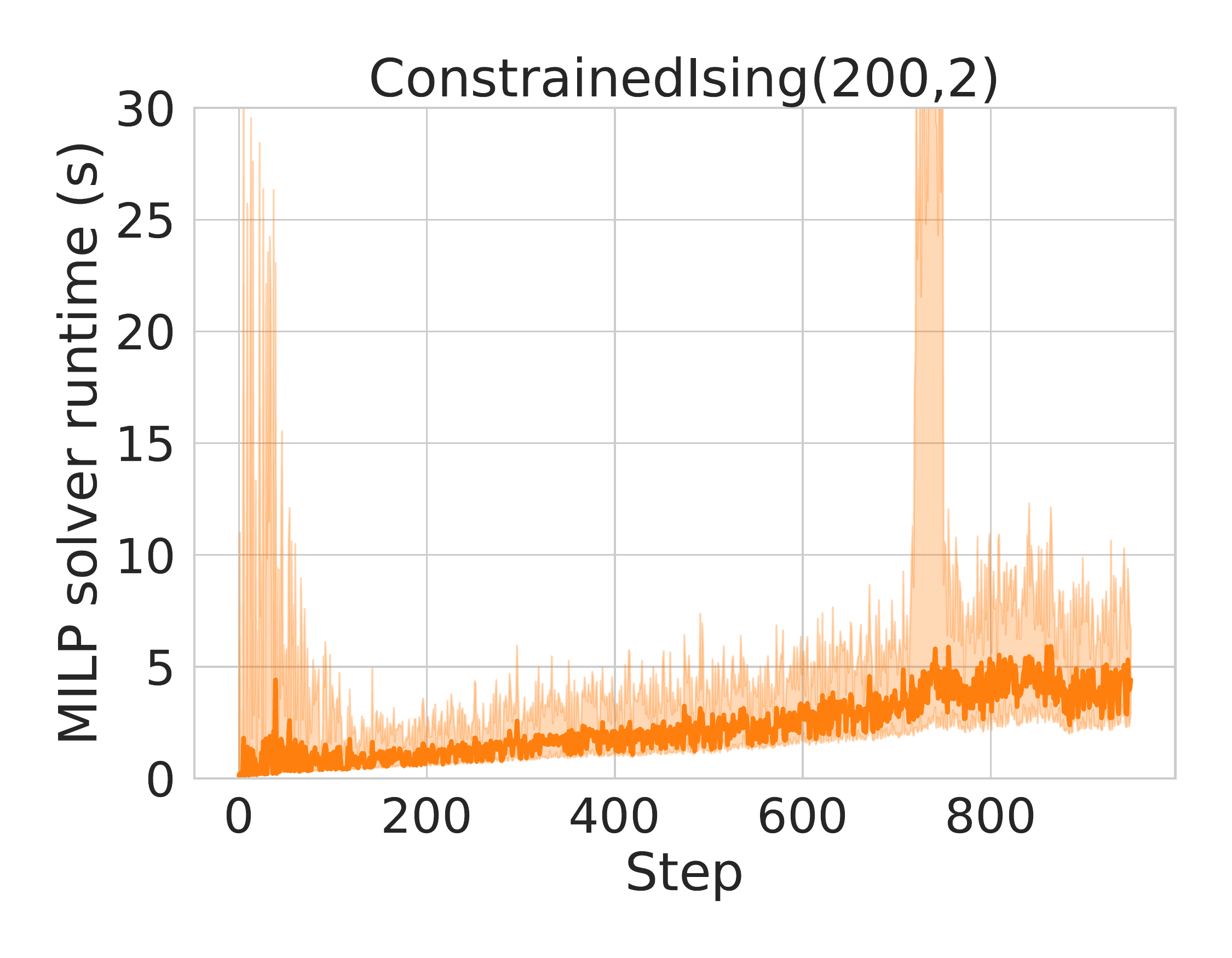}
\includegraphics[width=0.33\textwidth]{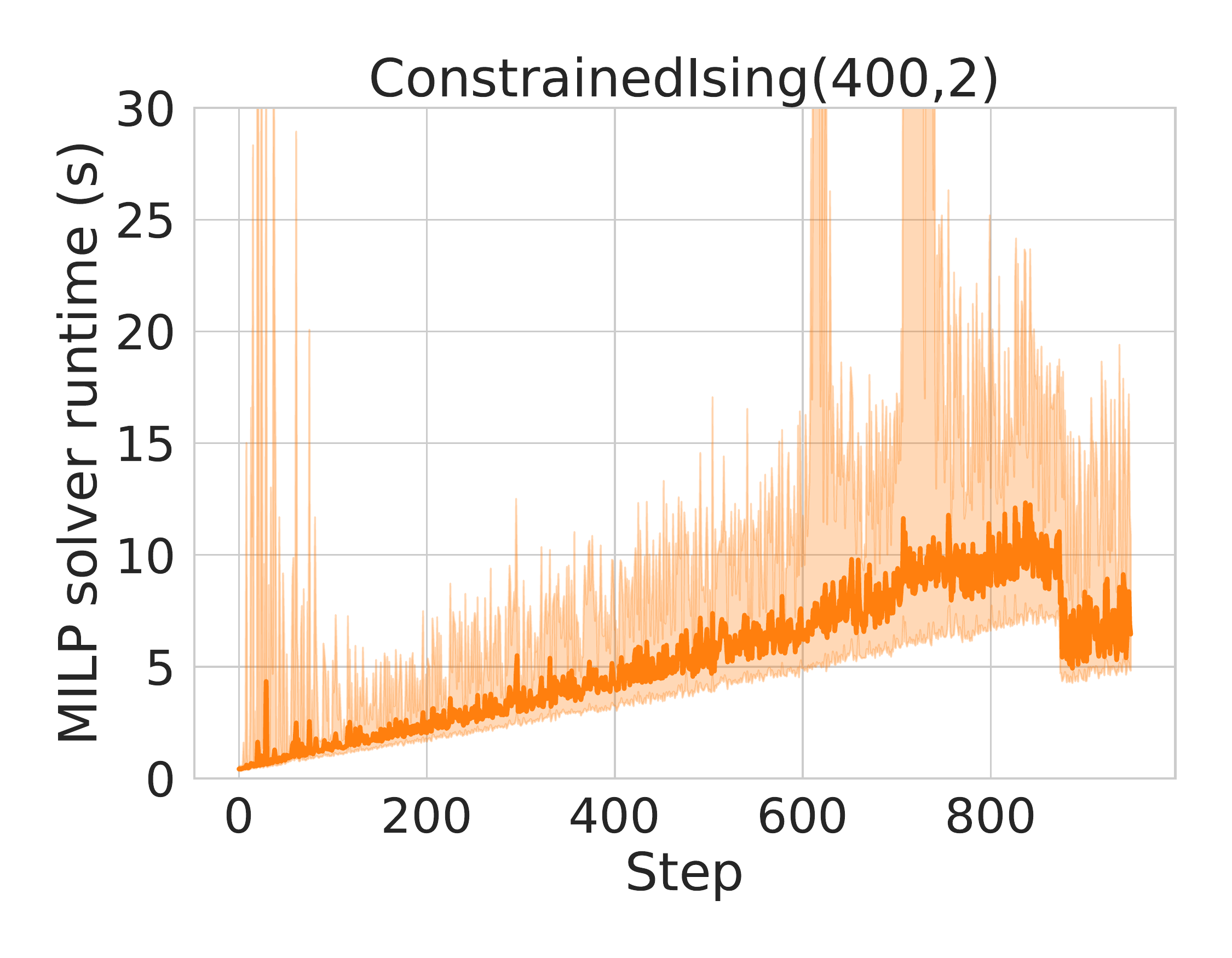}
\caption{Distribution of MILP acquisition problem solve times as a function of iteration split by objective class for unconstrained problems (Section~\ref{sec:exp_unconstr}) and constrained problems (Section~\ref{sec:exp_constr}). Line and bands show the median and 5th/95th percentile range over all trials of all problems in a class.}
\label{fig:appendix_runtimes}
\end{center}
\end{figure*}

In Figure~\ref{fig:appendix_runtimes} we plot the distribution of MILP acquisition problem solve times as a function of iteration for all experiments in Sections~\ref{sec:exp_unconstr} and~\ref{sec:exp_constr}, split by problem class (paralleling Figure~\ref{fig:unconstrained_runtimes} which included only the 12 unconstrained TfBind problems). As is often the case with MILP, the relationship between problem size and runtime can be unpredictable. For example, the lower-dimensional TfBind problems showed the highest mean and variance in solve time compared to the larger BBOB and RandomMLP problems in the unconstrained experiments. Moreover, even the largest constrained Ising problems ($n=400$) did not exhibit significantly higher average solve times than the unconstrained problems. 

We also observe a roughly linear increase in average solve as a function iteration, across all problem classes. This is presumably due to the increasing number of no-good constraints and the nature of surrogate models that have been fit on more data. 

\subsubsection{Scalability of MILP to Larger Networks}

We also perform an experiment to explore the impact of surrogate network size on MILP solve times. Here we vary \textit{NN+MILP}'s surrogate network architecture to use fully-connected (FC) networks with different numbers of hidden layers and neurons. We use parentheses to denote the number of neurons in each layer, e.g., FC(16) represents the single layer, 16-neuron network used throughout the main paper. We include ablations with FC(16), FC(32), FC(16,16), as well as a simple Linear model (no hidden layer), and run 20 trials of each, using different random initial datasets for each of the 12 unconstrained TfBind problems from Section~\ref{sec:exp_unconstr}. All other training and optimization hyper-parameters for \textit{NN+MILP} are the same as described in Section~\ref{sec:appendix_solvers_nnmilp}. 

Table~\ref{tab:architecture_runtimes} shows aggregate distribution statistics of acquisition solve time for the different architectures, across all steps of all trials of all problems. We note that solve times increase as the network size increases, but even for the largest network (two layers with 16 neurons each), the solver rarely times out and almost always terminates within a practical time limit. Furthermore, for larger networks we can improve scaling using advanced formulation techniques (e.g. Appendix~\ref{sec:relu_formulation}) or commercial MILP solvers (e.g., Gurobi). We also note that, in these experiments, there was no single architecture that consistently produced better optimization across different instances (though \textit{Linear} was almost always outperformed by the rest). 

\begin{table}[h]
\centering
\caption{Distribution of per-step MILP inner-optimization solve times in seconds for TfBind8 benchmarks when using different surrogate network architectures. The \textit{Network} column denotes the number of ReLUs in each fully-connected hidden layer. Runs were given a time limit (TL) of 300s. (*) means that the time limit was hit.}
\label{tab:architecture_runtimes}
\begin{tabular}{crrrrrrr}
\hline
Network & min & med & 95\% & 99\% & max & \%TL\\ \hline 
Linear & 0.004 & 0.4 & 1.4 & 2.9 & 16.9 & 0\%\\
FC(16) & 0.02 & 2.2 & 8.0 & 15.5 & 60.8 & 0\%\\
FC(32) & 0.04 & 11.7 & 49.2 & 85.5 & 300* & 0.1\%\\
FC(16,16) & 0.40 & 12.2 & 55.6 & 109.1 & 300* & 2.1\%\\
\hline
\end{tabular}
\end{table}

\section{Constrained binary quadratic problems from MINLPLib}
\label{sec:app_minlplib}

\begin{figure*}[ht]
\begin{center}
\includegraphics[width=\textwidth]{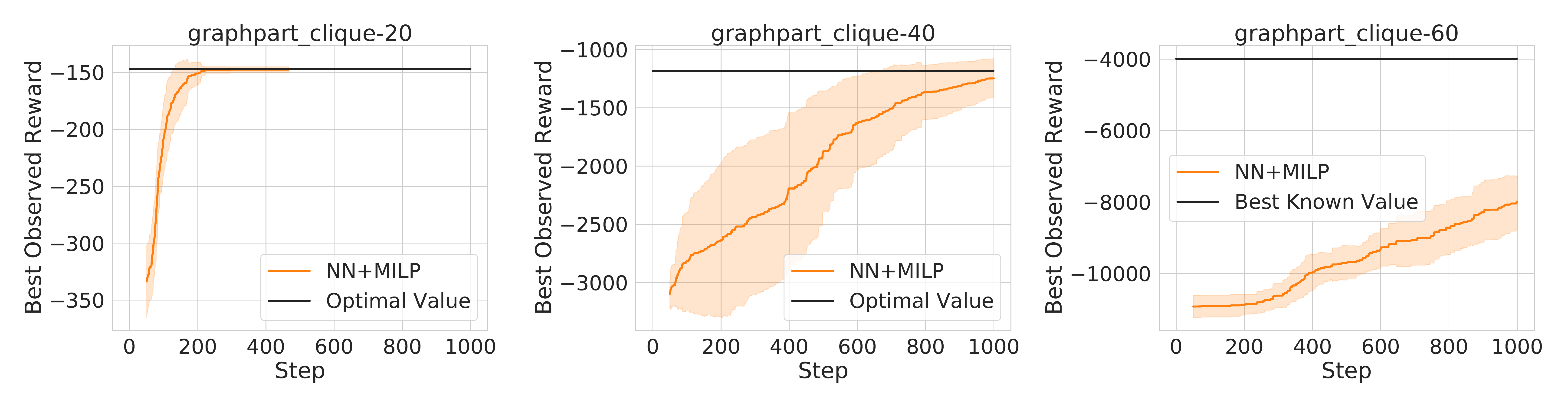}
\caption{Best observed reward as function of iteration for three graph partitioning instances from MINLPLib (negated for maximization), with 60, 120, and 180 binary variables respectively. Black lines show the best known feasible solution to the problem (as of October 1, 2021). Colored lines show the average over 20 trials, while bands indicate $\pm 1$ sd. Note that bands that exceed the black line are an artifact of the symmetric nature of standard deviation, and do not necessarily mean a trial found an improved solution.}
\label{fig:qp_minlp_example}
\end{center}
\end{figure*}

In this section, we show results for individual MINLPLib instances, in which we examine primal gaps (see Section~\ref{sec:minlp}) of \textit{NN+MILP} with respect to the best known primal feasible solution from the MINLPLib benchmark itself (as of October 1, 2021). We select all the instances of type ``BQP'' from MINLPLib with at least one linear constraint. We note that prefixes correspond to different classes of problems based on the constraints, which formed the basis of our categorization in Table~\ref{tab:minlp_summ}; the full set includes 31 graph partitioning (\texttt{graphpart}), 8 generalized assignment (\texttt{pb}), 6 quadratic shortest-path (\texttt{qspp}), and 16 additional (\texttt{other}) problems. For consistency with the remainder of the paper, we turn the problems into maximization problems by negating the objective function.

We run NN+MILP with the same settings as previous experiments (see Appendix~\ref{sec:appendix_solvers_nnmilp}). One difference here is that the feasible set may be too small to sample from using rejection sampling. Therefore, we build our initial set of 50 points by randomly choosing an objective direction and solving an MILP under the constraints of the problem, which is practically feasible since the scale of these problems is small in the context of MILP.

Table~\ref{tab:qp_minlp_results} lists the results of the 20 trials for each instance, omitting 11 instances prefixed by \texttt{celar6-sub0}, \texttt{color\_lab}, and \texttt{maxcsp} since our method was unable to find a solution with primal gap at most 10\% (although they are taken into account in the statistics in Section~\ref{sec:minlp}). Interestingly, in 20 out of the 61 instances, we match the best known objective value in at least one of the runs. This includes a number of instances that are considered to be large for black-box optimization, such as the general quadratic assignment problem \texttt{pb302095} which has 600 variables.

On the other hand, we also observe that the method has difficulties in finding a good solution for larger instances. This is more clearly illustrated by Figure~\ref{fig:qp_minlp_example}, in which we observe how the method scales with the graph partitioning problems denoted by \texttt{graphpart\_clique}. For the smaller instance (with 60 variables), our method finds an optimal solution in relatively few steps, but it has difficulties in reaching the best known solution for the larger instance (with 180 variables) with the same constraint structure, at least with the current parameters and network size.

\clearpage
\section{NAS-Bench-101 Case Study}
\label{sec:nasbench_formulation}

\subsection{Background}
\label{sec:app_nasbench_background}

In Section~\ref{sec:nas} we consider the NAS-Bench-101 \citep{nasbench101} neural architecture search (NAS) benchmark as a case study, to illustrate the power of MILP's declarative constraint language in modeling complex combinatorial domains. The optimization domain consists of directed acyclic graphs (DAGs) with a maximum of $V = 7$ nodes and $M = 9$ edges, representing the \textit{cell} in a neural architecture. The overall model is obtained by stacking multiple copies of the cell. Two nodes represent the input and output, and must be connected by a directed path, while the remaining nodes are each assigned to be 1x1 convolution, 3x3 convolution, or 3x3 max-pooling. Edges specify the flow of activations between nodes. The goal is to find the cell architecture that maximizes out-of-sample accuracy on a given image classification task.

NAS differs from the problem setting in Section \ref{sec:problem_setting} in three key ways. First, algorithms do not have access to the true objective (out-of-sample accuracy), but instead a correlated proxy (validation accuracy). Second, $f(x)$ is noisy due to the stochasticity of classifier training, and thus algorithms may benefit from repeated queries of the same point. Finally, algorithms may benefit by leveraging the validation accuracy at early epochs as a proxy to halt unpromising evaluations. 

Despite these differences, we apply \textit{NN+MILP} exactly as described in Section~\ref{sec:algo}/Appendix~\ref{sec:appendix_solvers_nnmilp} (including no-good constraints which prevent repeated queries). We reimplement the Regularized Evolution (\textit{RE}) and random search (\textit{RS}) baselines from the original NAS-Bench-101 paper \citep{nasbench101}, using the same hyper-parameters settings specified therein.  

The NAS-Bench-101 dataset contains pre-computed validation \textit{and} test accuracies for three independently trained replications of each architecture, as well as the training time of each. To simulate NAS, algorithms' observed reward after proposing an architecture is the validation accuracy of a randomly sampled replication from said architecture. This defines the notion of an ``incumbent'' proposal, namely the proposed architecture with the highest (observed) validation accuracy, which may not in fact be the best (unobserved) test accuracy. Instead of allowing algorithms a fixed budget of evaluations, we use a fixed budget of $T=5 \times 10^6$ seconds, and allow algorithms to query the objective until cumulative training time exceeds the budget. For evaluation purposes (e.g., Figure~\ref{fig:nasbench_results}) we plot the out-of-sample accuracy of the incumbent architecture as a function of cumulative architecture training time. 

\subsection{Domain Formulation}
\label{sec:app_nasbench_form}

To formulate the NAS-Bench-101 domain, we first define a representation of cell architectures as fixed-length binary vectors. We split the representation into two components; one set of variables encodes the presence or absence of each graph edge, while the second is a one-hot encoding of nodes' assigned operations. As all valid cell graphs are directed and acyclic, we limit the edge variables to the strict upper triangle of the adjacency matrix, which implicitly enforces a topological ordering of the nodes in any feasible solution and ensures acyclicity. The first- and last-indexed nodes are always assigned the \texttt{input} and \texttt{output} operations respectively, while intermediates nodes can be assigned any operation from the set $\mathcal{S} = \{\texttt{conv1x1}, \texttt{conv3x3}, \texttt{maxpool3x3}\}$.

To ensure a fixed-length set of decision variables while allowing for graphs with a variable number of nodes, we introduce a new \texttt{null} operation. Nodes assigned the null operation are not considered part of the computational graph of the cell. The algorithm then searches over the space of binary representations, constrained to yield feasible cell architectures. 

\newpage
Denoting by $V$ and $M$ the maximum number of allowable nodes and edges respectively, the decision variables (all binary) are:

\begin{itemize}
    \setlength{\itemsep}{0pt}
    \item $m_{i,j}$ for $1 \leq i < j \leq V$,  1 if there is an edge from node $i$ to node $j$, 0 otherwise.
    \item $w_{i,k}$ for $1 \leq i \leq V, 1 \leq k \leq |\mathcal{S}|$, 1 if node $i$ is assigned the $k$'th operation in $\mathcal{S}$, 0 otherwise.
    \item $z_i$ for $1 < i < V$, 1 if node $i$ is assigned the null operation,0 otherwise.
\end{itemize}

The feasible set of cell architectures can then be given in terms of linear constraints as follows:

\begin{flalign}
\setcounter{equation}{0}
    w_{1,k} = w_{V,k} = z_1 = z_V = 0  & \qquad \text{for } 1 \leq k \leq |\mathcal{S}| \label{eq:inout} \\
    z_i + \sum_{k=1}^{S} w_{i,k} = 1 & \qquad \text{for } 1 < i < V \label{eq:onehot} \\
    \sum_{i=1}^V \sum_{j=i+1}^V m_{i,j} \leq M & \label{eq:maxedges} \\
    m_{i,j} \leq 1 - z_j & \qquad \text{for } 1 \leq i < j \leq V \label{eq:nullin} \\
    m_{i,j} \leq 1 - z_i & \qquad \text{for } 1 \leq i < j \leq V \label{eq:nullout} \\
    \sum_{i=1}^{j-1} m_{i,j} \geq 1 - z_j & \qquad \text{for } 1 \leq j \leq V \label{eq:flowin} \\
    \sum_{j=i+1}^{V} m_{i,j} \geq 1 - z_i & \qquad \text{for } 1 \leq i \leq V \label{eq:flowout}\\
    z_i \leq z_{i+1} & \qquad \text{for } 1 < i < V - 1 \label{eq:symmetry}
\end{flalign}

Constraints~\ref{eq:inout} ensure that the input and output nodes are not assigned any operation from $\mathcal{S}$ or null, while \ref{eq:onehot} enforces the one-hot encoding of operations for intermediate nodes (including the possibility of a null operation). Constraint~\ref{eq:maxedges} imposes a limit on the number of edges in the graph, per the NAS-Bench-101 specifications. Constraints~\ref{eq:nullin}~\&~\ref{eq:nullout} assert that null nodes have no incoming or outgoing edges respectively, effectively disconnecting them from the remaining graph. Conversely, \ref{eq:flowin}~\&~\ref{eq:flowout} assert that non-null nodes have at least one ingoing and one outgoing edge. Crucially, due to the implicit topological sorting of nodes by the upper-triangular adjacency matrix, these also ensure that there is always a path from the input to the output node using only non-null nodes. Intuitively, all non-null nodes (including the input) have at least one outgoing edge -- which necessarily leads to a higher-indexed non-null node -- and all non-null nodes (including the output) have at least one-incoming edge -- which necessarily comes from a lower-indexed non-null node. The flow exiting the input node, must eventually enter the output node.

Finally, we focus on Constraints~\ref{eq:symmetry}, which we refer to as \textit{symmetry-breaking} constraints. These assert that a node can only be assigned the null operation if its topological successor has also been assigned it. While not necessary for feasibility, this constraint serves to eliminate symmetry by ensuring that all null nodes are topologically sorted after any non-null nodes. In essence, it introduces a ``canonical'' labeling of null vs. non-null nodes, whose isomorphic representations are excluded from the feasible region. We refer to the optimization algorithm that includes the symmetry-breaking constraints as \textit{NN+MILP (w/ symmetry)}, and the one that excludes them as simple \textit{NN+MILP}.

\subsection{Additional results}
\label{sec:app_nasbench_results}

\begin{figure}[ht]
    \centering
    \includegraphics[width=0.9\linewidth]{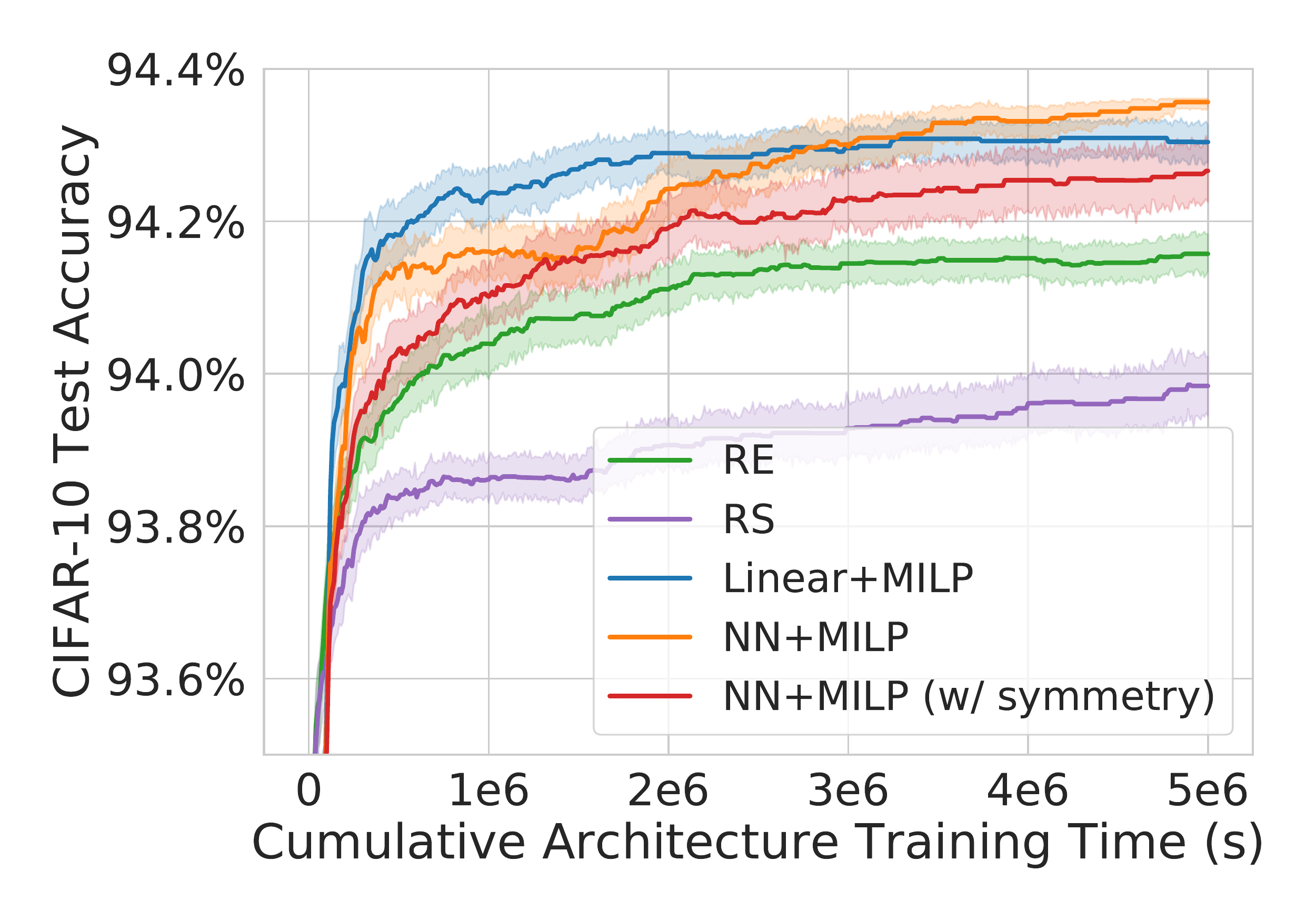}
    \caption{Test accuracy of algorithms' incumbent architecture as a function of cumulative training time on NAS-Bench-101, averaged over 100 trials, including the formulation with symmetry-breaking constraints. Bands indicate 95\% confidence interval for the mean.}
    \label{fig:nasbench_results_with_symmetry_breaking}
\end{figure}

In our experiments, we found that including symmetry-breaking constraints actually resulted in worse overall performance for the outer optimization problem (Figure~\ref{fig:nasbench_results_with_symmetry_breaking}). We hypothesize that this is due to a reduction in the exploration behaviour of \textit{NN+MILP}, as the surrogate's predictive distribution was more uncertain in the larger search space and the inner-loop optimizer thus more likely to propose points in unexplored areas. One possible future line of work could be to augment $\mathcal{D}_t$ with isomorphic representations before training, e.g., by random reordering of nodes in the representations of sampled points.

\clearpage
\begin{table*}[ht]
\centering
\caption{Results for a subset of binary quadratic problems from the MINLPLib benchmark, indicating the number of runs out of 20 for which the primal gap with respect to the best known primal feasible solution is at most 0\%, 1\%, and 10\%.}
\label{tab:qp_minlp_results}
\begin{tabular}{lrrrr}
\hline
&&\multicolumn{3}{c}{Number of runs with solution at}\\
Instance name & \# variables & 0\% gap & $\leq$ 1\% gap & $\leq$ 10\% gap\\ \hline
cardqp\_inlp & 50 & 5 & 14 & 20\\
cardqp\_iqp & 50 & 5 & 14 & 20\\
crossdock\_15x7 & 210 & 0 & 0 & 20\\
crossdock\_15x8 & 240 & 0 & 0 & 20\\
graphpart\_2g-0044-1601 & 48 & 17 & 17 & 20\\
graphpart\_2g-0055-0062 & 75 & 0 & 2 & 19\\
graphpart\_2g-0066-0066 & 108 & 0 & 0 & 17\\
graphpart\_2g-0077-0077 & 147 & 0 & 0 & 7\\
graphpart\_2g-0088-0088 & 192 & 0 & 0 & 7\\
graphpart\_2g-0099-9211 & 243 & 0 & 0 & 2\\
graphpart\_2g-1010-0824 & 300 & 0 & 0 & 0\\
graphpart\_2pm-0044-0044 & 48 & 20 & 20 & 20\\
graphpart\_2pm-0055-0055 & 75 & 13 & 13 & 20\\
graphpart\_2pm-0066-0066 & 108 & 6 & 6 & 16\\
graphpart\_2pm-0077-0777 & 147 & 0 & 0 & 7\\
graphpart\_2pm-0088-0888 & 192 & 0 & 0 & 6\\
graphpart\_2pm-0099-0999 & 243 & 0 & 0 & 3\\
graphpart\_3g-0234-0234 & 72 & 0 & 4 & 19\\
graphpart\_3g-0244-0244 & 96 & 0 & 1 & 17\\
graphpart\_3g-0333-0333 & 81 & 4 & 4 & 19\\
graphpart\_3g-0334-0334 & 108 & 0 & 0 & 16\\
graphpart\_3g-0344-0344 & 144 & 0 & 1 & 10\\
graphpart\_3g-0444-0444 & 192 & 0 & 0 & 8\\
graphpart\_3pm-0234-0234 & 72 & 7 & 7 & 19\\
graphpart\_3pm-0244-0244 & 96 & 0 & 0 & 17\\
graphpart\_3pm-0333-0333 & 81 & 1 & 1 & 15\\
graphpart\_3pm-0334-0334 & 108 & 0 & 0 & 12\\
graphpart\_3pm-0344-0344 & 144 & 0 & 0 & 5\\
graphpart\_3pm-0444-0444 & 192 & 0 & 0 & 6\\
graphpart\_clique-20 & 60 & 20 & 20 & 20\\
graphpart\_clique-30 & 90 & 20 & 20 & 20\\
graphpart\_clique-40 & 120 & 13 & 13 & 18\\
graphpart\_clique-50 & 150 & 0 & 0 & 0\\
graphpart\_clique-60 & 180 & 0 & 0 & 0\\
graphpart\_clique-70 & 210 & 0 & 0 & 0\\
pb302035 & 600 & 0 & 0 & 0\\
pb302055 & 600 & 0 & 0 & 20\\
pb302075 & 600 & 6 & 6 & 20\\
pb302095 & 600 & 16 & 20 & 20\\
pb351535 & 525 & 0 & 0 & 18\\
pb351555 & 525 & 1 & 4 & 20\\
pb351575 & 525 & 0 & 10 & 20\\
pb351595 & 525 & 7 & 20 & 20\\
qap & 225 & 0 & 0 & 7\\
qspp\_0\_10\_0\_1\_10\_1 & 180 & 5 & 5 & 20\\
qspp\_0\_11\_0\_1\_10\_1 & 220 & 9 & 20 & 20\\
qspp\_0\_12\_0\_1\_10\_1 & 264 & 13 & 13 & 20\\
qspp\_0\_13\_0\_1\_10\_1 & 312 & 0 & 19 & 20\\
qspp\_0\_14\_0\_1\_10\_1 & 364 & 0 & 16 & 20\\
qspp\_0\_15\_0\_1\_10\_1 & 420 & 12 & 13 & 20\\\hline
\end{tabular}
\end{table*}

\clearpage

\clearpage

\clearpage

\clearpage

\end{document}